\documentclass[final,3p,times]{elsarticle}

\usepackage{amsmath, amssymb}

\usepackage{color}

\usepackage{hyperref}       
\usepackage{url}            
\usepackage{amsfonts}       
\usepackage{nicefrac}       
\usepackage{microtype}      
\usepackage{doi}
\usepackage{dcolumn}
\usepackage{bm}
\usepackage{caption}
\usepackage{subcaption}
\usepackage{float}
\usepackage{array}
\newcolumntype{P}[1]{>{\centering\arraybackslash}p{#1}}
\usepackage{pdflscape}
\usepackage{amssymb}

\journal{Journal of Computational Physics}

\begin{document}

\begin{frontmatter}



\title{Multifidelity Deep Operator Networks For Data-Driven and Physics-Informed Problems}

\author[PNNL]{Amanda A. Howard}
\author[Sandia]{Mauro Perego}
\author[PNNL,Brown]{George Em Karniadakis}
\author[PNNL,cor]{Panos Stinis}

\cortext[cor1]{panagiotis.stinis@pnnl.gov}
\affiliation[PNNL]{organization={Pacific Northwest National Laboratory},
            addressline={P.O. Box 999}, 
            city={Richland},
            postcode={99352}, 
            state={WA},
            country={USA}}

\affiliation[Sandia]{organization={Sandia National Laboratories},
            addressline={P.O. Box 5800}, 
            city={Albuquerque},
            postcode={87185}, 
            state={NM},
            country={USA}}

\affiliation[Brown]{organization={Division of Applied Mathematics and School of Engineering, Brown University},
addressline={182 George Street}, 
city={Providence},
postcode={02912}, 
state={RI},
country={USA}}
            
\begin{abstract}
Operator learning for complex nonlinear systems is increasingly common in modeling multi-physics and multi-scale systems. However, training such high-dimensional operators requires a large amount of expensive, high-fidelity data, either from experiments or simulations. In this work, we present a composite Deep Operator Network (DeepONet) for learning using two datasets with different levels of fidelity to accurately learn complex operators when sufficient high-fidelity data is not available. Additionally, we demonstrate that the presence of low-fidelity data can improve the predictions of physics-informed learning with DeepONets. We demonstrate the new multi-fidelity training in diverse examples, including modeling of the ice-sheet dynamics of the Humboldt glacier, Greenland, using two different fidelity models and also using the same physical model at two different resolutions.
\end{abstract}




\end{frontmatter}


\section{Introduction}

In many applications across science and engineering it is common to have access to disparate types of data with different levels of fidelity. In general, low-fidelity data is easier to obtain in greater quantities, but it may be too inaccurate or not dense enough to accurately train a machine learning model. High-fidelity data is costly to obtain, so there may not be sufficient data to use in training, however, it is more accurate.  A small amount of high fidelity data, such as from measurements, combined with low fidelity data, can improve predictions when used together; this has motivated geophysicists to develop  {\em cokriging} \cite{wackernagel1995cokriging}, which is based on Gaussian process regression at two different fidelity levels by exploiting correlations- albeit only linear ones - between different levels. An example of cokriging for obtaining the sea surface temperature (as well as the associated uncertainty) is presented in  \cite{babaee2020multifidelity}, where satellite images are used as low-fidelity data whereas {\em in situ} measurements are used as high-fidelity data. To exploit nonlinear correlations at different levels of fidelity, a probabilistic framework based on Gaussian process regression and nonlinear autoregressive scheme was proposed in \cite{perdikaris2017nonlinear} that can learn complex nonlinear and space-dependent cross-correlations between multifidelity models. However, the limitation of this work is the high computational cost for big data sets, and to this end, the subsequent work in \cite{meng2020composite} was based on neural networks and provided the first method of multifidelity training  of deep neural networks.

Let us consider the motivating example of a low-order numerical model that is less computationally expensive, so it can be used to generate large amounts of low-fidelity data. However, low-order numerical models are less accurate and a trained deep learning model will only be as accurate as the model. In contrast, a higher-order numerical model can be too costly to run, or experimental data can be too expensive to generate in sufficient quantities for training. This situation is exacerbated when the machine learning model is required to approximate not just a function as in the aforementioned works, but an operator, \emph{e.g.}, the evolution operator of a partial differential equation (PDE). In the current work, we present a method that can combine low-fidelity and high-fidelity data to produce more accurate {\em operator approximations} than using either dataset alone. In addition, the proposed approach allows the incorporation of  information from physics in the training stage in order to further improve the approximation accuracy and possibly reduce the data generation cost. 

\textcolor{black}{Several methods have been developed for operator approximations, including Deep Operator Networks \cite{lu2021learning}, Fourier Neural Operators, \cite{li2020fourier}, Graph Kernel Networks \cite{li2020neural}, and Nonlocal Kernel Networks \cite{you2022nonlocal}.} The Deep Operator Network (DeepONet) framework \cite{lu2021learning}, inspired by the universal approximation theorem for operators \cite{chen1995universal, back2002universal}, allows for the learning of operators between infinite-dimensional spaces. DeepONets have been accurately applied to a wide range of applications \cite{sharma2021application, ranade2021generalized, goswami2021physics}, including bubble dynamics across a range of length scales \cite{lin2021seamless}, prediction of failure due to cracks \cite{goswami2021physics}, and prediction of linear instability waves in high-speed boundary layers \cite{di2021deeponet}. More recently, \cite{wang2021improved} introduced a modified DeepONet architecture, which is shown to lead to accurate results across a range of problems. \textcolor{black}{While in this work we will focus on DeepONets, the idea of multifidelity learning for operators is universal across all operator learning methods, and ideas from this paper can be applied to other methods.}

Physics-informed neural networks (PINNs) usually train for one set of input parameters \cite{raissi2019physics, sirignano2018dgm, karumuri2020simulator, sun2020surrogate, psichogios1992hybrid, lagaris1998artificial}. By construction, DeepONets are suited for training for a whole range of input parameters. In addition, because DeepONet outputs are differentiable with respect to the input coordinates, the same framework used in PINNs \cite{raissi2019physics}, which relies on automatic differentiation \cite{griewank1989automatic, baydin2018automatic}, can also be applied to DeepONets. Indeed, the work in \cite{wang2021learning} and \cite{wang2021long} extends the DeepONet framework by allowing the inclusion of physics-informed terms in the loss function. This can potentially eliminate the need for training data in the form of input-output pairs required to find \emph{e.g.,} the solution to known or parametrized PDEs. 

Previous work with PINNs has enabled multifidelity learning of functions with both data and physics-informed training \cite{meng2020composite, lu2020extraction, penwarden2022multifidelity, jagtap2022deep, regazzoni2021physics}. This work can broadly be split into three categories: transfer learning, where the network is first trained for the low-fidelity data, then a correction is found to correct the low-fidelity output;  simultaneous training; and consecutive training. Transfer learning can be trained with or without physics \cite{song2021transfer, de2020transfer, de2022neural}. The approaches in \cite{chakraborty2021transfer} and \cite{harada2022application} use transfer learning and a two-step training process to first enforce (approximate) physics and then use a small amount of high-fidelity data to correct the trained PINN output. This process applies when the physics is either not known exactly or the cost to generate low-fidelity data with a solver is too high. In simultaneous training, \cite{meng2020composite} learns both the linear and nonlinear correlations between the low-fidelity and high-fidelity data, by simultaneously training the low- and high-fidelity networks. This method allows for problems with complex correlations between the two datasets \cite{meng2021fast}. Simultaneous training can also be applied successfully without physics \cite{meng2020composite, zhang2021multi}. In the consecutive training category, \cite{liu2019multi} trains three networks, two low- and high-fidelity physics-constrained neural networks and a third neural network to learn the correlation between the low- and high-fidelity output. \textcolor{black}{We note that very recent work considers bifidelity data with DeepONets \cite{lu2022multifidelity, de2022bi}. 
\cite{lu2022multifidelity} considers several different architectures for multifidelity DeepONets, which allows the learned low-fidelity DeepONet to inform training of the high-fidelity DeepONet through input augmentation, where the low-fidelity output is used as an additional input for the high-fidelity DeepONet, and residual learning. The authors do not explicitly learn the linear and nonlinear correlations between the low- and high-fidelity datasets, nor do they consider physics-informed losses. }

In the current work, we consider the case where we have a large amount of low-fidelity data and either a smaller amount of higher quality high-fidelity data or knowledge of the physical laws the system obeys. We learn both the linear and nonlinear correlation between the high- and low-fidelity data, which allows for learning of complex correlations between the data sets. The low-fidelity and high-fidelity datasets do not need to be for the same set of input functions, giving additional flexibility to this method. 
We note that the framework presented here is adaptable and that the physical laws can instead be applied on the low-fidelity network, with corrections from high-fidelity data. Additionally, the framework can enforce both data and physics within a single level of fidelity. For example, if one has a low-order numerical solver that approximates a PDE and sparse measurements of the true system, both the measurements and the PDE can be enforced at the high-fidelity level and the low-order numerical data can be enforced at the low-fidelity level. 

The paper is organized as follows. We introduce the architecture and notation in Section \ref{sec:note}. We divide multifidelity DeepONets into data-driven and physics-informed. In Section \ref{sec:comp}, we consider illustrative one- and two-dimensional data-driven examples. As a motivating application to a complex problem, we apply in Section \ref{sec:ice} the data-driven multifidelity DeepONet framework with low- and high-fidelity simulations for ice-sheet dynamics. In Section \ref{sec:PI}, we investigate the performance of the physics-informed multifidelity DeepONet framework for the case of enforcing physics in the absence of high-fidelity data, \textcolor{black}{focusing on applications to Burgers equation}. Section \ref{sec:discussion} concludes with a brief discussion and suggestions for future work.

\section{Multifidelity DeepONets}\label{sec:note}

\subsection{Architecture}
A ``standard'' unstacked single fidelity DeepONet consists of two neural networks, the branch and the trunk, which are trained simultaneously \cite{lu2021learning}. The input to the branch network is a function $\mathbf{u}$ discretized at points $\{x_i\}_{i=1}^M$, and the output is  $[b_1, b_2, \ldots, b_p]^T \in \mathbb{R}^Q.$
The input to the trunk net is the coordinates $\mathbf{x} \in \mathbb{R}^n$, and the trunk output is $[t_1, t_2, \ldots, t_p]^T \in \mathbb{R}^Q.$ The DeepONet output can be expressed as
\begin{equation}
    \mathcal{G}_{SF}^\theta(\mathbf{u})(\mathbf{x}) = \sum_{k=1}^p b_kt_k
\end{equation}
where $\theta$ denotes the trainable parameters of the unstacked single fidelity DeepONet (for a detailed description see \cite{lu2019deeponet, lu2021learning}). 

In this work we use ``modified'' DeepONets, proposed in \cite{wang2021improved}. Modified DeepONets include the addition of two encoder networks, one each for the branch and trunk networks, which are included in each hidden layer of the branch and trunk networks through a convex combination. The modified DeepONet has been shown to improve performance of the method \cite{wang2021improved}. Following the notation in \cite{wang2021improved}, we denote the branch encoder by $\mathbf{U}$ with weights $\mathbf{W}_u$ and biases $\mathbf{b}_u$ and denote the trunk encoder by $\mathbf{V}$ with weights $\mathbf{W}_x$ and biases $\mathbf{b}_x$. Thus, $\mathbf{U} = \phi(\mathbf{W}_u\mathbf{u} + \mathbf{b}_u)$ and $\mathbf{V} = \phi(\mathbf{W}_x\mathbf{x} + \mathbf{b}_x)$, where $\phi$ is an activation function. We consider a modified DeepONet with $L$ layers and denote the weights and biases of the branch and trunk networks by 
$\left\{\mathbf{W}_u^{(l)}, \mathbf{b}_u^{(l)}\right\}_{l=1}^L$ and $\left\{\mathbf{W}_x^{(l)}, \mathbf{b}_x^{(l)}\right\}_{l=1}^L$, respectively. Then, the forward pass of the modified DeepONet is given by:
\begin{align}
    \mathbf{Z}_u^{(1)} &= \phi\left(\mathbf{W}_u^{(1)}\mathbf{u} +  \mathbf{b}_u^{(1)}\right) , \quad \mathbf{Z}_x^{(1)} = \phi\left(\mathbf{W}_x^{(1)}\mathbf{x} +  \mathbf{b}_x^{(1)}\right) \\ 
                \mathbf{H}_u^{(l)} &= (1-   \mathbf{Z}_u^{(l)})\odot \mathbf{U} +  \mathbf{Z}_u^{(l)}\odot \mathbf{V}, \quad
        \mathbf{H}_x^{(l)} = (1-   \mathbf{Z}_x^{(l)})\odot \mathbf{U} +  \mathbf{Z}_x^{(l)}\odot \mathbf{V}, \quad l = 1, \ldots, L-1 \\
        \mathbf{Z}_u^{(l)} &= \phi\left(\mathbf{W}_u^{(l)} \mathbf{H}_u^{(l-1)}  +  \mathbf{b}_u^{(l)}\right) , \quad \mathbf{Z}_x^{(l)} = \phi\left(\mathbf{W}_x^{(l)}\mathbf{H}_x^{(l-1)}  +  \mathbf{b}_x^{(l)}\right), \quad l = 2, \ldots, L-1 \\
        \mathbf{H}_u^{(L)} &= \mathbf{W}_u^{(L)} \mathbf{H}_u^{(L-1)}  +  \mathbf{b}_u^{(L)} , \quad \mathbf{H}_x^{(L)} = \mathbf{W}_x^{(L)} \mathbf{H}_x^{(L-1)}  +  \mathbf{b}_x^{(L)}  \\
G^\theta(\mathbf{u})(\mathbf{y}) &= \left\langle \mathbf{H}_u^{(L)} , \mathbf{H}_x^{(L)}  \right\rangle
\end{align}
Here, $\odot$ represents point-wise multiplication and $\langle \cdot, \cdot \rangle$ represents an inner product. 
Note that for single fidelity DeepONets no activation function is applied to the last layer. 

\begin{figure}
\centering
\includegraphics[width=\textwidth]{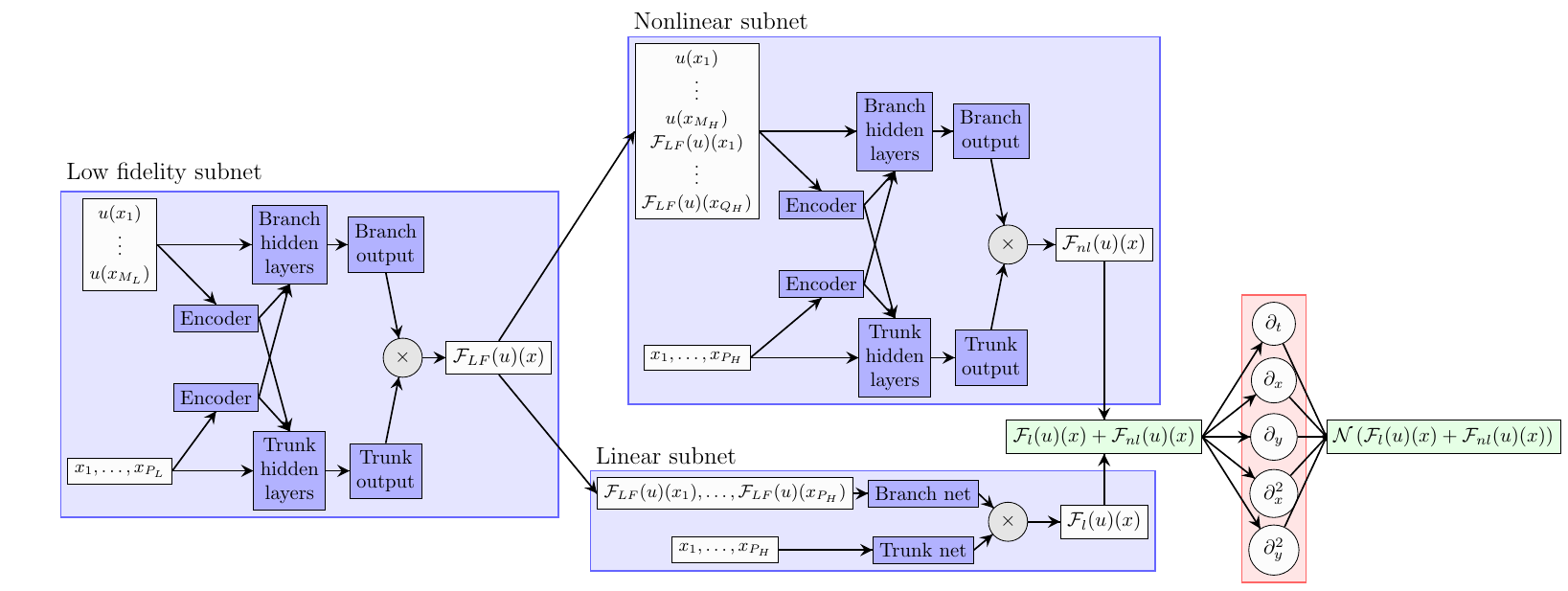}\label{fig:PI-setup}
\caption{Schematic of the composite physics-informed multifidelity DeepONet setup. $\mathcal{F}_{LF}({u})(x),$ $\mathcal{F}_{nl}({u})(x)$, and $\mathcal{F}_{l}({u})(x)$ are the outputs of the low-fidelity, nonlinear, and linear DeepONet subnets. } \label{fig:setup}
\end{figure}

The multifidelity DeepONet framework consists of three blocks, trained simultaneously. 
A schematic of the multifidelity DeepONet architecture is given in Fig. \ref{fig:setup} for both data-driven and physics-informed cases. The low-fidelity block is a standard modified DeepONet \cite{wang2021improved}.
It represents an approximation of the low-fidelity data. The nonlinear block encodes the nonlinear correlation between the output of the low-fidelity network and the high-fidelity data or physics. The nonlinear block is a modified DeepONet with the activation function applied on every layer. The linear block is a standard DeepONet with no activation functions and is used to approximate the linear correlation between the output of the low-fidelity block and the high-fidelity data or physics. We note that the outputs of the linear and nonlinear DeepONets are continuously differentiable with respect to their input coordinates, and therefore we can use automatic differentiation on the outputs \cite{griewank1989automatic, baydin2018automatic}. If a low-fidelity solver is available, the solver can be used in place of the low-fidelity block. A discussion of this case is given in \ref{sec:noncomp}. 

\subsection{Notation}
Consider a nonlinear operator $G$ mapping from one space of functions to another space of functions, $G : \mathcal{U}  \rightarrow \mathcal{S}$. It has been shown that a neural network with a single hidden layer can accurately approximate the operator $G$ \cite{chen1995approximation, chen1995universal, lu2021learning}. Denote the input function to the operator by $u \in \mathcal{U}$, and denote the output function by $G(u) \in \mathcal{S}$. For any point $x$ in the domain our network takes input $(u, x)$ and outputs $G(u)(x).$

Now, following the notation in \cite{wang2021improved}, we consider a linear or nonlinear differential operator $\mathcal{N} : \mathcal{U} \times \mathcal{S} \rightarrow \mathcal{V} $, where $(\mathcal{U}, \mathcal{S}, \mathcal{V})$ is a triplet of Banach spaces, and a general parametric PDE of the form $\mathcal{N} (u, s) = 0$ with boundary conditions $\mathcal{B}(u, s) = 0$. Here,  $u\in \mathcal{U} $ is the input function, and $s \in \mathcal{S}$ is a function satisfying the differential operator subject to the boundary conditions. $\mathcal{B}$ can represent Dirichlet, Neumann, Robin, or periodic boundary conditions. With this notation, if for all $u \in \mathcal{U}$ there exists a unique solution $s = s(u)$ to $\mathcal{N} (u, s) = 0$ and $\mathcal{B} (u, s) = 0$, the solution can be represented as an operator $G : \mathcal{U} \rightarrow \mathcal{S}$, where $G(u) = s(u)$ \cite{wang2021improved}.

A schematic of the multifidelity DeepONet architecture is given in Fig. \ref{fig:setup} for both data-driven and physics-informed cases. In both cases, we train three DeepONets simultaneously. The low-fidelity subnet approximates the low-fidelity operator. We assume that we have low-fidelity data with inputs to the operators given by $\mathbf{u}^j \in \mathcal{U}$ for $j = 1, \ldots, N_L.$ While each $\mathbf{u}^j$ can be a continuous function, we need to discretize $\mathbf{u}^j$ and therefore evaluate $\mathbf{u}^j$ on a set of values called sensors for input to the low-fidelity branch net, given by $\left(\mathbf{u}^j(x_{1}), \ldots, \mathbf{u}^j(x_{M_L})\right).$ The output values are given by $\left(y_L(\mathbf{u}^j)(x_{1}), \ldots, y_L(\mathbf{u}^j)(x_{P_L})\right)_{j = 1, \ldots, N_L}  \in \mathcal{S}$. We denote the output of the low-fidelity, nonlinear, and linear subnets by $\mathcal{F}_{LF}(\mathbf{u})(x),$ $\mathcal{F}_{nl}(\mathbf{u})(x)$, and $\mathcal{F}_{l}(\mathbf{u})(x)$, respectively. 
The nonlinear and linear subnets take the output of the low-fidelity subnet as part or all of their branch networks, and the sum of the linear and nonlinear subnet outputs approximates the high-fidelity operator. The input to the linear branch net is given by $[\mathcal{F}_{LF}^\theta(\mathbf{u})(\mathbf{x}_1), \ldots, \mathcal{F}_{LF}^\theta(\mathbf{u})(\mathbf{x}_{P_H})]$, that is, the input is the low-fidelity network evaluated at all points used in the input to the linear trunk net. The input to the nonlinear branch net is given by $[\mathbf{u}(\mathbf{x}_1), \ldots, \mathbf{u}(\mathbf{x}_{M_H}), \mathcal{F}_{LF}^\theta(\mathbf{u})(\mathbf{x}_1), \ldots, \mathcal{F}_{LF}^\theta(\mathbf{u})(\mathbf{x}_{Q_H})]^T$.
In this way, the nonlinear subnet learns the nonlinear correlation between the low- and high-fidelity operators, and the linear subnet learns the linear correlation between the low- and high-fidelity operators. Typically, we take $Q_H = P_L$ and augment the nonlinear branch net input with the output of the low-fidelity network at all low-fidelity training points, however, in some cases choosing a different set of output points can improve performance by reducing the number of points needed in the evaluation and therefore increasing speed, or by highlighting important features. 
We denote by $\mathbf{\theta} = (\theta_{nl}, \theta_{l}, \theta_{LF})$ the set of all trainable parameters of the three subnetworks. 

The low-fidelity loss for a given input function $\mathbf{u}^j$ is given by 
\begin{equation}
\mathcal{L}_{LF}(\mathbf{u}^j, \theta_{LF}) = \frac{1}{P_L} \sum_{k=1}^{P_L} \left| y_L(  \mathbf{u}^j)(x_k) - \mathcal{F}_{LF}^\theta(\mathbf{u}^j)(x_k)\right|^2
\end{equation}
and the full low-fidelity loss is then 
\begin{equation}
    \mathcal{L}_{LF}(\theta_{LF}) = \frac{1}{N_L} \sum_{j=1}^{N_L} \mathcal{L}_{LF}(\mathbf{u}^j, \theta_{LF}) = \frac{1}{N_LP_L} \sum_{j=1}^{N_L} \sum_{k=1}^{P_L} \left| y_L(  \mathbf{u}^j)(x_k) - \mathcal{F}_{LF}^\theta(\mathbf{u}^j)(x_k)\right|^2. 
\end{equation}

For the high-fidelity networks, we divide into two cases depending on whether high-fidelity data or physical knowledge of the system is enforced. 

\subsubsection{Data-driven multifidelity notation} \label{sec:DD_note}

In the data-driven case, we assume we have high-fidelity data with input values to the operator $\mathbf{u}^k \in \mathcal{U}$ for $k = 1, \ldots, N_H$ evaluated on sensors $\left(\mathbf{u}^k(x_{1}), \ldots, \mathbf{u}^k(x_{M_H})\right).$ We note that we do not need the high-fidelity input values to be the same as the low-fidelity values, that is $\{\mathbf{u}^k\}_{k=1}^{N_H} $ and $\{\mathbf{u}^j\}_{j=1}^{N_L} $ can be independent sets. The output values are given by $\left(y_H(\mathbf{u}^k)(x_1), \ldots, y_H(\mathbf{u}^k)(x_{P_H})\right)\in \mathcal{S}$. We simultaneously train three DeepONets to learn the low-fidelity operator and the linear and nonlinear correlations between the low-fidelity operator and the high-fidelity operator. 
The full high-fidelity loss is 
\begin{equation}
    \mathcal{L}_{HF}(\theta_{nl}, \theta_{l}) = \frac{1}{N_H P_H} \sum_{j=1}^{N_H} \sum_{k=1}^{P_H} \left| y_H(  \mathbf{u}^j)(x_k) - \mathcal{F}_{nl}^\theta(\mathbf{u}^j)(x_k) - \mathcal{F}_{l}^\theta(\mathbf{u}^j)(x_k)\right|^2.
\end{equation}
The full data-driven loss is then given by: 
\begin{equation}
\mathcal{L}_{DD}(\theta) = \lambda_1 \mathcal{L}_{HF}(\theta_{nl}, \theta_{l}) + \lambda_2 \mathcal{L}_{LF}(\theta_{LF}) + \lambda_3 \left( \sum w_{nl}^2 + \sum b_{nl}^2 \right) + \lambda_4 \left( \sum w_{LF}^2 + \sum b_{LF}^2 \right) . \label{eq:loss_data}
\end{equation}
Here, $w_{nl}$ and $b_{nl}$ are the weights and biases from the nonlinear branch net and $w_{LF}$ and $b_{LF}$ are the weights and biases from the low-fidelity branch net. The regularization term on the nonlinear branch net serves to minimize the nonlinear correlation, forcing the network to learn a linear correlation if appropriate. The regularization term on the low-fidelity branch net prevents over-fitting of the low-fidelity data. $\lambda_i$, $i = 1, 2, 3, 4$ are weights that can be chosen for each case. 

We compare our results with a single fidelity modified DeepONet trained on the high-fidelity data. 
The full single fidelity data-driven loss is 
\begin{equation}
    \mathcal{L}_{SF}(\theta_{SF}) = \frac{1}{N_HP_H} \sum_{j=1}^{N_H} \sum_{k=1}^{P_H} \left| y_H(  \mathbf{u}^j)(x_k) - \mathcal{G}^\theta_{SF}(\mathbf{u}^j)(x_k) \right|^2,
\end{equation}
where $\mathcal{G}^\theta_{SF}(\mathbf{u}^j)(x_k)$ denotes the output from the single fidelity modified DeepONet with parameters $\theta_{SF}$.

\subsubsection{Physics-informed multifidelity formulation}
In the second case, we assume we have no high-fidelity training data in the form of  input-output pairs for the high-fidelity network, but that the output does satisfy a PDE with appropriate boundary conditions. 
The boundary condition loss for a given input function $\left\{\mathbf{u}^i\right\}_{i=1}^{N_H}$
can be written as: 
\begin{equation}
\mathcal{L}_{\mathcal{B}}(\mathbf{u}^i, \theta_{nl}, \theta_{l})  = \frac{1}{P_{BC}} \sum_{j=1}^{P_{BC}}    \left| \mathcal{B}\left(\mathbf{u}^i, \mathcal{F}^\theta_{nl}(\mathbf{u}^i)(\mathbf{x}_j)+ \mathcal{F}^\theta_{l}(\mathbf{u}^i)(\mathbf{x}_j)\right) \right|
\end{equation}
where the points $\left\{\mathbf{x}_j\right\}_{j=1}^{P_{BC}}$ are randomly chosen on the boundary of the domain. 
The full boundary condition loss is then: 
\begin{equation}
\mathcal{L}_{\mathcal{B}}( \theta_{nl}, \theta_{l})  = \frac{1}{{N_H}P_{BC}} \sum_{i=1}^{N_H}\sum_{j=1}^{P_{BC}}    \left| \mathcal{B}\left(\mathbf{u}^i, \mathcal{F}^\theta_{nl}(\mathbf{u}^i)(\mathbf{x}_j)+ \mathcal{F}^\theta_{l}(\mathbf{u}^i)(\mathbf{x}_j)\right) \right|.
\end{equation}
We treat initial conditions as a special case of boundary conditions. 

We also consider the loss in satisfying the parametric PDE, given by:
\begin{equation}
\mathcal{L}_{physics}( \theta_{nl}, \theta_{l})  = \frac{1}{{N_H}P_p} \sum_{i=1}^{N_H}\sum_{j=1}^{P_p}    \left| \mathcal{N}\left(\mathbf{u}^i, \mathcal{F}^\theta_{nl}(\mathbf{u}^i)(\mathbf{x}_j)+\mathcal{F}^\theta_{l}(\mathbf{u}^i)(\mathbf{x}_j)\right) \right|,
\end{equation}
where the points $\left\{\mathbf{x}_j\right\}_{j=1}^{P_p}$ are randomly chosen on the interior of the domain. 

For ease of notation, we split the boundary condition operator into two parts,  the part representing the boundary conditions ($\mathcal{L}_{BC}$) and the part representing the initial condition ($\mathcal{L}_{IC}$). Then, the full multifidelity physics-informed loss can be written as: 
\begin{align}
\mathcal{L}_{PI}(\theta) = &\lambda_1 \mathcal{L}_{physics}(\theta_{nl}, \theta_{l}) + \lambda_2 \mathcal{L}_{LF}(\theta_{LF}) +\lambda_5 \mathcal{L}_{IC}(\theta_{nl}, \theta_{l})+\lambda_6 \mathcal{L}_{BC}(\theta_{nl}, \theta_{l}) \nonumber \\
&+ \lambda_3 \left( \sum w_{nl}^2 + \sum b_{nl}^2 \right) + \lambda_4 \left( \sum w_{LF}^2 + \sum b_{LF}^2 \right) . \label{eq:loss_physics}
\end{align}

We compare our results with a single fidelity physics-informed DeepONet. 
The full single fidelity physics-informed loss is 
\begin{align}
\mathcal{L}_{PI}(\theta) = &\lambda_1 \mathcal{L}_{physics}(\theta) +\lambda_5 \mathcal{L}_{IC}(\theta)+\lambda_6 \mathcal{L}_{BC}(\theta).  \label{eq:loss_physics_SF}
\end{align}
While we consider only the case where physics is applied as a high-fidelity correction to low-fidelity data, it is also possible to consider a case where physics represents a low-fidelity model, such as when the exact physics of the system is not known. In that case, the low-fidelity physics-informed DeepONet could use high-fidelity data from, \emph{e.g.} experiments, together as a training set.

\section{Data-driven multifidelity DeepONets}\label{sec:comp}
\textcolor{black}{In this section we discuss data-driven multifidelity training cases for data-driven problems in Subsections \ref{1d_jump} and \ref{2d_nonlinear}, and for the two-dimensional problem of ice-sheet modeling in Subsection \ref{sec:ice}. For reference, two additional examples are given in \ref{sec:add_data-drive}. In data-driven training, we consider cases with low- and high-fidelity data. The low fidelity data is abundant, but has lower accuracy. In traditional numerical approaches to PDEs, highly accurate data is generally generated on a finer mesh, so high fidelity data has high resolution. However, this is not necessarily the case in other applications. We can consider a case where it is possible to cheaply collect many data points with lower accuracy, perhaps using a less accurate instrument that is cheap to operate, and then gather fewer high fidelity data points with very expensive measurements. Thus, the low-fidelity data may be at higher resolution, but with lower accuracy, as in the cases considered in Subsections \ref{1d_jump} and \ref{1d_corr}. The task of the multifidelity DeepONet is to learn the correlations between the low-fidelity inaccurate data and the sparse, but accurate, high-fidelity data.}

\subsection{One-dimensional, jump function}\label{1d_jump}
\begin{figure}[ht]
\centering
\begin{subfigure}{0.4\textwidth}
\centering
\caption{$a=11.1526$}
\includegraphics[width=\textwidth]{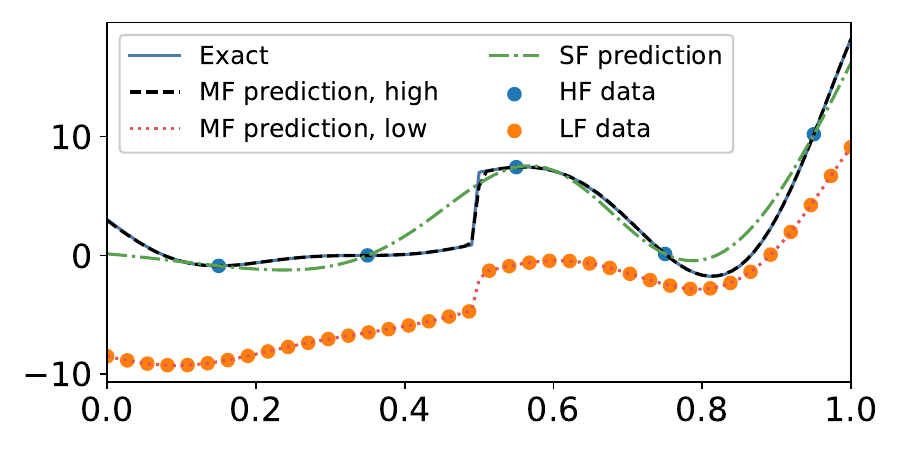}
\end{subfigure}\hspace{0.02\textwidth}
\begin{subfigure}{0.4\textwidth}
\centering
\caption{$a=13.2579$}
\includegraphics[width=\textwidth]{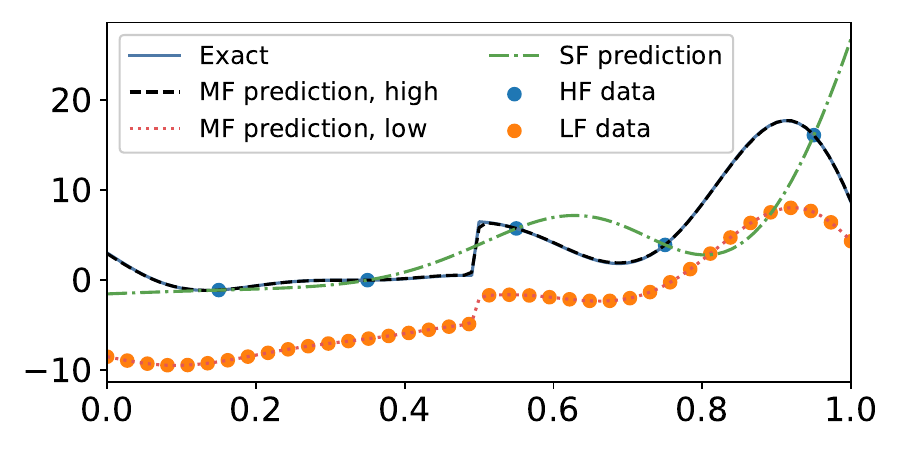}
\end{subfigure}\\
\begin{subfigure}{.3\textwidth}
\centering
\caption{SF error}
\includegraphics[width=\textwidth]{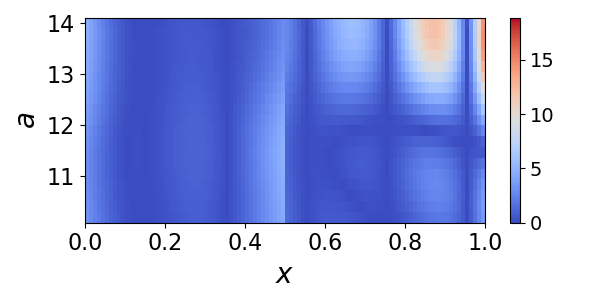}\label{fig:jumpSF}
\end{subfigure}
\begin{subfigure}{.3\textwidth}
\centering
\caption{MF error, high-fidelity data}
\includegraphics[width=\textwidth]{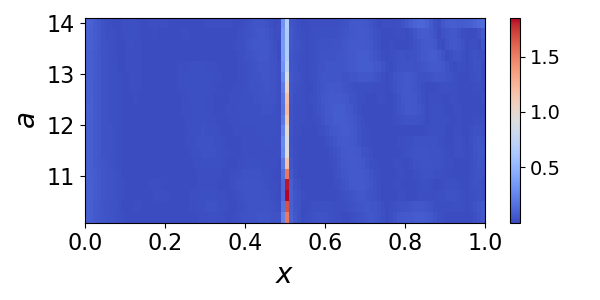}
\end{subfigure}
\begin{subfigure}{.3\textwidth}
\centering
\caption{MF error, low-fidelity data}
\includegraphics[width=\textwidth]{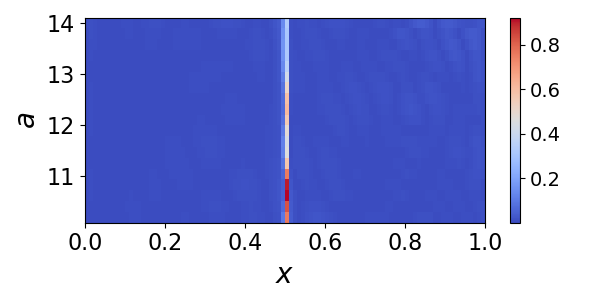}
\end{subfigure}
\caption{Data-driven multifidelity: one-dimensional, jump function. (a-b) Results of the single fidelity and multifidelity predictions of the high- and low-fidelity data.  (c) Single-fidelity error as a function of $a$ and $x$, (d) multifidelity high-fidelity prediction error as a function of $a$ and $x$, and (e) multifidelity low-fidelity prediction error as a function of $a$ and $x$.} \label{fig:1d_jump_fig}
\end{figure}

We first consider a case where the low- and high-fidelity data are represented by jump functions with a linear correlation. We show that we can recover both the jump and the linear correlation accurately. The low- and high-fidelity data are given by: 
\begin{align}
    y_L(u)(x) &= \begin{cases} 
      0.5(6x-2)^2\sin(u) + 10(x-0.5) - 5 & \; x\leq  0.5\\
      0.5(6x-2)^2\sin(u) + 10(x-0.5) - 2 & \; x > 0.5 
   \end{cases} \\
    y_H(u)(x) &=2y_L(u)(x) -20x + 20  \\
    u &= ax-4
\end{align}
for $x \in [0, 1]$ and $a \in [10, 14]$.  Parameters are given in Tab. \ref{tab:train_params_all} (see \ref{sec:training_params} and \ref{1d_jump_app}) and results in Fig. \ref{fig:1d_jump_fig}. The learned linear correlation is: 
\begin{equation}
    \mathcal{F}_l(u)(x) = 1.9479 \mathcal{F}_{LF}(u)(x)-19.1719x+19.3459 - 0.04870 x \mathcal{F}_{LF}(u)(x). 
\end{equation}
The learned linear correlation accurately captures the exact correlation. Fig. \ref{fig:jumpSF} shows that the single fidelity method fails to capture the jump, with a large error at $x = 0.5$, and the absolute error is quite large across the domain. With the multifidelity method, both the high-fidelity and low-fidelity errors are smaller. The errors are concentrated at the jump due to the limitations of the resolution of the low-fidelity training set. 

\subsection{Two-dimensional, nonlinear correlation}\label{2d_nonlinear}
\begin{figure}[ht]
\begin{subfigure}{\textwidth}
\centering
\caption{$a = 8.5211$}
\includegraphics[width=.25\textwidth]{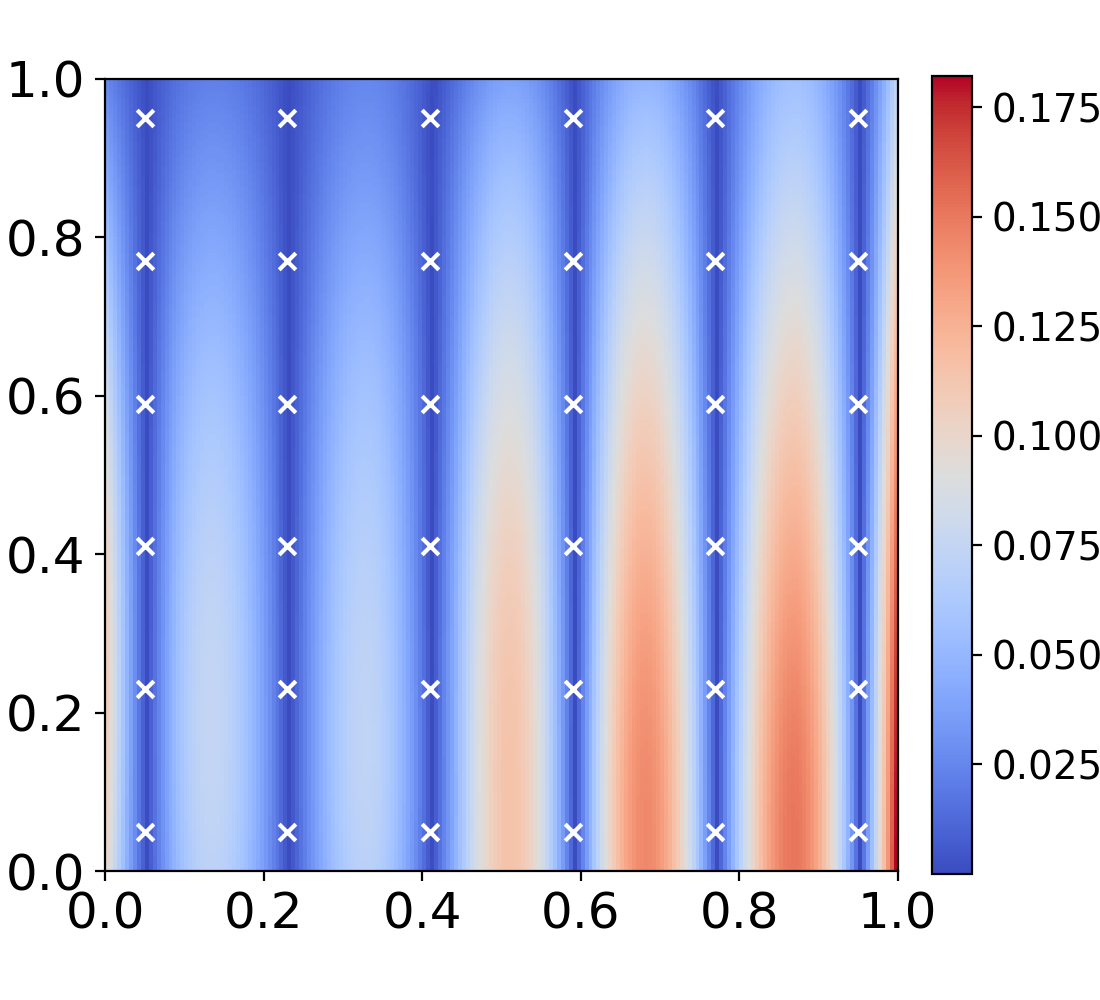}
\hspace{.05\textwidth}
\includegraphics[width=.25\textwidth]{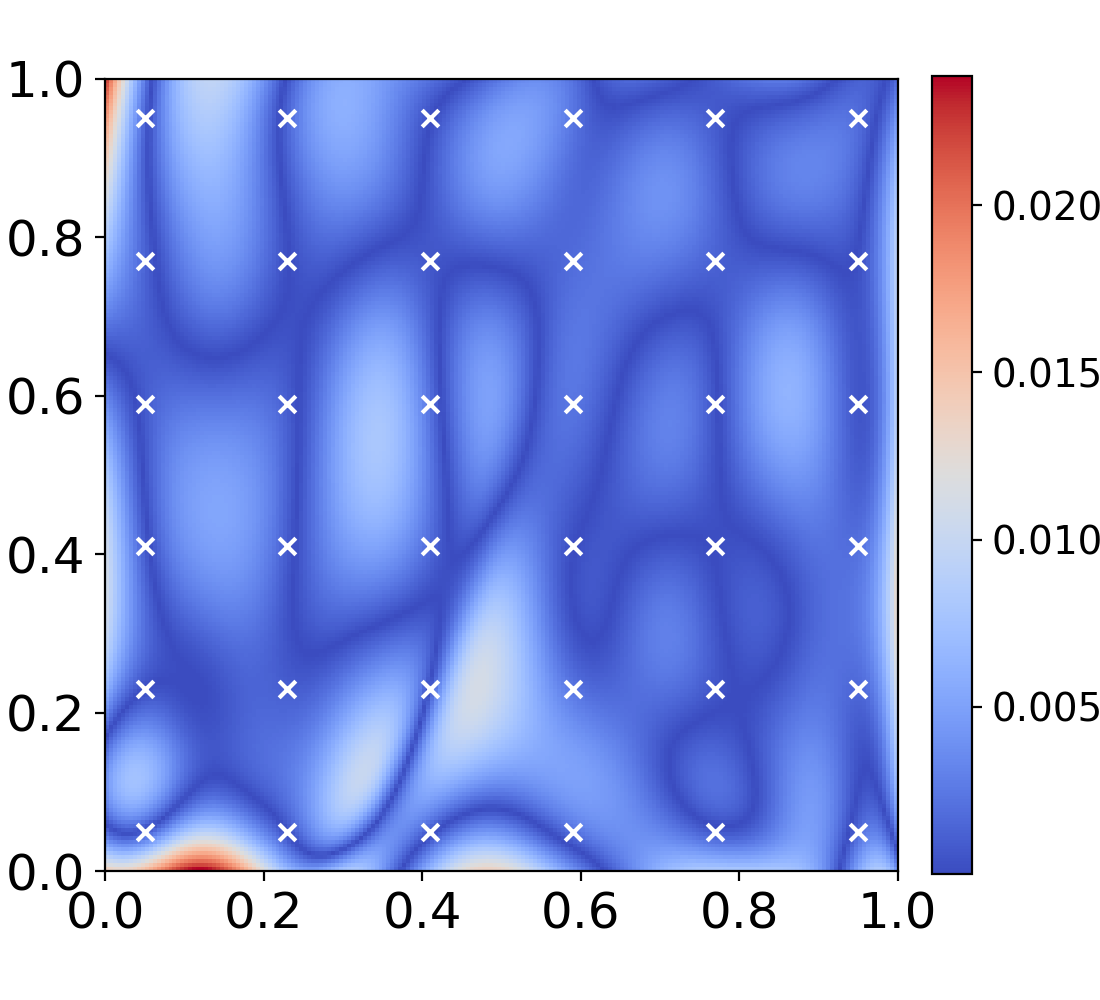}
\hspace{.05\textwidth}
\includegraphics[width=.25\textwidth]{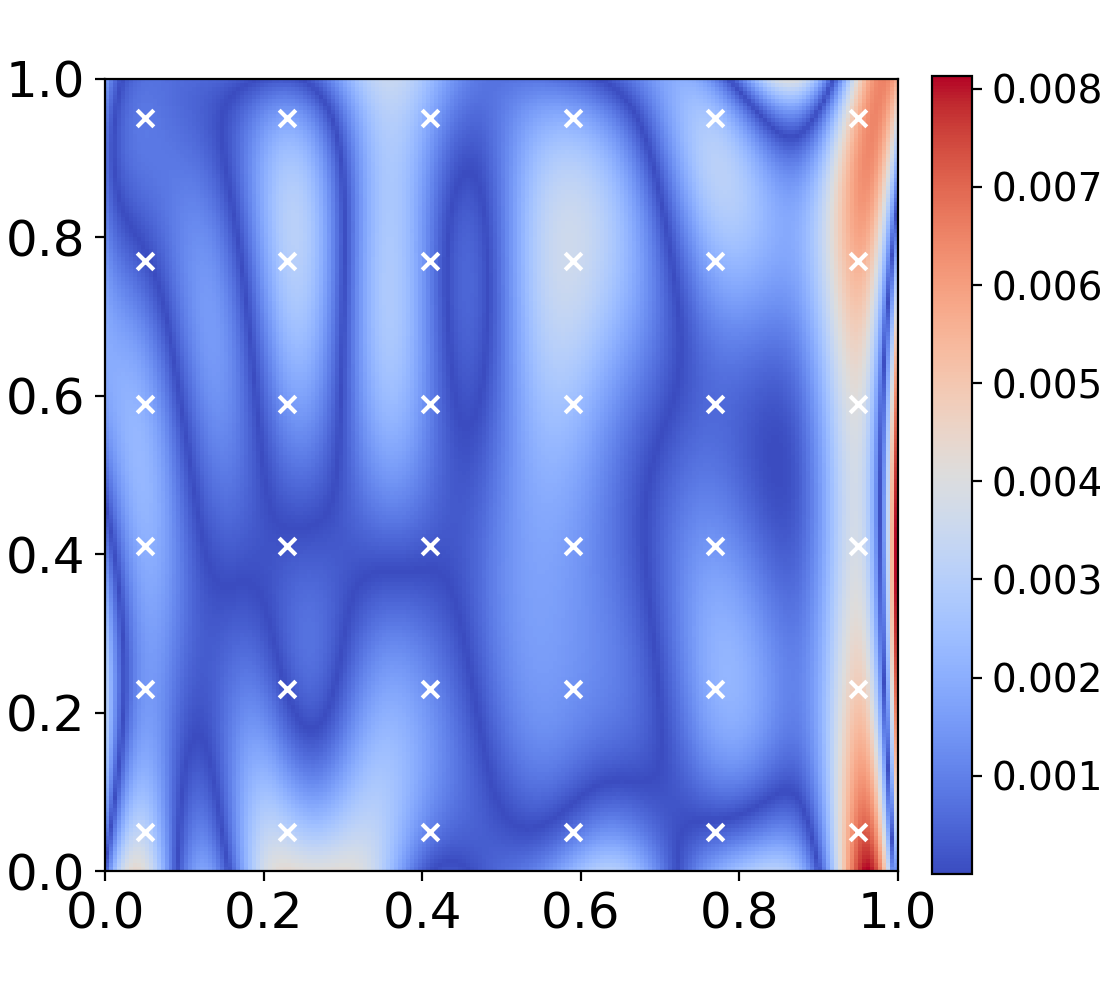}
\end{subfigure}
\begin{subfigure}{\textwidth}
\centering
\caption{$a = 9.5737$}
\includegraphics[width=.25\textwidth]{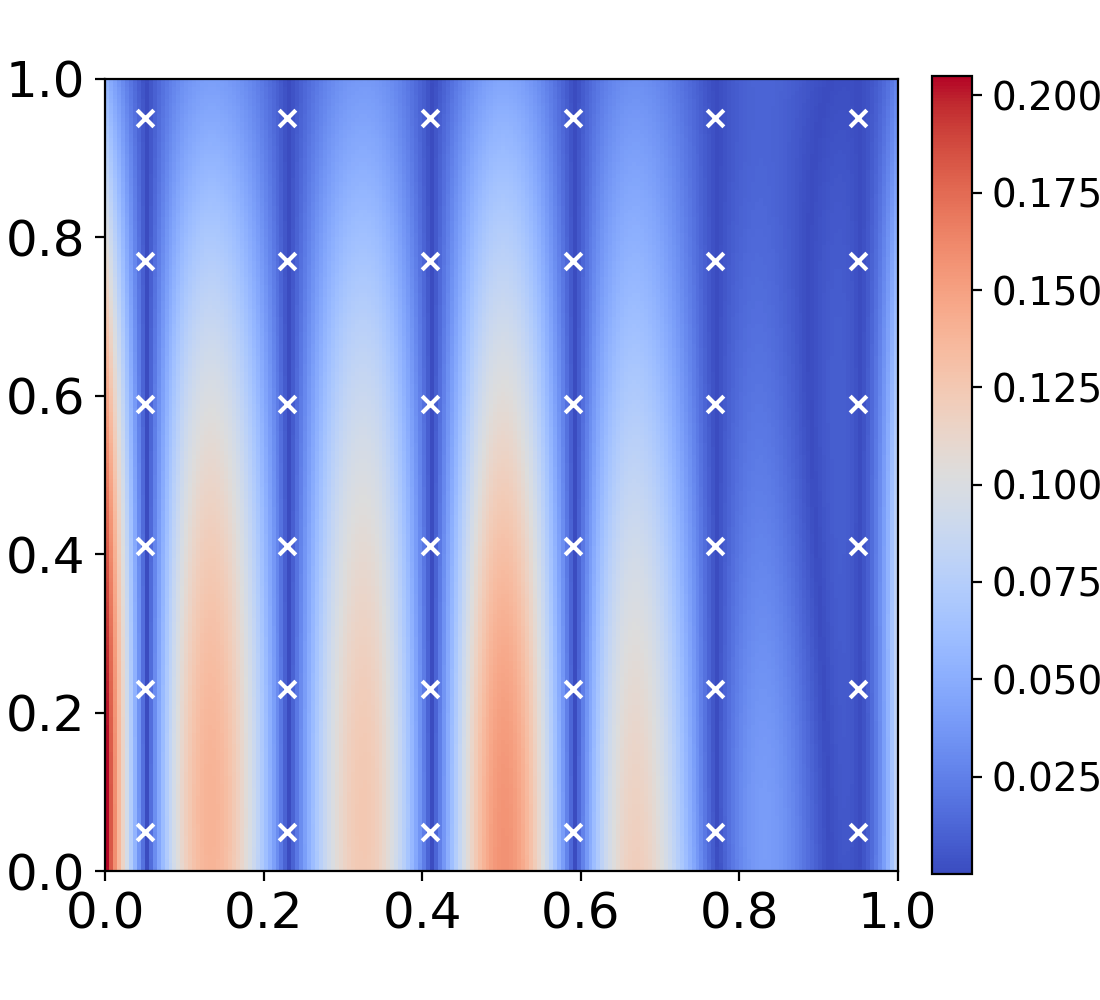}
\hspace{.05\textwidth}
\includegraphics[width=.25\textwidth]{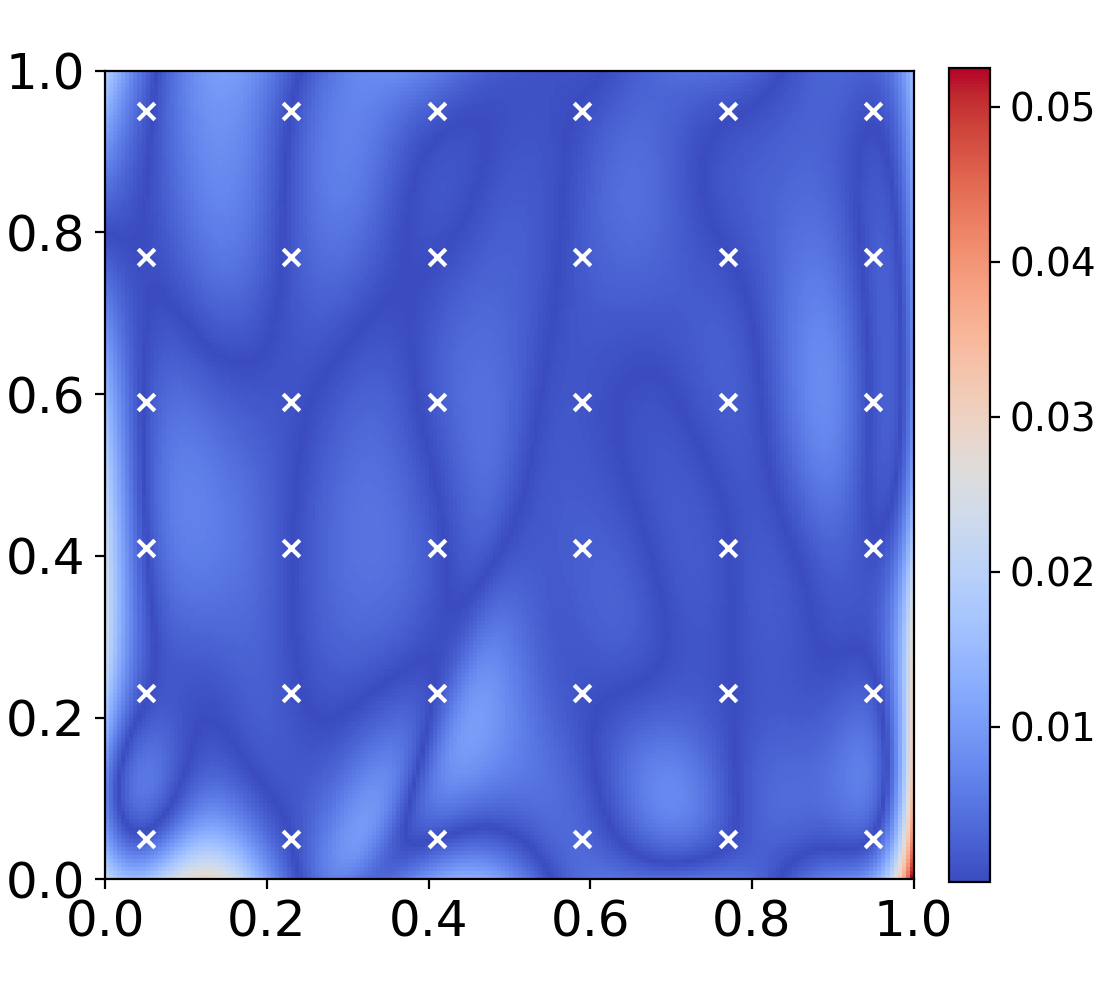}\hspace{.05\textwidth}
\includegraphics[width=.25\textwidth]{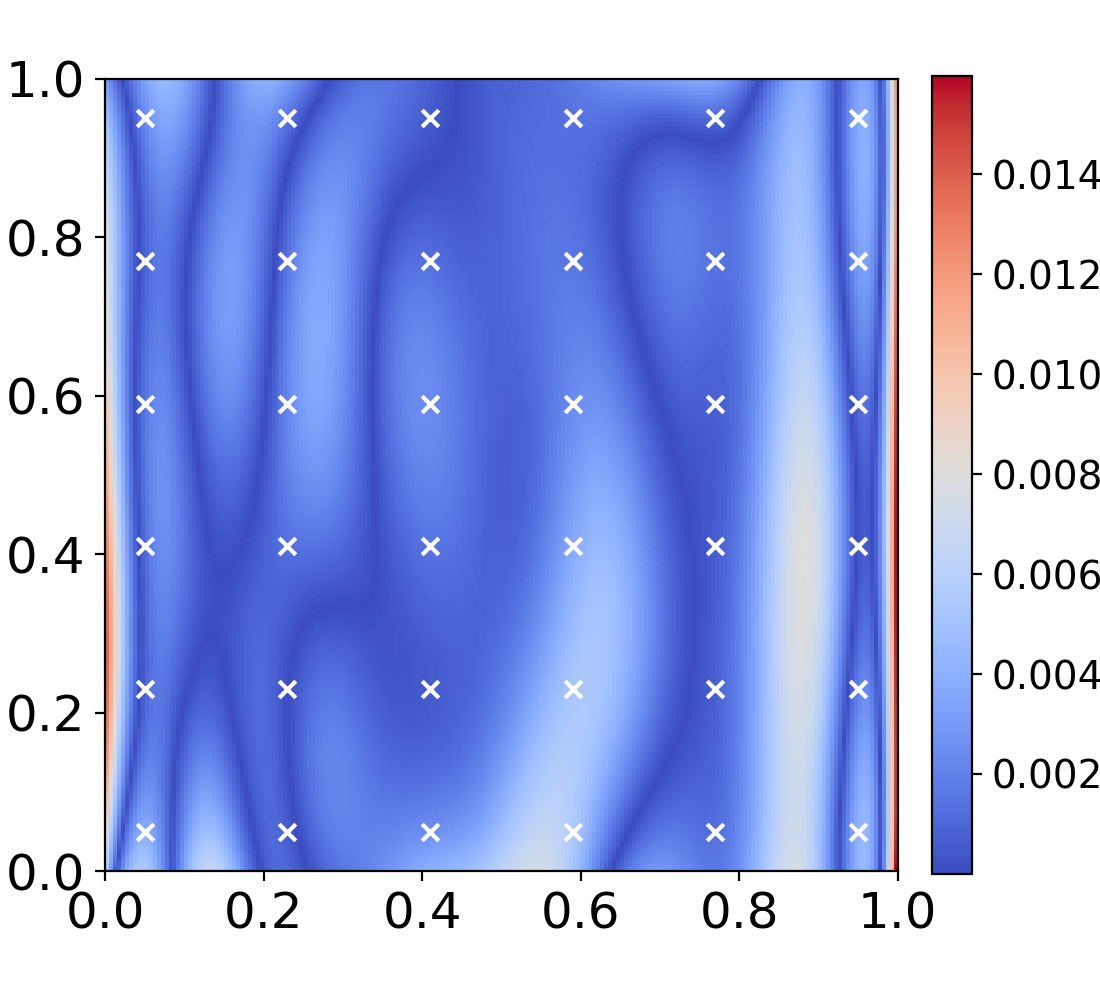}
\end{subfigure}
\caption{Data-driven multifidelity: two-dimensional, nonlinear correlation. (a) Absolute error of the high-fidelity prediction, multifidelity prediction of the high-fidelity data, and multifidelity prediction of the low-fidelity data for $a = 8.5211$. (b) Absolute error of the high-fidelity prediction, multifidelity prediction of the high-fidelity data, and multifidelity prediction of the low-fidelity data for $a = 9.5737$.   The high-fidelity data points are shown in white for clarity.} \label{fig:2d_nonlin_results} 
\end{figure}

\begin{figure}[h!]
\begin{subfigure}{\textwidth}
\centering
\caption{$a = 8.5211$}
\includegraphics[width=.22\textwidth]{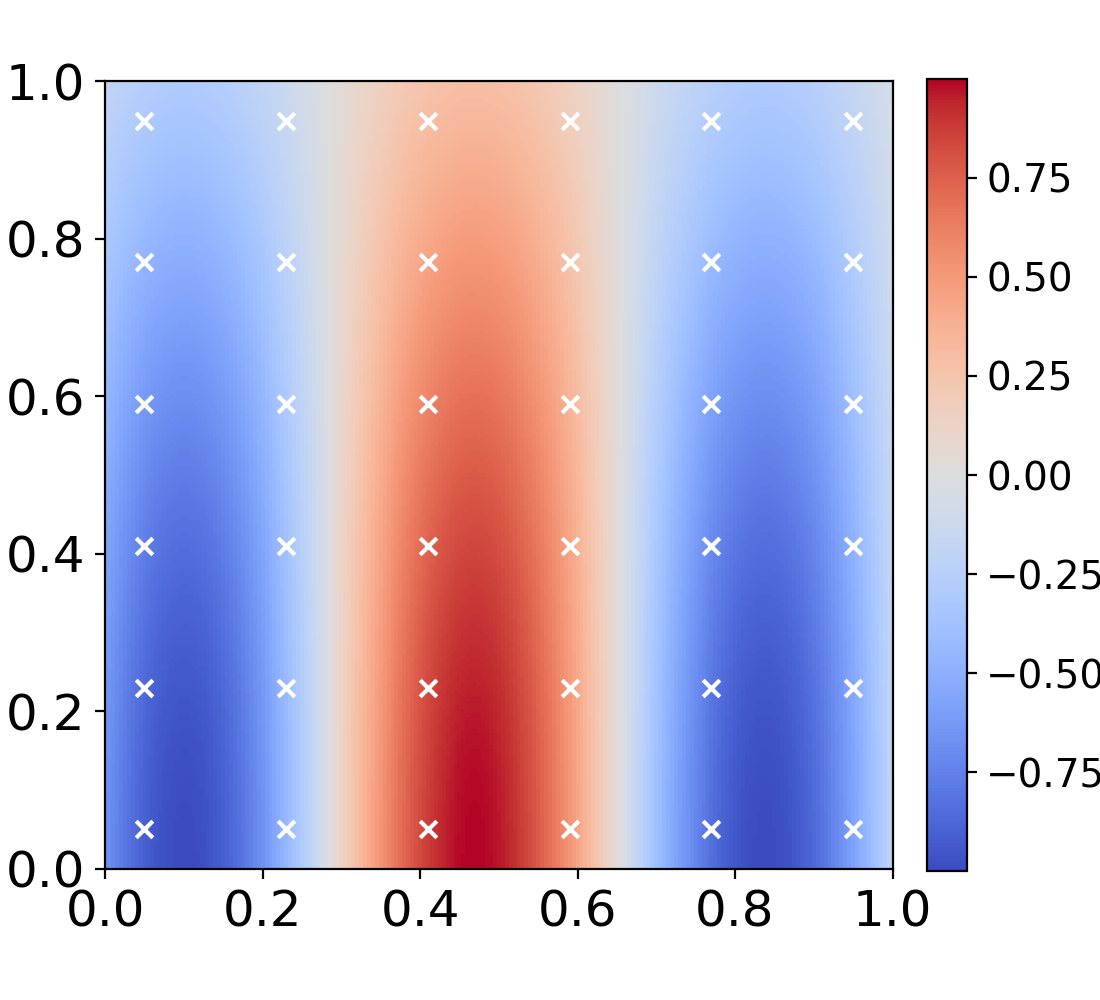}
\includegraphics[width=.22\textwidth]{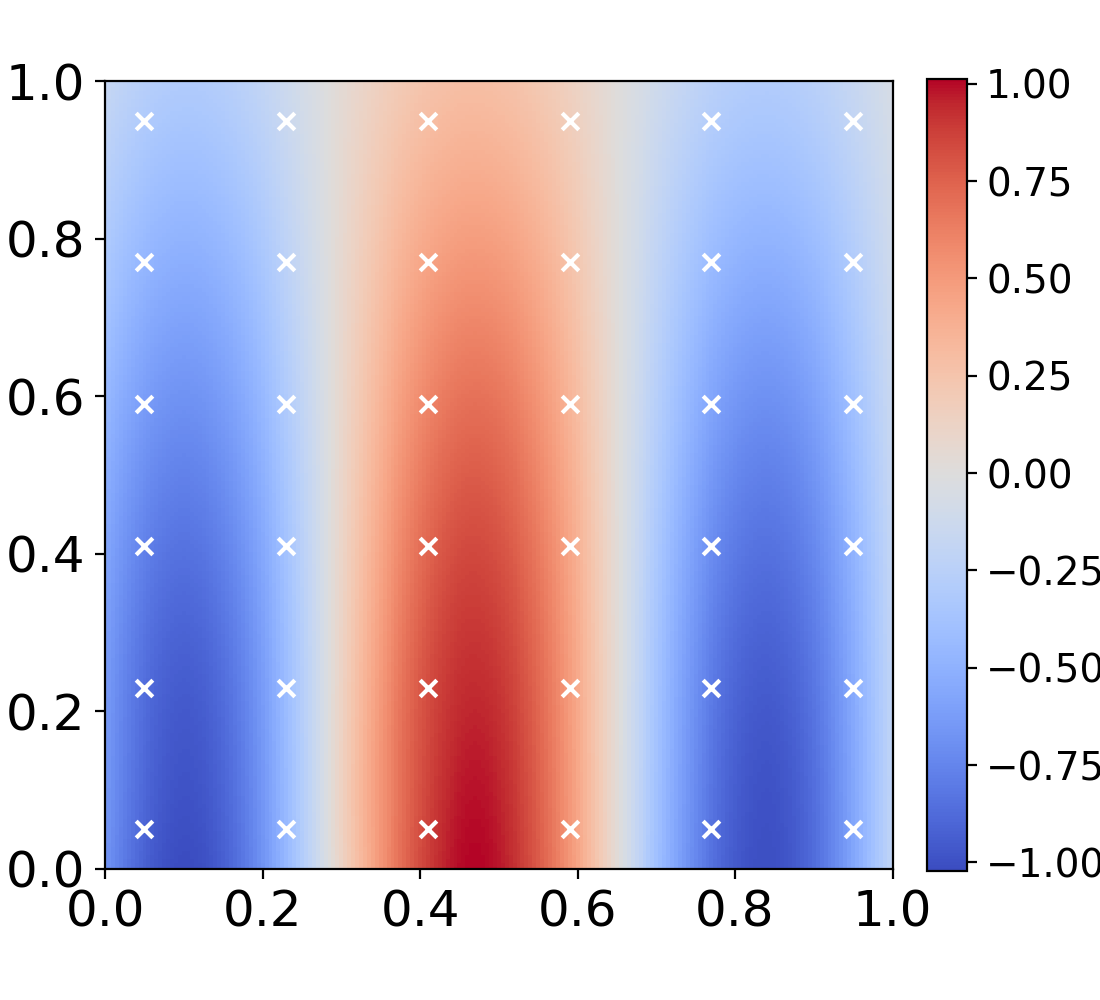}
\includegraphics[width=.22\textwidth]{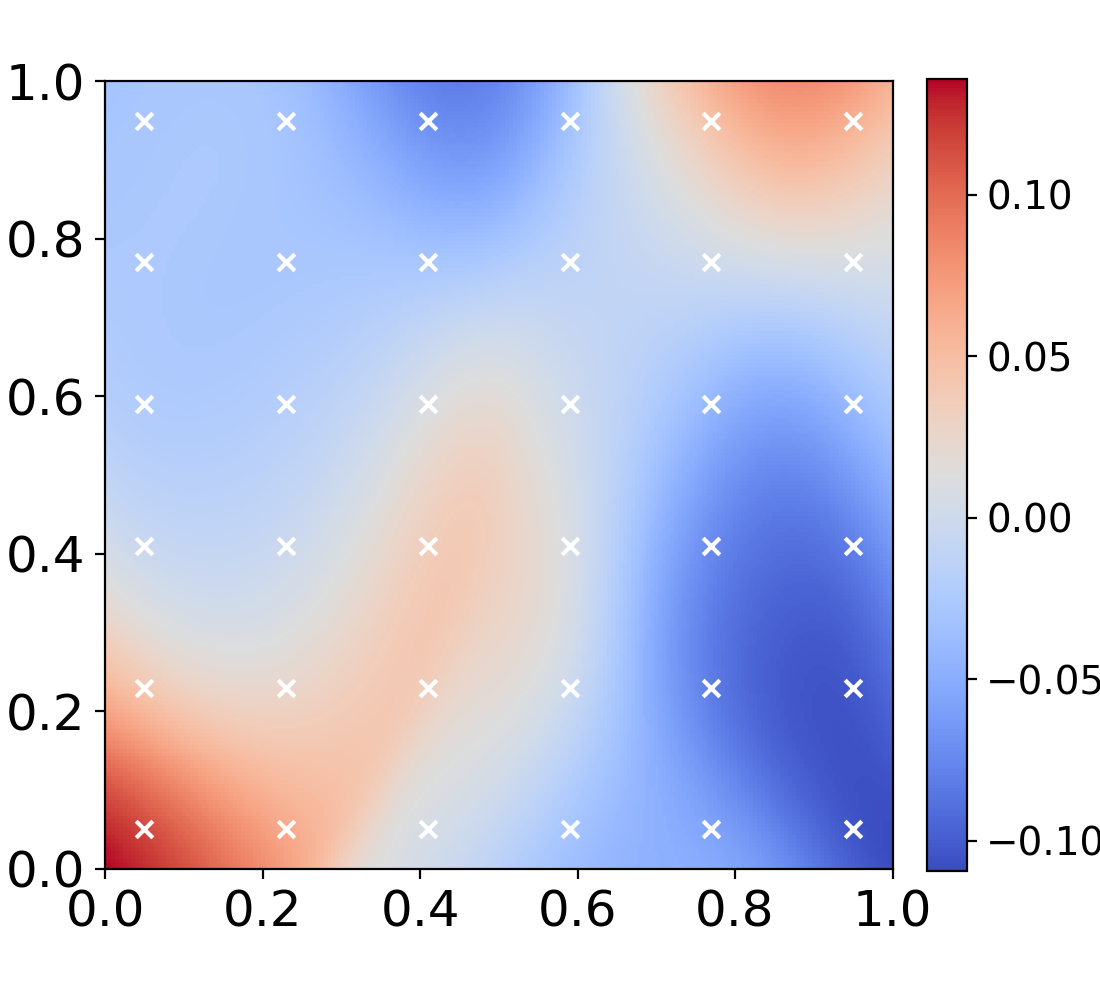}
\includegraphics[width=.22\textwidth]{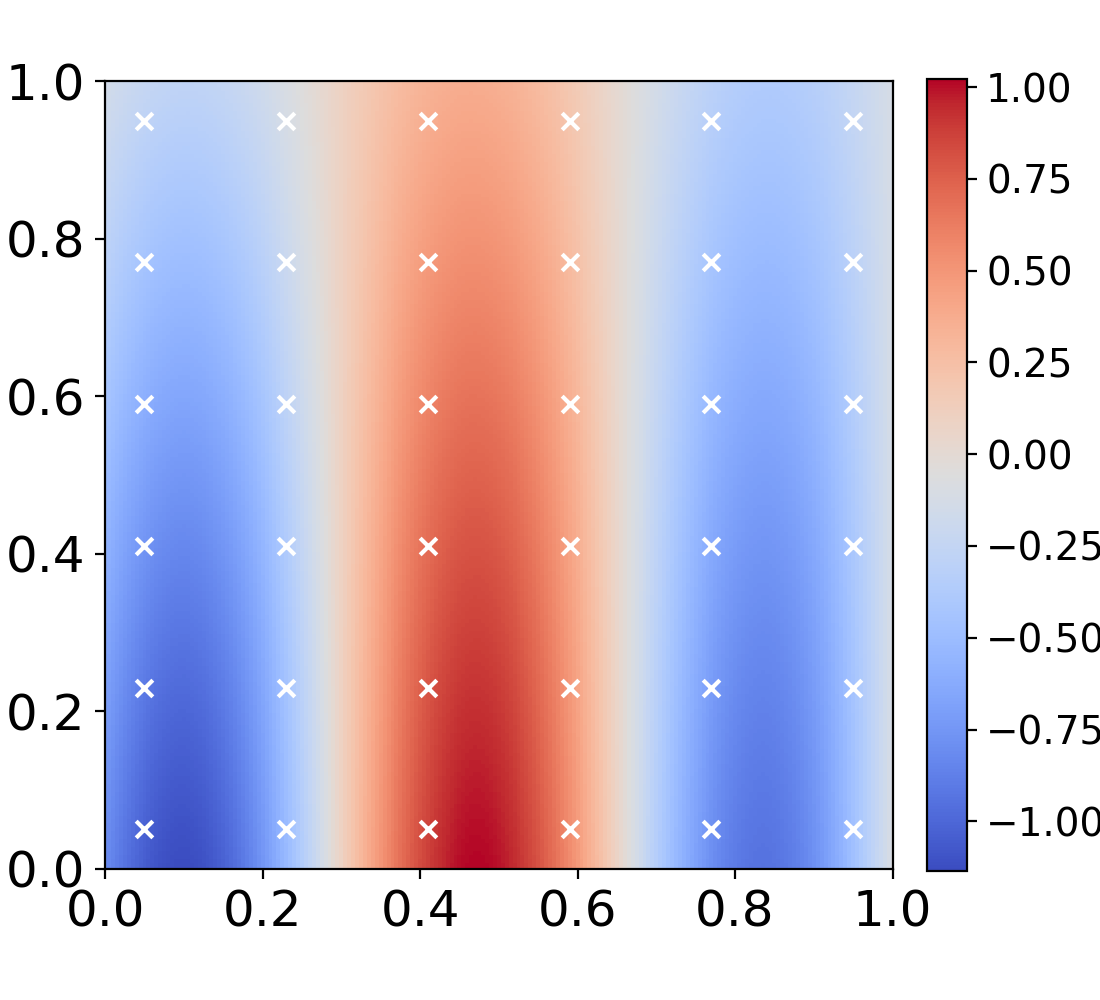}
\end{subfigure}
\begin{subfigure}{\textwidth}
\centering
\caption{$a = 9.5737$}
\includegraphics[width=.22\textwidth]{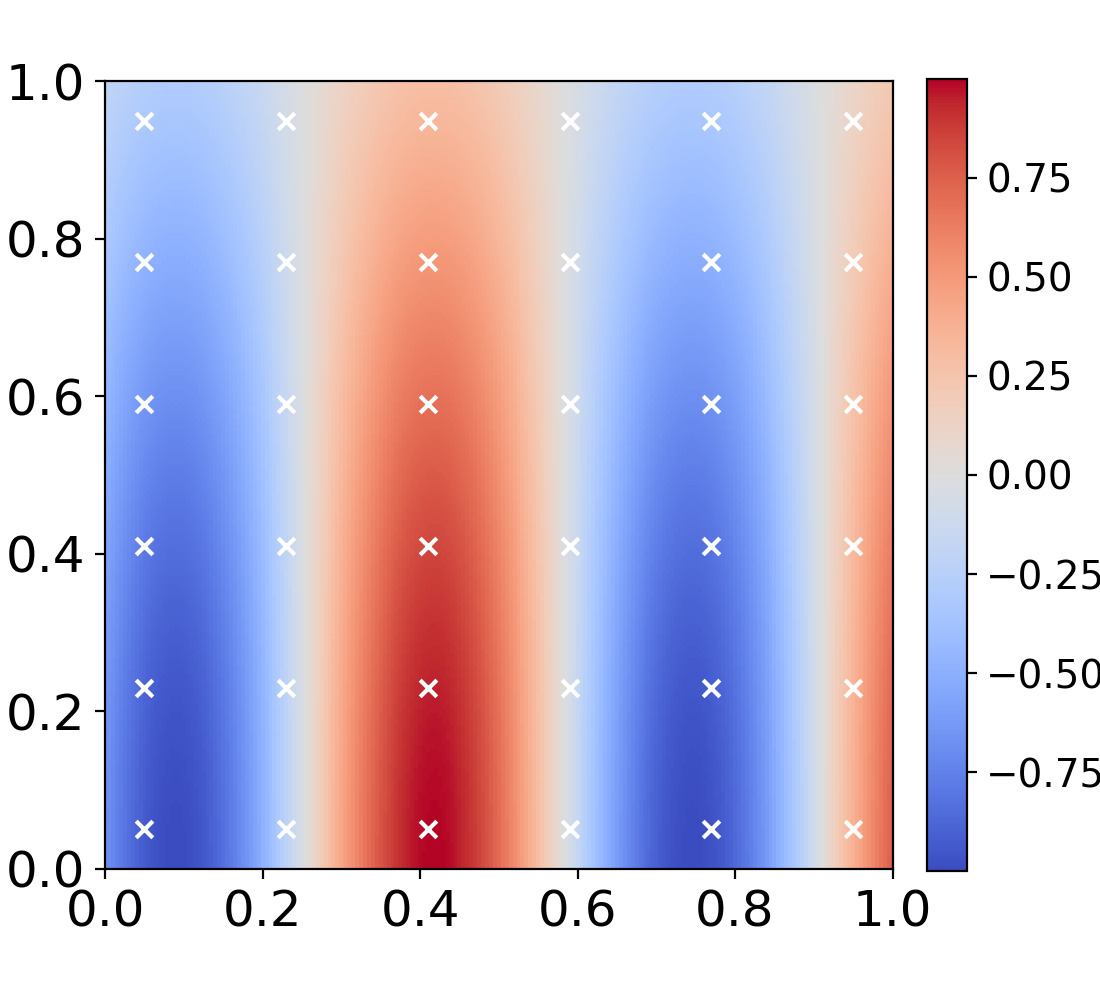}
\includegraphics[width=.22\textwidth]{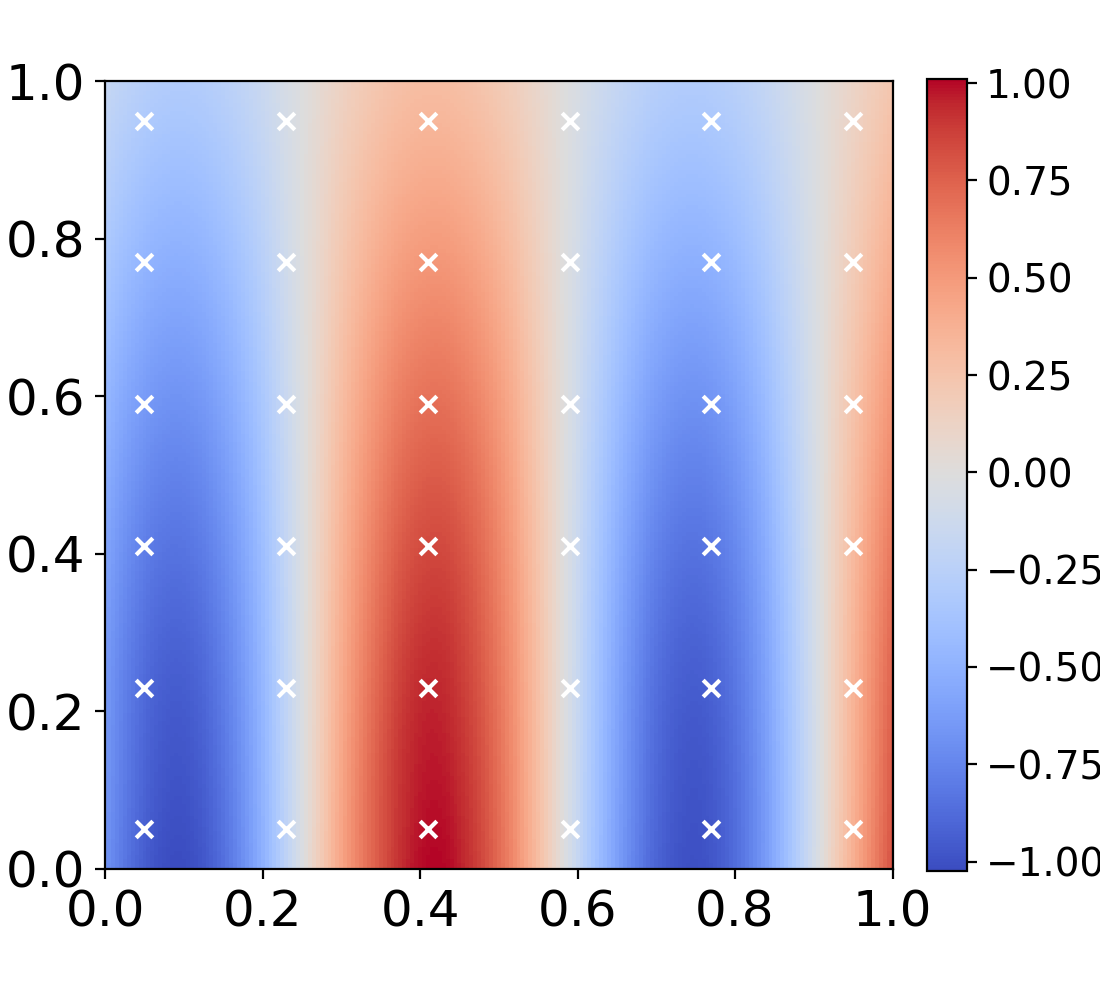}
\includegraphics[width=.22\textwidth]{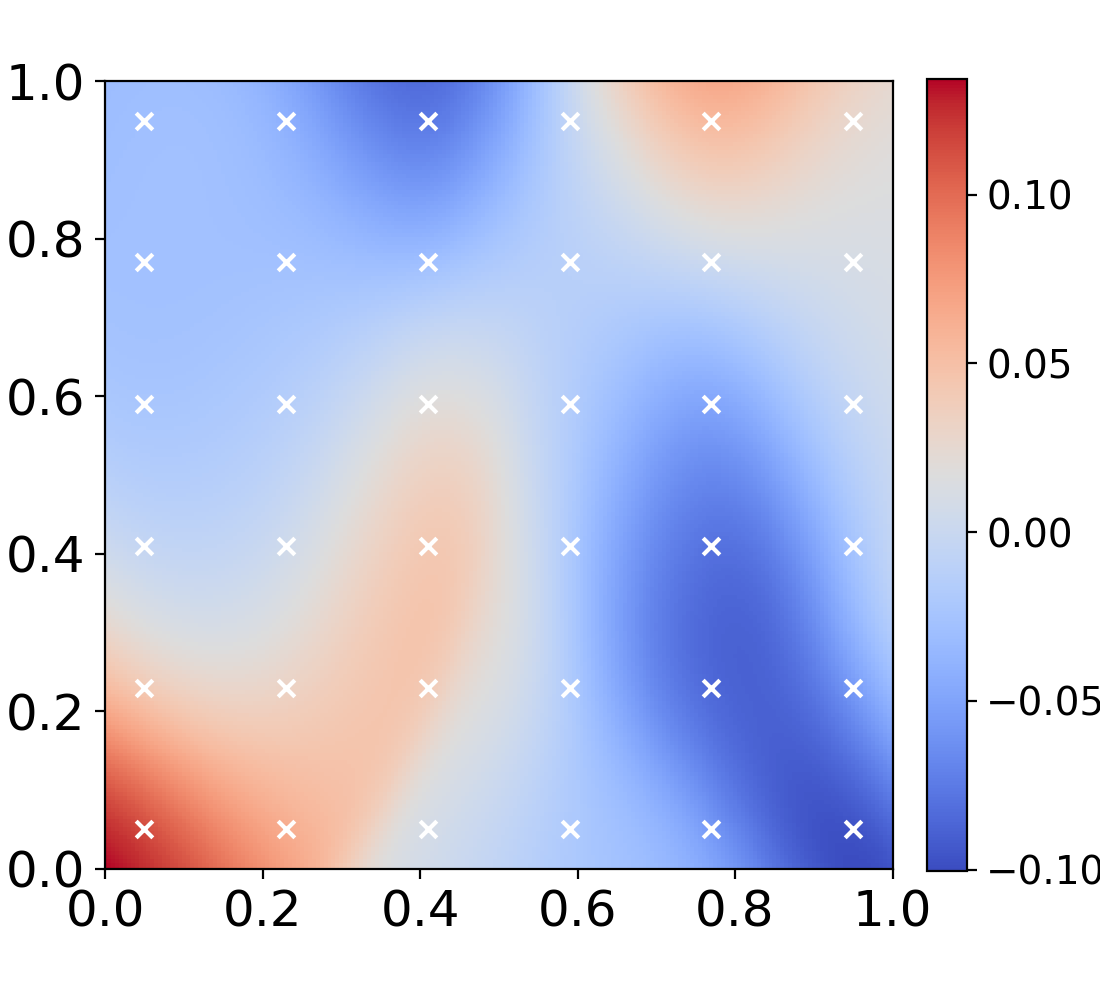}
\includegraphics[width=.22\textwidth]{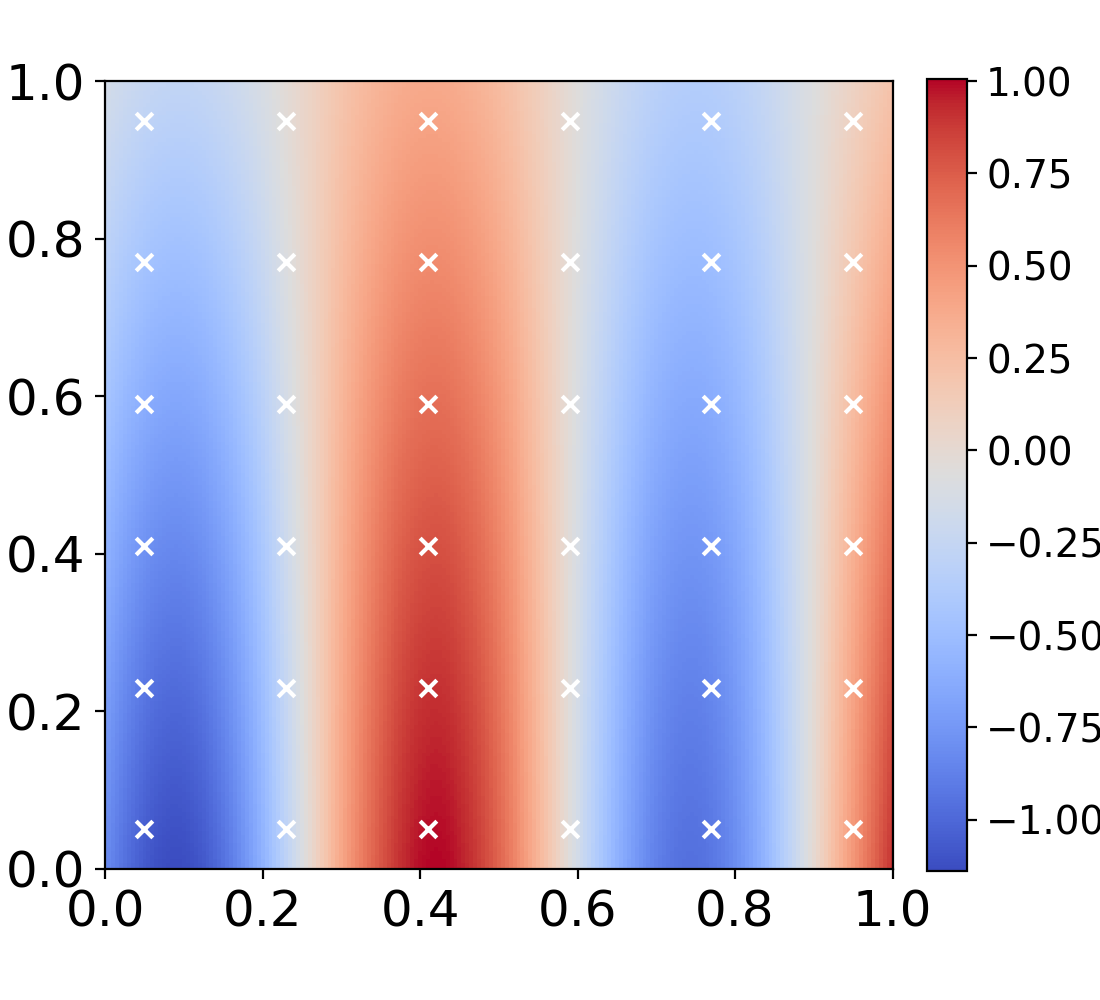}
\end{subfigure}
\caption{Data-driven multifidelity: two-dimensional nonlinear correlation. From left to right: exact high-fidelity solution, multifidelity DeepONet prediction, multifidelity DeepONet nonlinear correlation, and  multifidelity DeepONet linear correlation. The high-fidelity data points are shown in white for clarity.}\label{fig:2d_nonlin_cor}
\end{figure}

We consider a two-dimensional problem  with a nonlinear correlation between the low-fidelity and high-fidelity data: 
\begin{align}
    z_L(u)(x, y) &= \cos(u)\cos(y) + x \\
    z_H(u)(x, y) &= \cos(u)\cos(y)^2  \\
    u &= ax-4
\end{align}
for $x, y \in [0, 1]$ and $a \in [8, 10]$. The training parameters are given in Tab.  \ref{tab:train_params_all} (\ref{sec:training_params}), and the low- and high-fidelity functions are plotted in Fig. \ref{fig:2d_nonlin_exact} (see \ref{2d_nonlinear_app}). The results are given in Fig. \ref{fig:2d_nonlin_results}. Even though the correlation is complex and nonlinear, the composite multifidelity DeepONet improves the predictions by up to an order of magnitude. While the single fidelity method agrees well at locations where training data is provided, overfitting results in large errors in areas where there is no training data. 
In Fig. \ref{fig:2d_nonlin_cor} we show the outputs of the linear and nonlinear DeepONets for two input functions. The correction learned by the nonlinear DeepONet is smaller in magnitude than the output of the linear DeepONet, representing a small nonlinear correction to the linear correlation. 

\subsection{Choice of loss function hyperparameters} \label{sec:choice_lambda}

\textcolor{black}{The loss function introduced in Eq. \ref{eq:loss_data} contains four scaling parameters, $\lambda_i$ for $i=1, 2, 3, 4$. In this section we will discuss how to chose the relative magnitudes of these terms using the jump function in Sec. \ref{1d_jump} and the two-dimensional non-linear correlation in Sec. \ref{2d_nonlinear} as demonstrative examples. We  will fix $\lambda_2$, the coefficient in front of the low-fidelity MSE, and $\lambda_4$, the coefficient that prevents overfitting of the low-fidelity data, and vary $\lambda_1$ and $\lambda_3$. Roughly, $\lambda_1$ denotes the importance of matching the high fidelity training data exactly, and $\lambda_3$ denotes the weight with which the multifidelity DeepONet will learn a bilinear correlation instead of a nonlinear correlation. }

\begin{table}[h!]
    \centering
    \begin{tabular}{c |c|c|c}
    \hline
 $\lambda_1/\lambda_2 $ & $\lambda_3/\lambda_2 $        &  Mean MSE (Eq. \ref{eq:mean_MSE}), HF data & Mean MSE (Eq. \ref{eq:mean_MSE}), LF data  \\ \hline
$0.001$& $1\times 10^{-1}$ &  0.022 & $1.72 \times 10^{-4}$\\
$ 0.01$& $ 1\times 10^{-1}$ & 0.128  & $1.60 \times 10^{-4}$\\
$  0.1$& $ 1\times 10^{-1}$ & 0.014 & $1.53 \times 10^{-4}$\\
$1$& $ 1\times 10^{-1}$ & 0.777  & $1.21 \times 10^{-3}$\\
$  10$& $1\times 10^{-1}$ &  0.832 & $1.52 \times 10^{-3}$\\ 
$  100$& $1\times 10^{-1}$ &  0.927 &  $7.27 \times 10^{-2}$\\
$0.1$& $ 1\times 10^{-4}$ & 0.637&  $1.04 \times 10^{-4}$\\
$0.1$& $ 1\times 10^{-3}$ & 0.751&  $2.18 \times 10^{-4}$\\
$0.1$& $ 1\times 10^{-2}$&  0.722&  $9.87 \times 10^{-5}$\\
$ 0.1$& $ 1\times 10^{0}$ & 0.028&  $2.24 \times 10^{-4}$\\
$ 0.1$ & $ 1\times 10^{1}$ &  0.013 &  $1.21 \times 10^{-4}$\\
        \hline
    \end{tabular}
    \caption{\textcolor{black}{Mean square errors for the trained multifidelity model with varying hyperparameters in the loss function, Eq. \ref{eq:loss_data}, \textcolor{black}{for the one-dimensional jump function from Sec. \ref{1d_jump}}. We fix $\lambda_2 = 1$ and $\lambda_4 = 10^{-4}$. The results in Sec. \ref{1d_jump} take $\lambda_1/\lambda_2 = 0.1$, and $\lambda_3/ \lambda_2 = 1\times 10^{-1}$.}}
    \label{tab:lambda_1d}
\end{table}

\textcolor{black}{For the jump function in Sec. \ref{1d_jump}, we have a very small amount of high fidelity data, but there is an exactly linear correlation between the high fidelity and low fidelity data. Therefore, we expect that it is important for $\lambda_3$ to be relatively large, to emphasize the linear correlation, and $\lambda_1$ to be relatively small, to emphasize learning the low fidelity data well, and therefore a more accurate linear correlation between the low- and high-fidelity data. This is demonstrated in the results in Tab. \ref{tab:lambda_1d}. When  $\lambda_1 / \lambda_2 $ is small, the MSE between the learned model and the high-fidelity data is small, and the error grows with  $\lambda_1 / \lambda_2 $. In contrast, the error decreases as  $\lambda_3 / \lambda_2 $ increases. This can be see\textcolor{black}{n} clearly in the learned linear correlations, given in Appendix \ref{1d_jump_app}.}

\textcolor{black}{For the two-dimensional nonlinear correlation, the trends in Tab. \ref{tab:lambda_2d} are less clear because there is no large linear correlation to learn. When $\lambda_1 / \lambda_2 $ is very small, the high fidelity MSE is larger, due to increased relative values of the other terms in the loss function. A lower weight is placed on learning the high fidelity training set. Our results show that the training is robust to changes in $\lambda_3$, for very nonlinear correlations.}

\begin{table}[h!]
    \centering
    \begin{tabular}{c |c|c|c}
    \hline
 $\lambda_1/\lambda_2 $ & $\lambda_3/\lambda_2 $        &  Mean MSE (Eq. \ref{eq:mean_MSE}), HF data & Mean MSE (Eq. \ref{eq:mean_MSE}), LF data  \\ \hline
$0.01$& $1\times 10^{-3}$ & $1.74\times 10^{-3} $&$ 1.59\times 10^{-4} $ \\
$ 0.1$& $ 1\times 10^{-3}$ &$2.38\times 10^{-4} $ & $ 4.53\times 10^{-5}$ \\
$  1$& $ 1\times 10^{-3}$ & $3.54\times 10^{-5} $&  $1.76\times 10^{-5} $\\
$10$& $ 1\times 10^{-3}$ & $4.15\times 10^{-5} $& $2.18\times 10^{-5} $ \\
$  100$& $1\times 10^{-3}$ & $2.85\times 10^{-5} $&$  1.94\times 10^{-5}$ \\
$1$& $ 1\times 10^{-5}$ & $7.47\times 10^{-5}$&$2.53\times 10^{-5}$ \\
$1$& $ 1\times 10^{-4}$& $3.92\times 10^{-5}$&$1.69\times 10^{-5}$ \\
$ 1$& $ 1\times 10^{-2}$ & $2.98\times 10^{-5}$&$2.14\times 10^{-5}$ \\
$ 1$ & $ 1\times 10^{-1}$ & $2.85\times 10^{-5}$&$1.94\times 10^{-5}$ \\
        \hline
    \end{tabular}
    \caption{Mean square errors for the trained multifidelity model with varying hyperparameters in the loss function, Eq. \ref{eq:loss_data}, \textcolor{black}{for the two-dimensional nonlinear correlation problem from Sec. \ref{2d_nonlinear}}. We fix $\lambda_2 = 1$ and $\lambda_4 = 10^{-4}$. The results in Sec. \ref{2d_nonlinear} take $\lambda_1/\lambda_2 = 1$, and $\lambda_3/ \lambda_2 = 1\times 10^{-3}$.}
    \label{tab:lambda_2d}
\end{table}

\textcolor{black}{In practice, the weights $\lambda_i$ can be picked according to knowledge about the problem if such knowledge is available. If, for example, the low fidelity data is very noisy, $\lambda_4$ can be increased to prevent overfitting to a noisy training set. If there is expected to be a strong linear correlation between the high fidelity and low fidelity training sets, $\lambda_3$ can be increased. If only a very small amount of high fidelity data are available, $\lambda_1$ can be decreased. Many methods have been developed recently to adaptively choose weighting terms in neural network loss functions, and these methods can be applied well to multifidelity DeepONets, as well. For example, while it is outside the scope of this paper, we have found the soft attention mechanism weights from \cite{mcclenny2020self} perform well for multifidelity DeepONets} \textcolor{black}{for limited cases tested. However, the soft attention mechanism weights were designed for PINNs to learn weights that are adaptive in space for each term of the loss function, and may not be as applicable to operator training where the relative areas in which higher weight terms are needed vary depending on the initial condition.}

\subsection{Ice-sheet modeling} \label{sec:ice}
To provide an example of the data-driven multifidelity method applied to complex problems, here we consider ice-sheet modeling. 
Ice-sheet models are an important component of earth system models and are critical for computing projections of sea-level rise. In order to quantify the uncertainty of sea-level projections, the ice-sheet models need to be evaluated a large number of times. \textcolor{black}{The basal friction field, $\beta(x, y)$, is the friction coefficient measured between the ice-sheet and the bedrock under the ice-sheet. It is one of the biggest controls on ice velocity, and cannot be measured directly so it is typically estimated by solving a PDE-constrained optimization problem \cite[e.g.][]{perego2014} to assimilate observation of the surface ice velocity. As a result, the basal friction field is affected by both uncertainties in the observations and in the uncertainties in the model. One research goal is to perform uncertainty quantification in the ice-sheet evolution and mass change as a result \textcolor{black}{of} the uncertainty of $\beta$. This requires the solution of the ice-sheet models for a large number of samples of $\beta$. At present, the high computational cost of these models hinders the ability to perform uncertainty quantification, and DeepONets can be used to overcome this issue by generating inexpensive surrogates of ice-sheet models. In our recent work, \cite{he2023hybrid}, we have trained single fidelity DeepONets to create an efficient hybrid method, which accelerates calculations of the ice-sheet velocity for a given $\beta$ field and allows for fast computations for an ensemble of basal friction fields. However, training a single fidelity DeepONet as in \cite{he2023hybrid} requires a great deal of outputs from existing computational models, which is expensive to generate. In \textcolor{black}{the present} work, we will discuss how multifidelity DeepONets can use a small amount of data generated by more accurate high-fidelity models, coupled with large amounts of data generated by less accurate models, to accurately predict the ice-sheet velocity.} 
 The proposed multifidelity DeepOnet framework can be used to leverage the available hierarchy of ice-sheet models~\cite[e.g.][]{Dukowicz2010,perego2012,robinson2022} of different fidelity and cost, and significantly reduce the cost of generating model data. This hierarchy of models is based on different approximations of the Stokes equation that exploit the shallow nature of ice-sheets (see \ref{sec:ice_models}). In this work we will focus on two ice-sheet models, the low-fidelity (zeroth-order\footnote{Here the order of an approximation refers to the order of the terms, with respect to the aspect-ratio (thickness over horizontal length) of an ice-sheet, that are retained in the approximation}) \emph{Shallow Shelf Approximation} (SSA)~\cite{morland_johnson_1980}, and the higher-fidelity model, \emph{MOno-Layer Higher-Order model} (MOLHO)~\cite{dosSantos2022}. In addition to considering different models, we will also consider the same model at different resolutions. \textcolor{black}{xA more general discussion of applying single fidelity DeepONets to ice-sheet modeling, including the computational models used in this work, is available in \cite{he2023hybrid}.} 

\subsubsection{Data generation}
While it is possible to characterize the probability distribution for $\beta$ using a Bayesian inference approach \cite[e.g.][]{petra2014}, here we adopt a simplified log-normal distribution for $\beta$. We write the basal friction field as $\beta = \exp(\gamma)$, where $\gamma$ is normally distributed as
\begin{equation} \label{distribution}
\gamma \sim \mathcal G\left(\log(\beta_{\text{opt}}), k_l\right),\; \text{ and } \; k_l(\mathbf x_1,\mathbf x_2) = a \exp\left(-\frac{|\mathbf x_1-\mathbf x_2|^2}{2l^2}\right).
\end{equation}
 Here $\beta_{\text{opt}}$ is the nominal value of $\beta$, often obtained by assimilating the observed velocities \cite{perego2014}, $l$ the correlation length and $a$ a scaling factor. \textcolor{black}{The correlation length is associated to the smoothness of the samples: larger correlations lead to smoother fields. In this work we choose the correlation lengths so that there is enough spatial variability in the basal friction field to showcase our method while avoiding too small correlation lengths that would require larger training sets.}. \textcolor{black}{For a discussion of varying the correlation length in the ice-sheet model, see \cite{he2023hybrid}.}
 
\begin{figure}[ht]
\centering
\includegraphics[width=\textwidth]{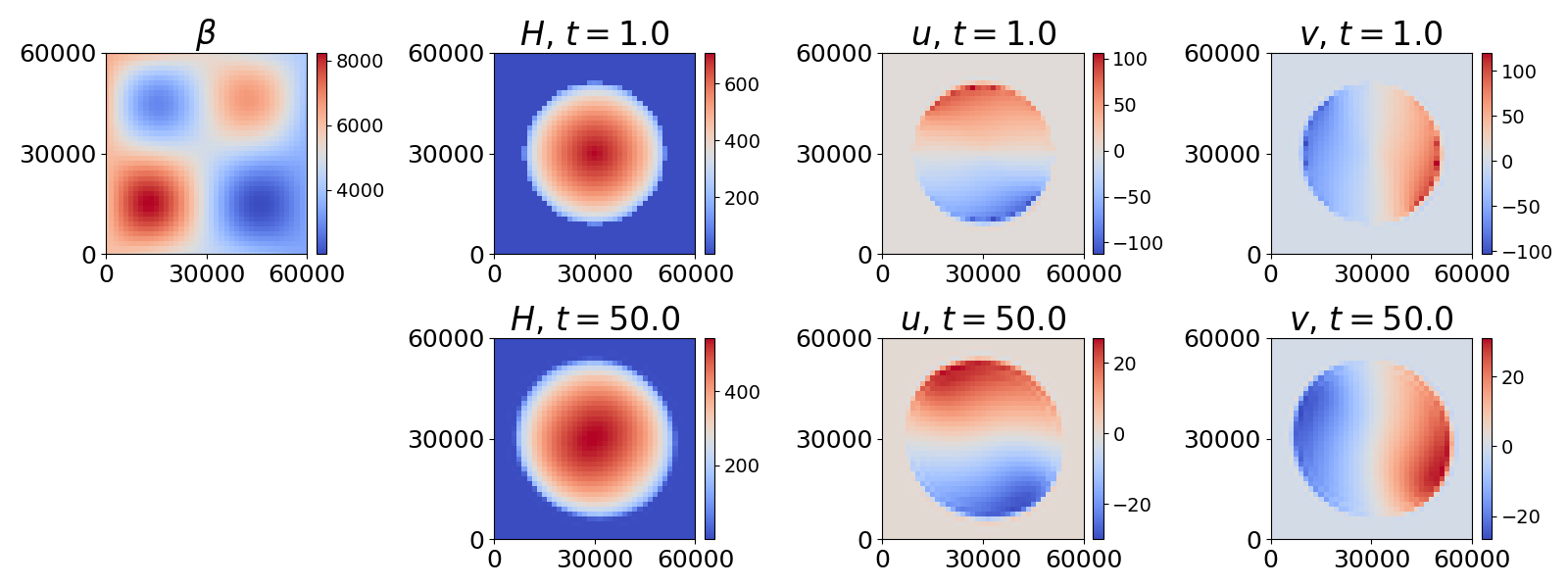}
\caption{Data-driven multifidelity: multiresolution ice-sheet dynamics. Example of the ice-sheet dynamics from an MOLHO simulation with resolution 41x41 at times $t = 1.0$~yr and $t = 50.0$~yr. The basal friction is denoted by $\beta$ (units: [Pa yr\,/\,m]), the ice thickness by $H$ (units: [m]), and the $x-$ and $y-$ depth-averaged velocities by $u$ and $v$ (units: [m\,/\,yr]). The basal friction $\beta$ is constant in time.  } \label{fig:Ice-sheets-exact}
\end{figure}

As a step towards enabling efficient uncertainty quantification, we use a DeepONet to find the depth-averaged velocity $\mathbf {\bar u}$ for an ice-sheet as a function of the basal friction and ice-sheet thickness.
In order to create data for training the DeepONet we sample values for beta according to \eqref{distribution}. For each sample, $\beta_i$, one can solve, with a numerical method, the coupled thickness-velocity problem, and obtain values for the ice thickness $H_i^k(x,y)$ and velocities $\mathbf{\bar u}_i^k(x,y)$ at time $t^k$ that can be used to train the DeepONet (see Fig. \ref{fig:Ice-sheets-exact}). \textcolor{black}{The  Mono-Layer Higher-Order (MOLHO) model and the Shallow Shelf approximation (SSA) are solved with a finite element discretization implemented in FEniCS \cite{alnaes2015fenics}. The model uses continuous piece-wise linear finite elements
for both the thickness and the velocity fields, and we solve the discretized problem with PETSc \cite{koric2023data} SNES nonlinear
solvers. }
 
\subsubsection{Fixed MOLHO Model with Multiresolution}\label{sec:multiresolution}

We first consider the multiresolution case for the so-called Halfar dome~\cite{halfar_1983} deformed under no external forcing.
 We generate two training sets using meshes of different resolution and running the problem forward in time for $50$~yrs for each $\beta$ sampled from \eqref{distribution} with $l=54$~km, $a=1.0$, and $\beta_{\text{opt}} = 3000 \left(1 + \frac12 \sin\left(\frac{\pi x}{30\text{km}}\right) \sin\left(\frac{\pi y}{30\text{km}}\right)\right)$ . The low-fidelity training set uses the MOLHO model on a coarse $P_L = 15\times15$ mesh. The high-fidelity training set uses the same MOLHO model but on a much finer $P_H = 41\times41$ mesh. Due to the size of the mesh, it is time demanding both to generate and to train on the full high-resolution training set. 

We consider multiple cases of the single fidelity training with $N_H = 20$ and $N_H = 50$ high-fidelity 41x41 datasets, and the multifidelity case trained on $N_L = 100$ 15x15 low-fidelity datasets and $N_H = 20$ 41x41 high-fidelity datasets.
The training parameters are given in Tab. \ref{tab:train_params_all} (\ref{sec:training_params}) and the computational cost is given in Tab. \ref{tab:Ice-sheets-comp-cost} in \ref{sec:multiresolution_app}. \textcolor{black}{The hyperparameters are chosen as $\lambda_2 =1$, $\lambda_3 = 1\times 10^{-3}$ because we do not expect a strong linear correlation, and through testing we found that $\lambda_1 =10$ provides the most accurate fit to the high fidelity data.}

\begin{figure}[h!]
\begin{subfigure}{0.95\textwidth}
\centering
\caption{Single fidelity, $N_H=20$}
\includegraphics[width=\textwidth]{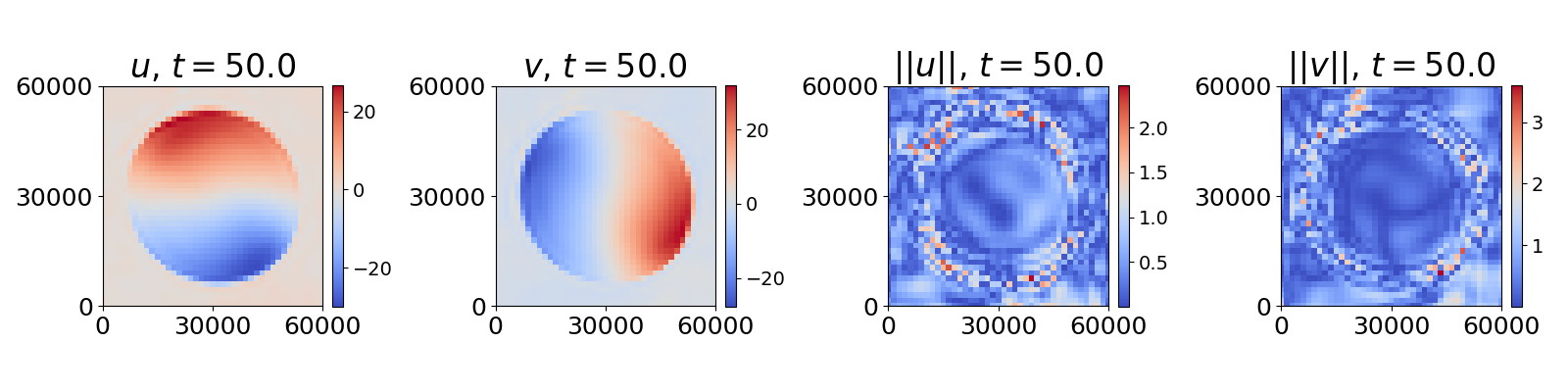}
\end{subfigure}
\begin{subfigure}{0.95\textwidth}
\centering
\caption{Single fidelity, $N_H=50$}
\includegraphics[width=\textwidth]{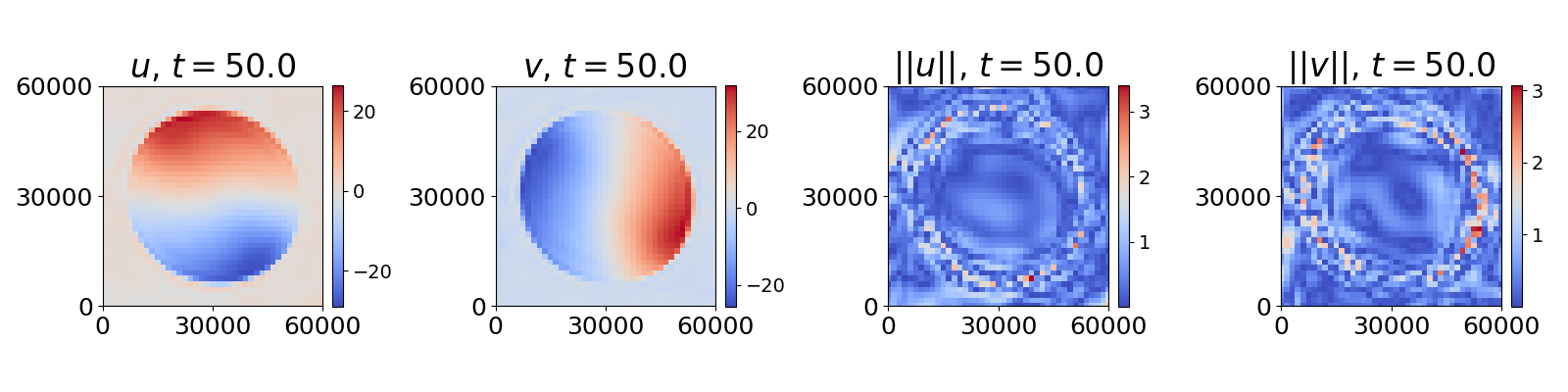}
\end{subfigure}
\begin{subfigure}{0.95\textwidth}
\centering
\caption{Multifidelity, $N_H=10$, $N_L=100$}
\includegraphics[width=\textwidth]{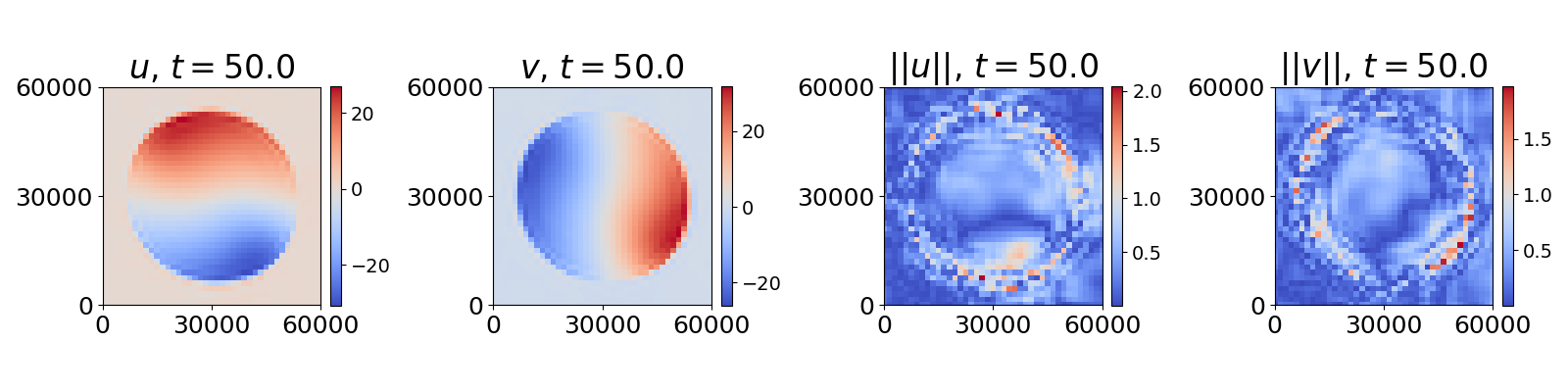}
\end{subfigure}
\caption{Data-driven multifidelity: multiresolution ice-sheet dynamics. Predictions at time $t = 50.0$~yr for the depth-averaged velocity components $u$ and $v$, and pointwise error for $u$ and $v$ are shown from left to right, units: [m / yr]. (a) DeepONet single fidelity prediction with 20 samples in the training set. (b) DeepONet single fidelity prediction with 50 samples in the training set. (b) DeepONet multifidelity prediction with 100 samples in the low-fidelity training set and 20 samples in the high-fidelity training set.  } \label{fig:Ice-sheets-results}. 
\end{figure}

The testing errors are calculated by taking the mean squared error (MSE), defined as
\begin{equation}
    E_{1} = \frac{1}{N_H^T}\sum_{i = 1}^{N_H^T}\frac{1}{P_H} \sum_{j = 1}^{P_H}\left[\bar u(\beta_i, H_i)(\mathbf{x}_j) - \mathcal{F}_l(\beta_i, H_i)(\mathbf{x}_j)- \mathcal{F}_{nl}(\beta_i, H_i)(\mathbf{x}_j)\right]^2 \label{eq:mean_MSE}
\end{equation}
and mean relative $L_2$ error, 
\begin{equation}
    E_2 = \frac{1}{N_H^T}\sum_{i = 1}^{N_H^T}\sqrt{\frac{\sum_{j = 1}^{P_H}\left[\bar u(\beta_i, H_i)(\mathbf{x}_j) - \mathcal{F}_l(\beta_i, H_i)(\mathbf{x}_j)- \mathcal{F}_{nl}(\beta_i, H_i)(\mathbf{x}_j)\right]^2 }{\sum_{j = 1}^{P_H}\left[\bar u(\beta_i, H_i)(\mathbf{x}_j)\right]^2 }}, \label{eq:mean_rel_L2}
\end{equation}
for $N_H^T = 50$ high-resolution 41x41 datasets not used in training, and are shown in Tab. \ref{tab:Ice-sheets-results}. The errors from the single fidelity case with $N_H=50$ datasets match the multifidelity case, although we note that this case requires more than twice the more expensive high-resolution simulations to generate the training data. In comparison, the multifidelity approach achieves accurate results with only $N_H=20$ high-resolution simulations (see Fig. \ref{fig:Ice-sheets-results}). 

\begin{table}[h!]
    \centering
    \begin{tabular}{c|c|c}
    \hline
      Method   &  Mean MSE (Eq. \ref{eq:mean_MSE})& Mean relative $L_2$ error (Eq. \ref{eq:mean_rel_L2}) \\ \hline
    Single fidelity, $N_H=20$ & $1.26$ & 0.0479 \\
    Single fidelity, $N_H=50$ & $1.13$  & 0.0455 \\
    Multifidelity, $N_H=20$, $N_L=100$ & $0.71$ &0.0312\\
        \hline
    \end{tabular}
    \caption{Data-driven multifidelity: multiresolution ice-sheet dynamics. Mean relative $L_2$ errors and mean MSEs for each case, tested over $N_H^T = 50$ testing sets.}
    \label{tab:Ice-sheets-results}
\end{table}

\subsubsection{Multifidelity in Physical Models}\label{sec:multiorder}

We can also consider the case where we have two different models, a low-order and a high-order model, applied to the Humboldt glacier, Greenland. In this case, the low-order model is faster to run, so it is easier to generate a large amount of low-order data to represent the low-fidelity dataset. The high-order method is more time consuming, so we use a small amount of high-order data in the high-fidelity dataset. 

\begin{figure}[ht]
\centering
\includegraphics[width=.22\textwidth]{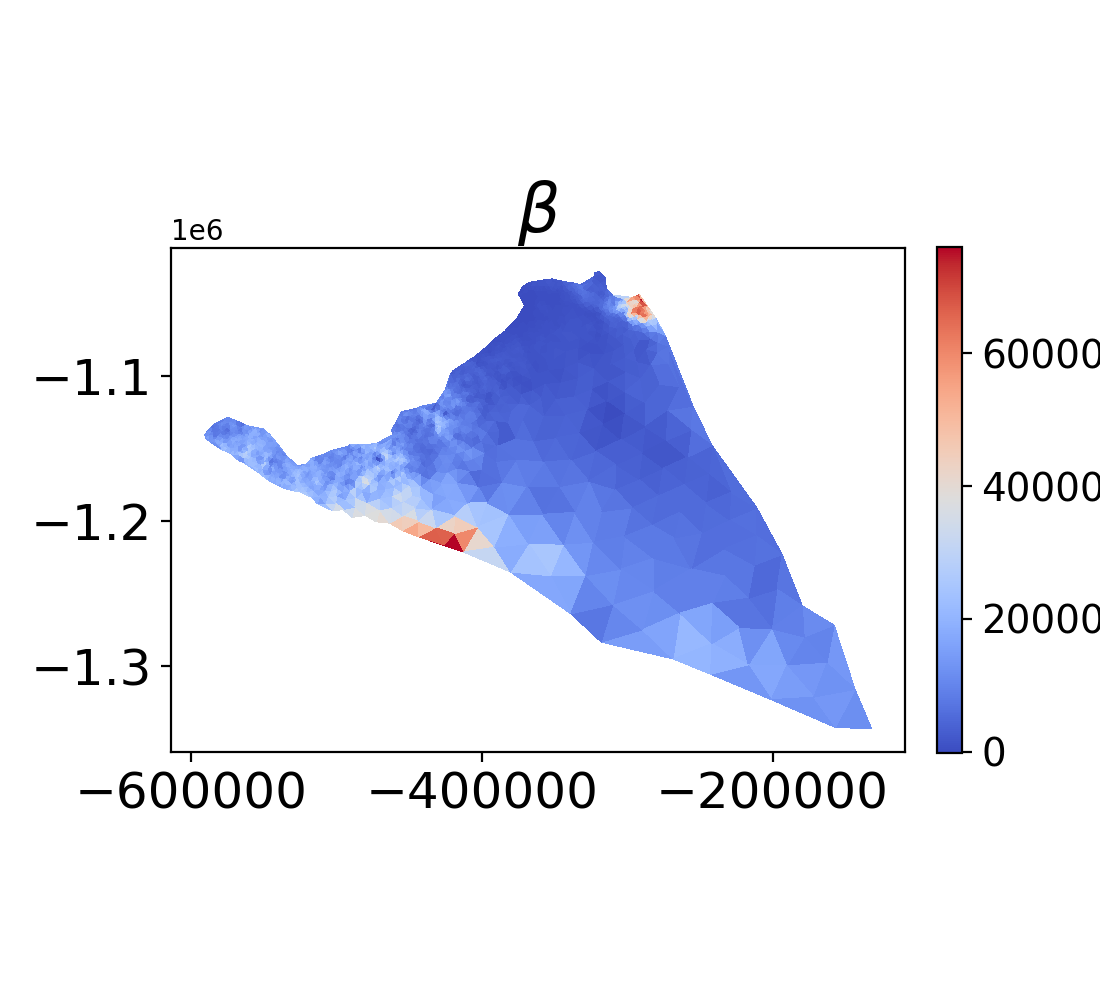}
\includegraphics[width=.22\textwidth]{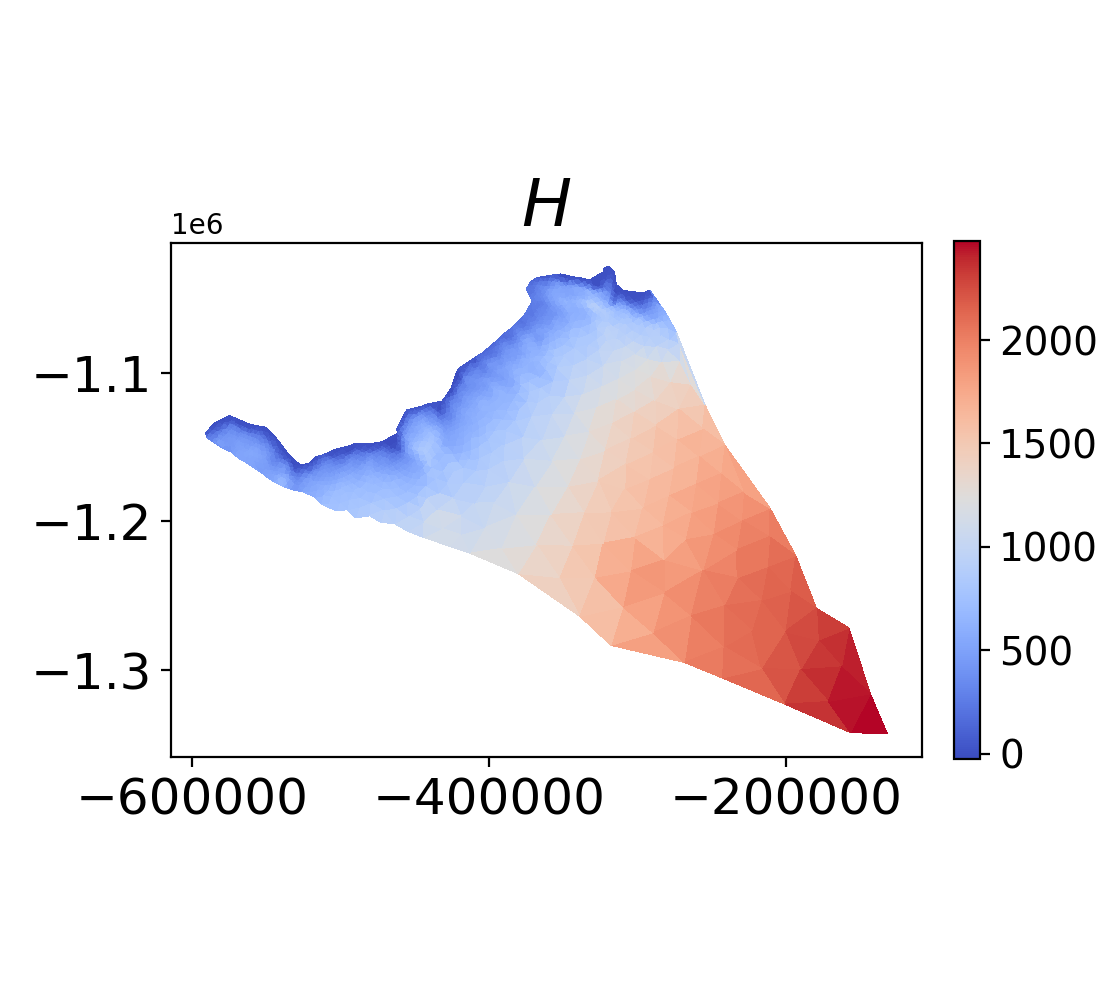}
\includegraphics[width=.22\textwidth]{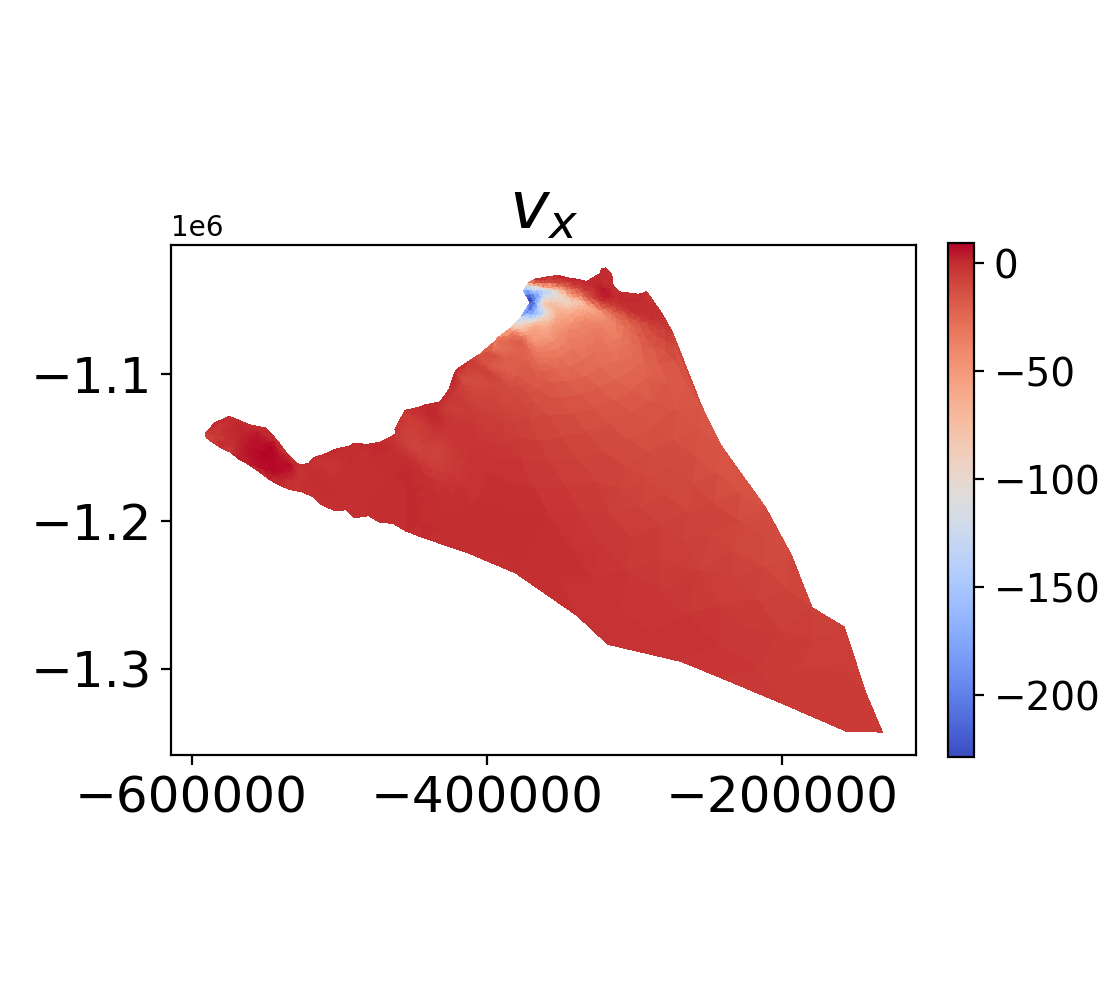}
\includegraphics[width=.22\textwidth]{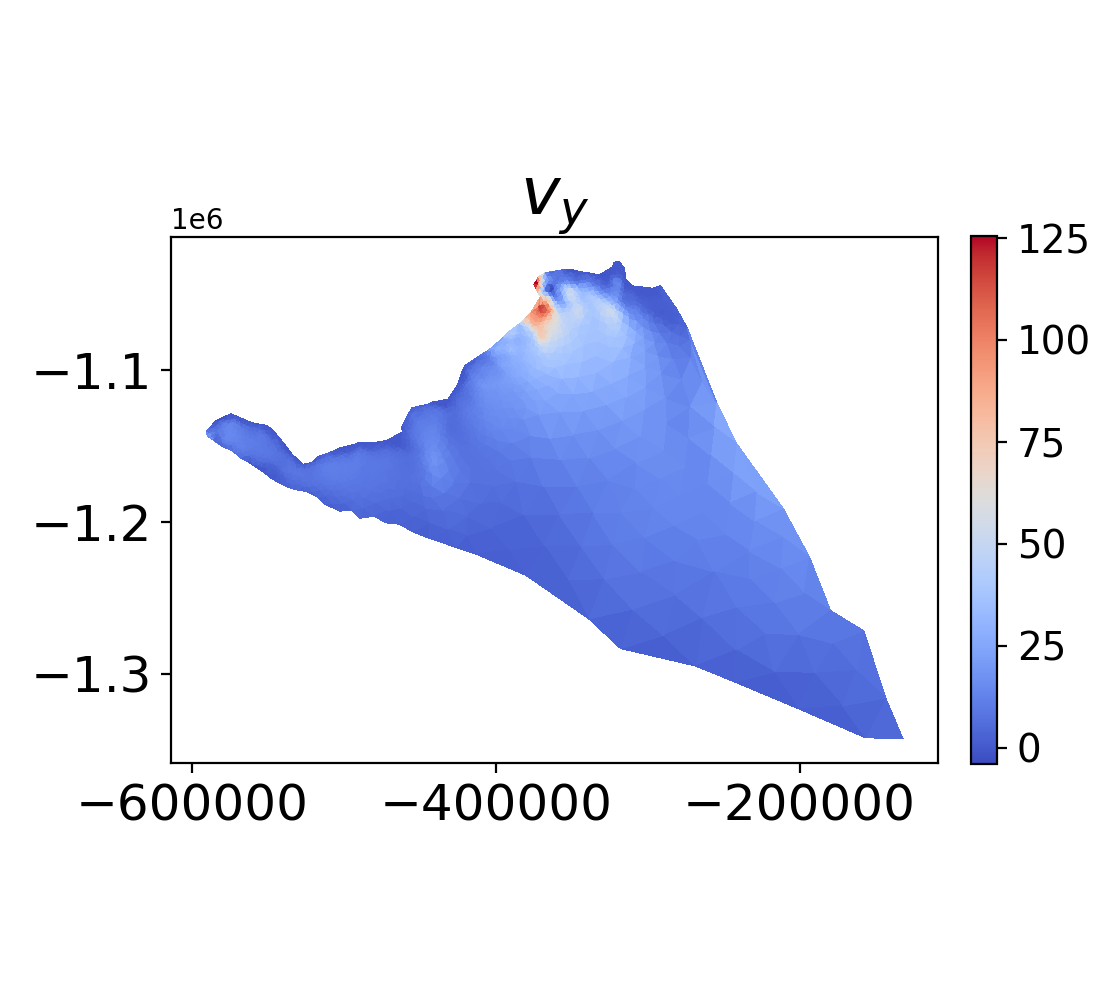}
\caption{Data-driven multifidelity: multiorder ice-sheet dynamics. Example of the ice-sheet dynamics from an MOLHO simulation with 1426 points at times $t = 99.0$~yr. From left to right: the basal friction $\beta$ [Pa yr / m], the ice thickness $H$ [m], and the depth-averaged velocity components $v_x$ and $v_y$ [m / yr].} \label{fig:Ice-sheets-exact-hmbldt}. 
\end{figure}

\begin{figure}[ht]
\begin{subfigure}{\textwidth}
\centering
\caption{Single fidelity}
\includegraphics[width=.22\textwidth]{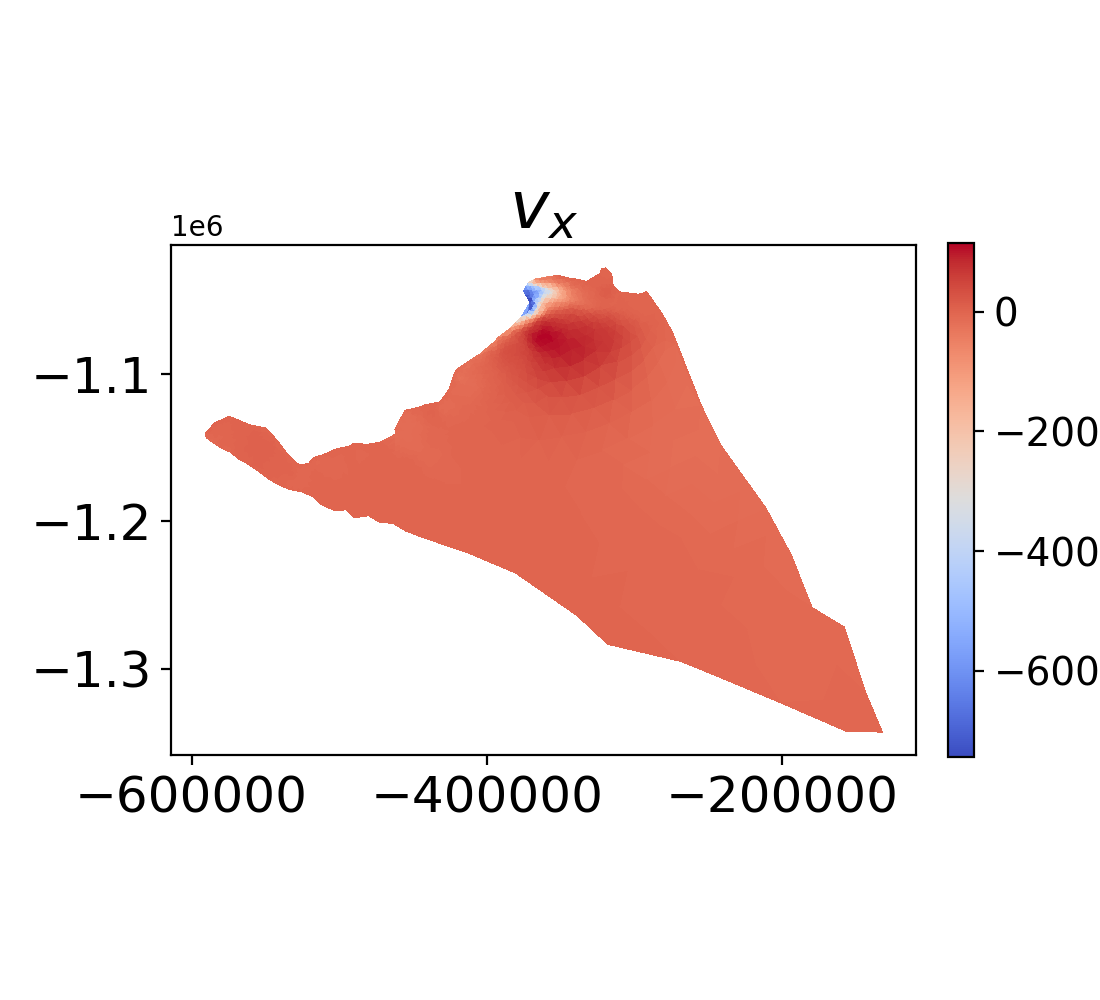}
\includegraphics[width=.22\textwidth]{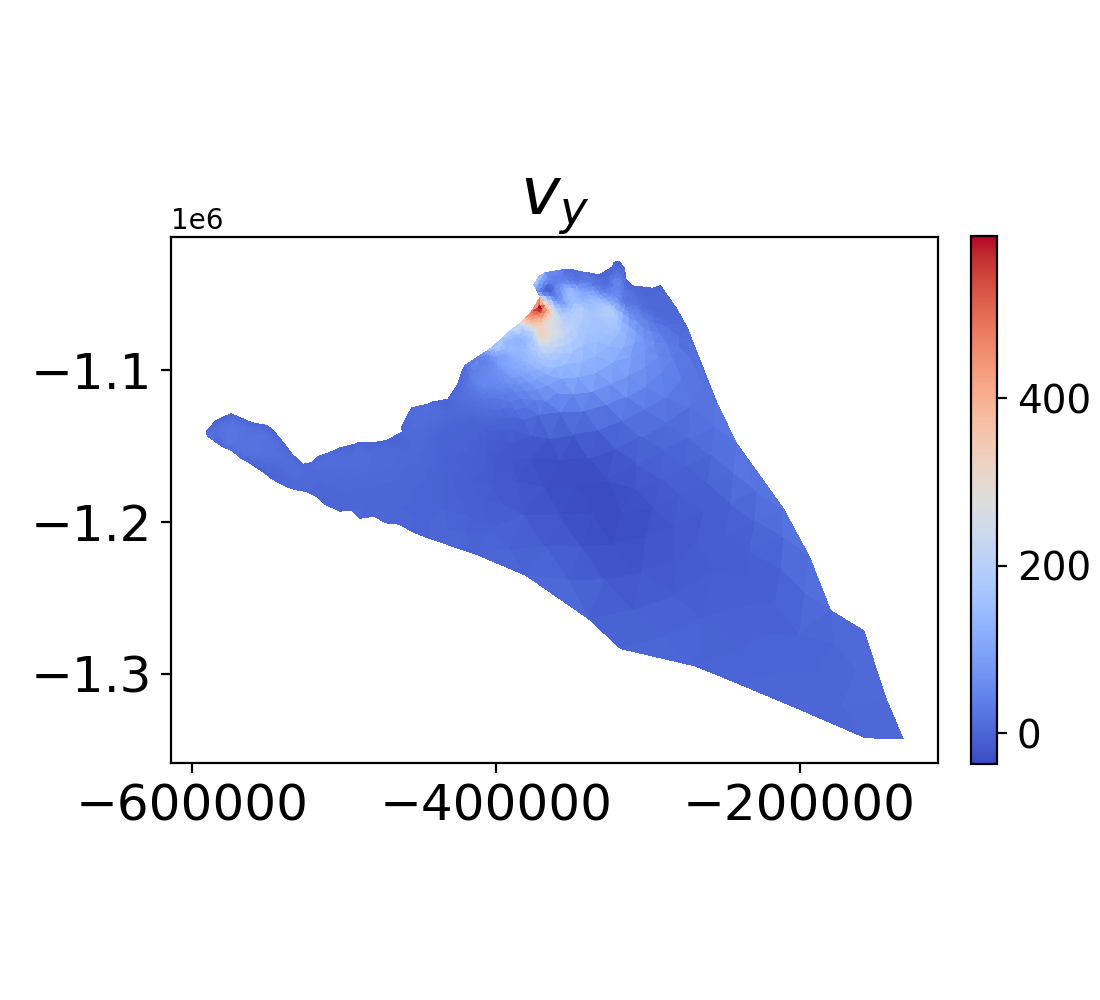}
\includegraphics[width=.22\textwidth]{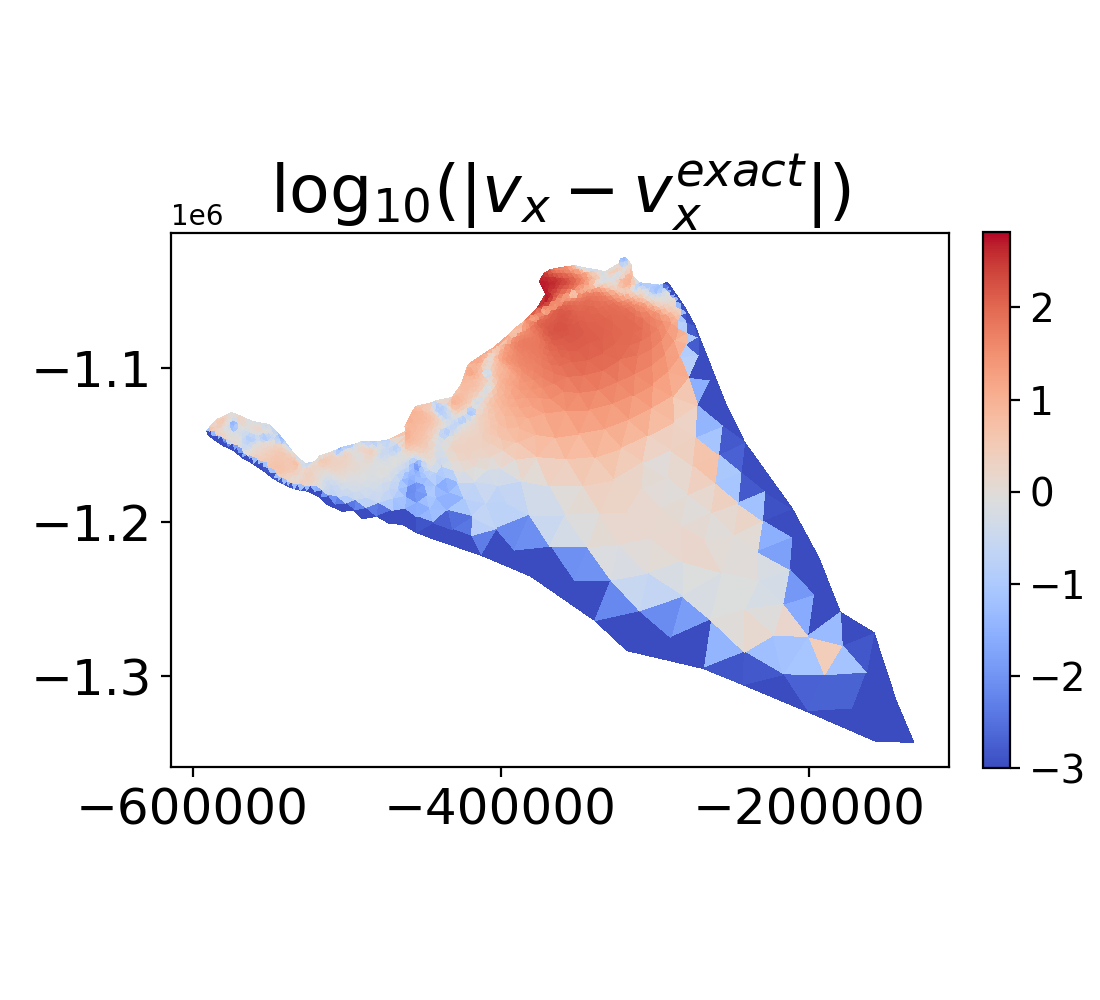}
\includegraphics[width=.22\textwidth]{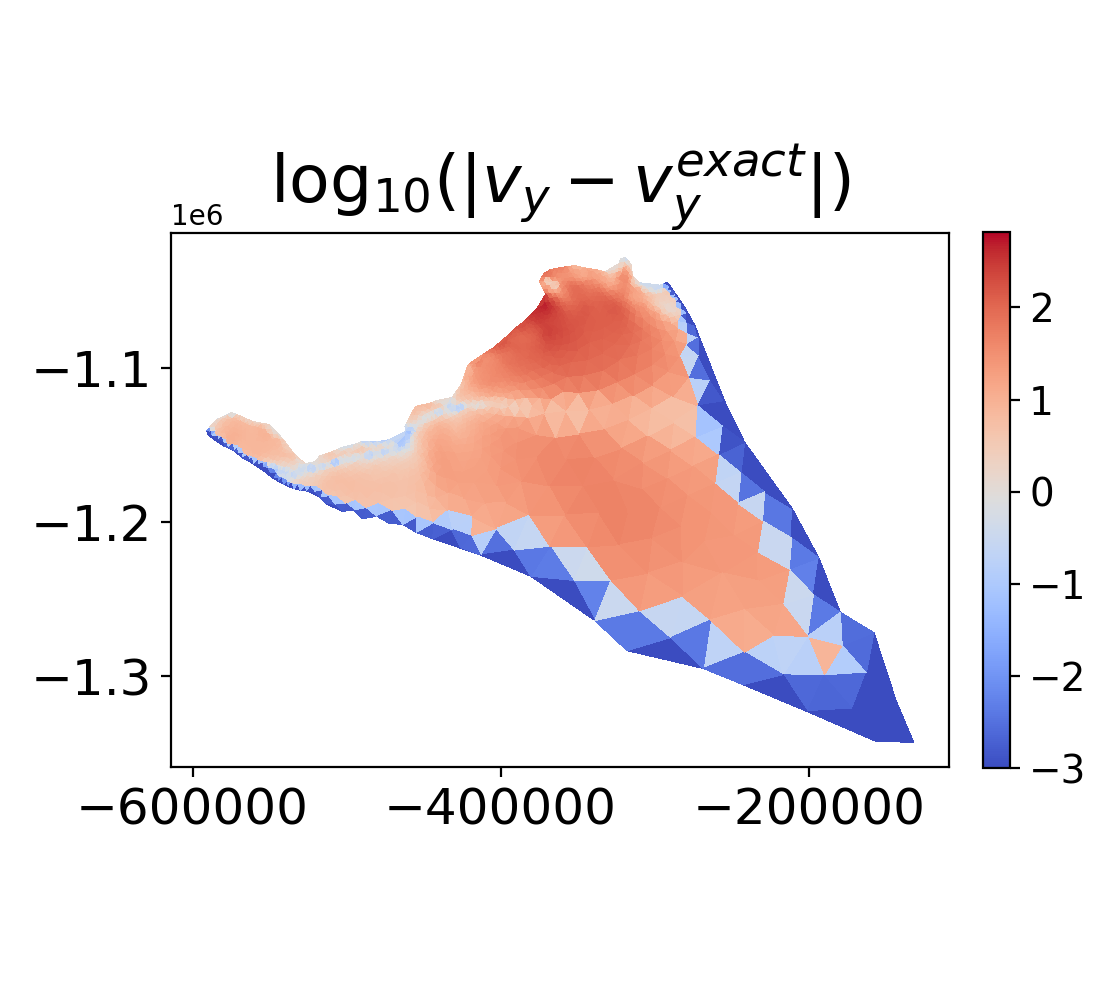}
\end{subfigure}
\begin{subfigure}{\textwidth}
\centering
\caption{Multifidelity}
\includegraphics[width=.22\textwidth]{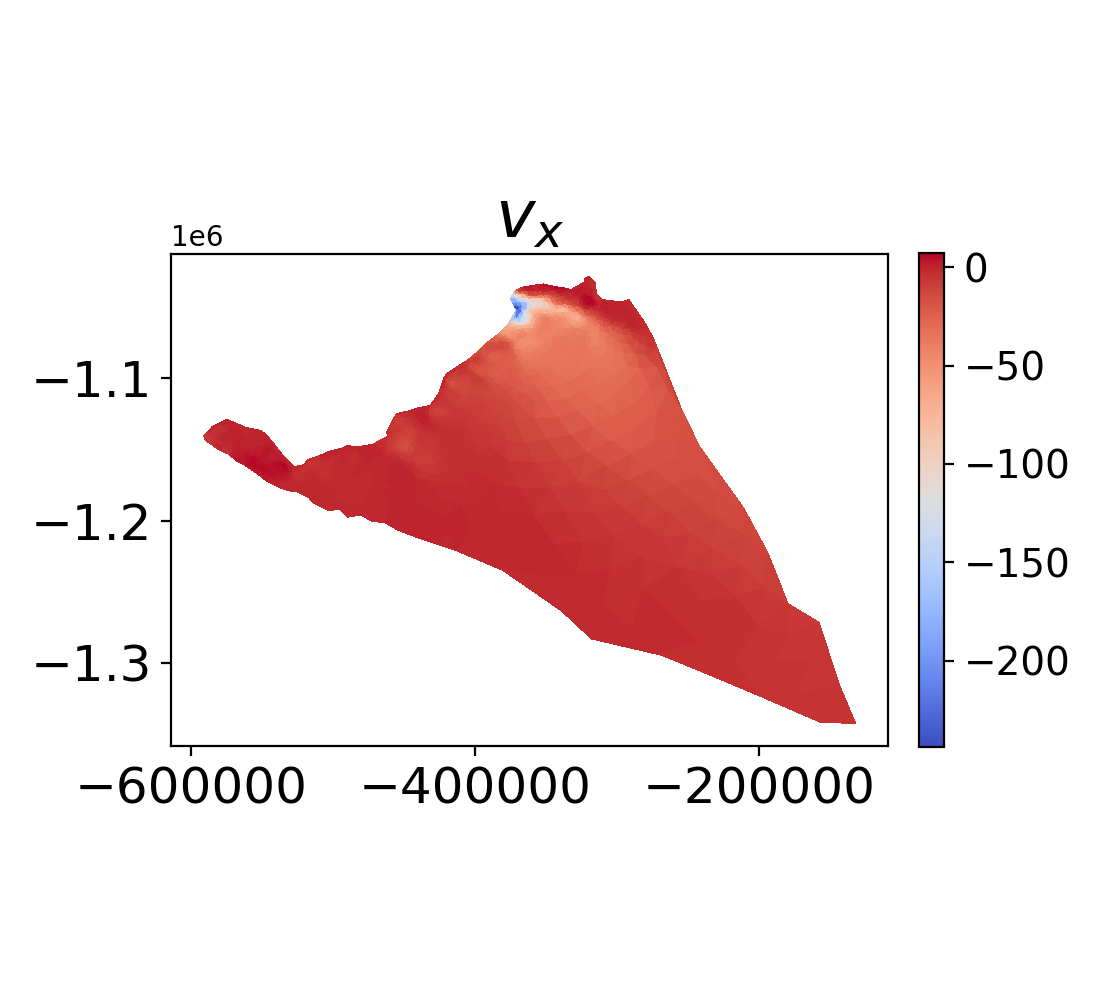}
\includegraphics[width=.22\textwidth]{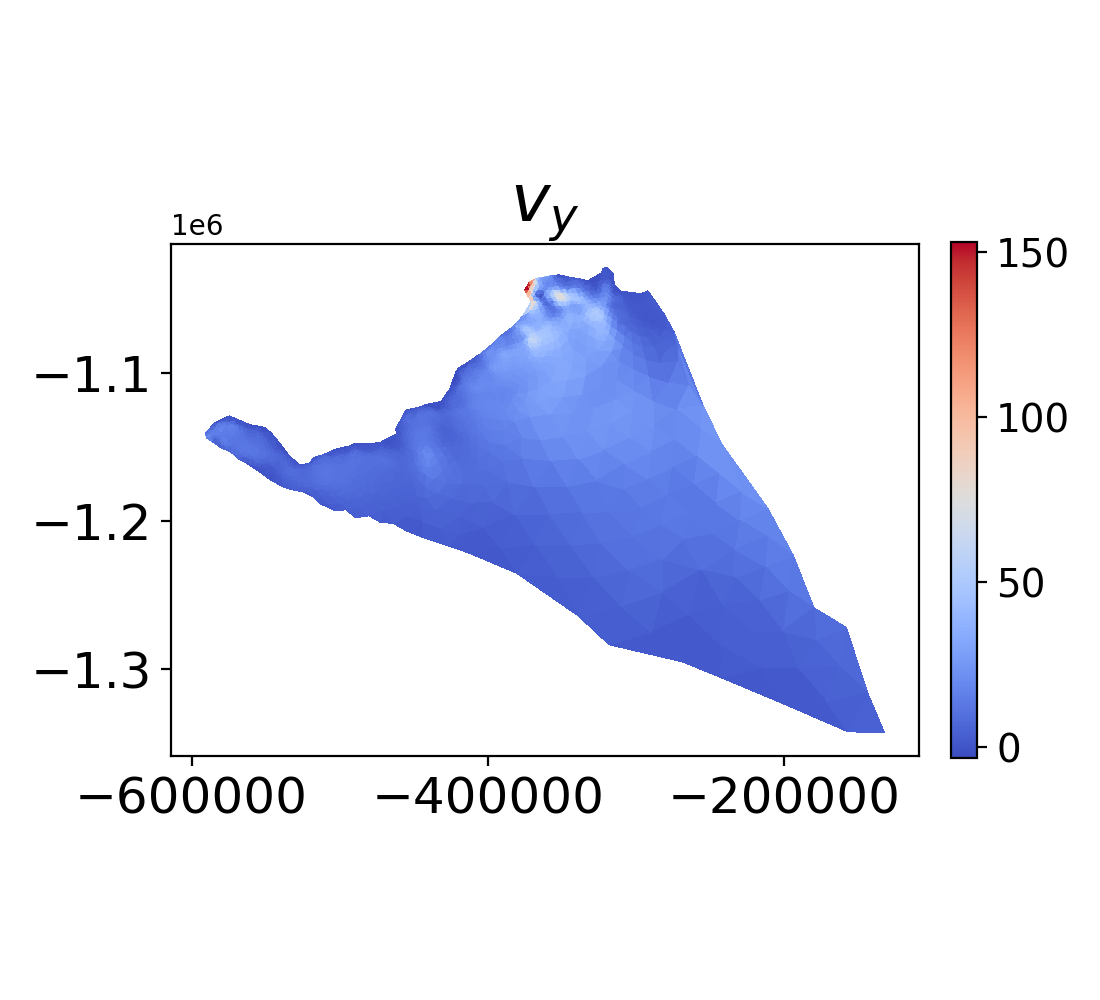}
\includegraphics[width=.22\textwidth]{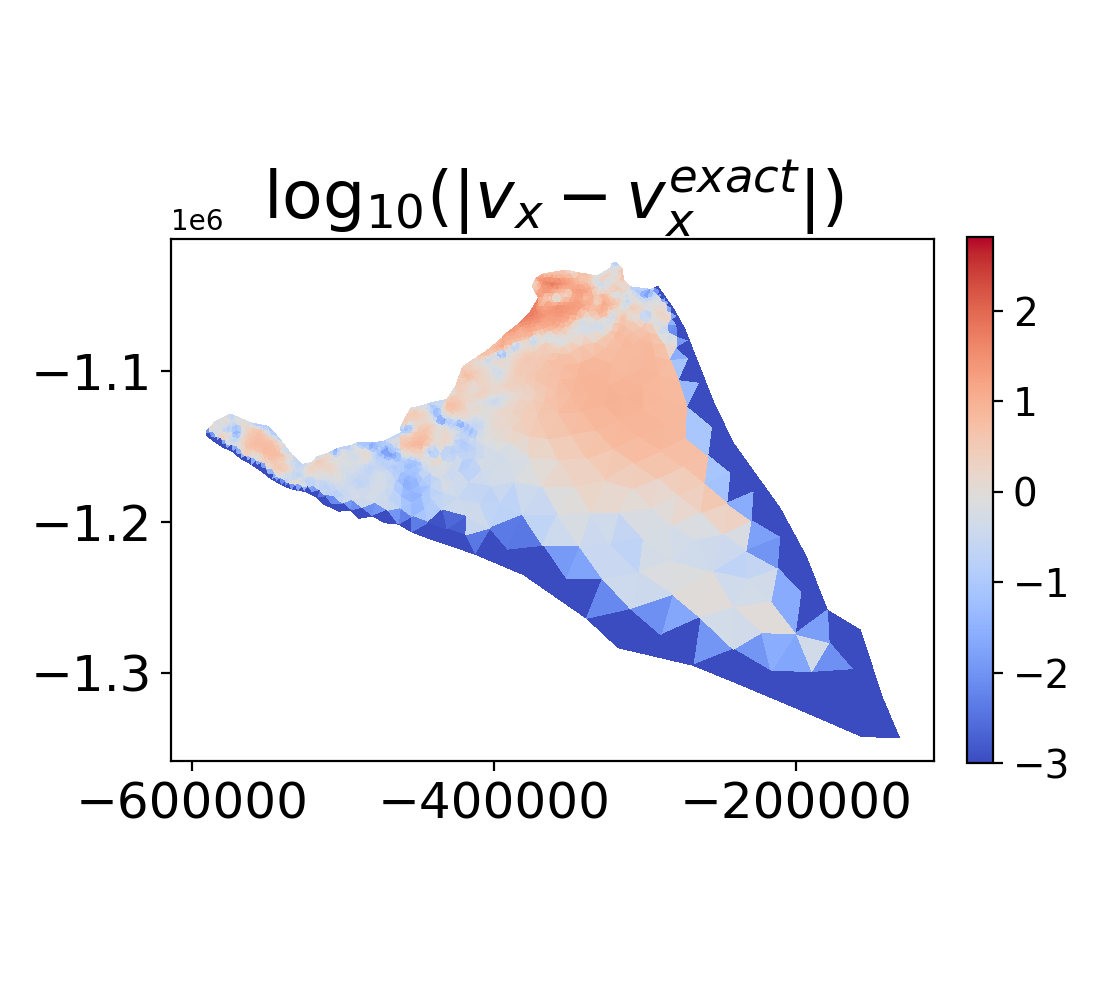}
\includegraphics[width=.22\textwidth]{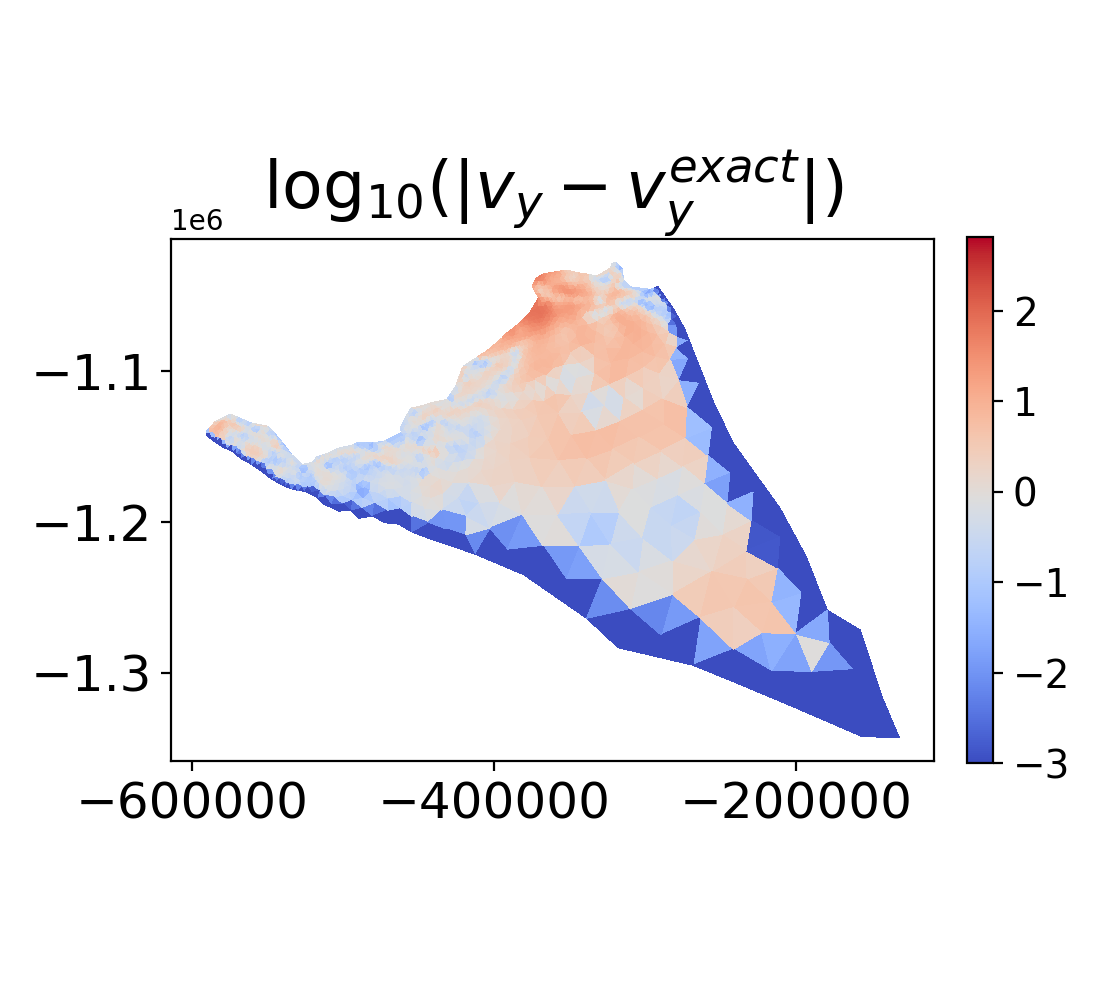}
\end{subfigure}
\caption{Data-driven multifidelity: multiorder ice-sheet dynamics. Output from the test set for the single fidelity (a) and multifidelity (b) training at time $t=99.0$~yr. For both (a) and (b), from left to right we show  $v_x$, $v_y$,  the \textcolor{black}{log of the pointwise error in $v_x$, and the log of the pointwise error in $v_y$, (units: [m / yr]). The colorbars for the logs of the pointwise error are standardized between (a) and (b) for clarity and ease of comparison.}} \label{fig:Ice-sheets-output-hmbldt}. 
\end{figure}

\textcolor{black}{To generate each dataset, we first generate a $\beta$ sample by sampling \ref{distribution} with correlation length $l=80$~km and scaling $a = 0.2$. Then, the ice flow model is run forward in time for $100$~years, using the sampled $\beta$ and a climate forcing $f_H$ (see \ref{thickness}) generated accordingly to the \emph{Representative Concentration Pathway 2.6} (see \cite{Hillebrand2022} for the problem definition and the data used including the nominal basal friction $\beta_{\text{opt}}$).}
We use $N_H = 20$ runs (corresponding to 20 samples of $\beta$) using the MOLHO model as our high-fidelity dataset. The low-fidelity dataset has $N_L = 80$ runs of the SSA model. Both the high- and low-fidelity datasets use the same nonuniform mesh. Example output from the MOLHO simulations is shown in Fig. \ref{fig:Ice-sheets-exact-hmbldt}. The training parameters are given in Tab. \ref{tab:train_params_all} (\ref{sec:training_params}) and the computational cost is given in Tab. \ref{tab:ice_multiord-comp-cost} in \ref{sec:multiorder_app}. \textcolor{black}{The hyperparameters are chosen as $\lambda_2 =1$, $\lambda_3 = 1\times 10^{-3}$ because we do not expect a strong linear correlation, and $\lambda_1 =1$ is chosen to balance the loss between the low- and high-fidelity training sets.} The output from the single fidelity and multifidelity training is shown in Fig. \ref{fig:Ice-sheets-output-hmbldt}. \textcolor{black}{For this test case, the single fidelity training, shown in  Fig. \ref{fig:Ice-sheets-output-hmbldt}a, under- and overestimates the ice sheet velocities, resulting in errors on the same order of magnitude as the velocity values. The multifidelity errors shown in the right two panels of Fig. \ref{fig:Ice-sheets-output-hmbldt}b, are an order of magnitude smaller.} The single fidelity method has a mean relative $L_2$ error of 1.1005, while the multifidelity method has a mean relative $L_2$ error of 0.3676 across five simulations in the high-fidelity test set. 

\section{Physics-informed multifidelity DeepONets}\label{sec:PI}

In this section, we consider the application of the multifidelity method to cases where we have low-fidelity data and knowledge of the physics of the system. In particular, the high-fidelity network is trained without paired input-output observations. In practice, it is often inexpensive to generate low-fidelity solutions to a PDE with a numerical method, but the numerical method may be too expensive to run at high resolution or may be low order and miss important features of the full solution. Here, we show we can use the low-fidelity numerical data as the input to the low-fidelity network, and then correct the output using physics-informed DeepONets to result in more accurate outputs. We consider first a one-dimensional case with a nonlinear correlation between the low-fidelity data and the ODE solution. We then consider the viscous Burgers equation, where the low-fidelity data is provided by a numerical solver both with and without noise added. \textcolor{black}{Single-fidelity physics-informed DeepONets can be difficult to train, and the viscous Burgers equation is an equation that has been shown to present particular challenges due to the sharp gradients that form at later times, similar to the sharp gradients that form in the ice sheet velocity at the edges of the ice sheets in Sec. \ref{sec:multiresolution}.}

\subsection{One-dimensional, nonlinear correlation}\label{sec:1d_pi}
\begin{figure}[ht]
\centering
\begin{subfigure}{0.45\textwidth}
\centering
\caption{$a =0.79605263$}
\includegraphics[width=\textwidth]{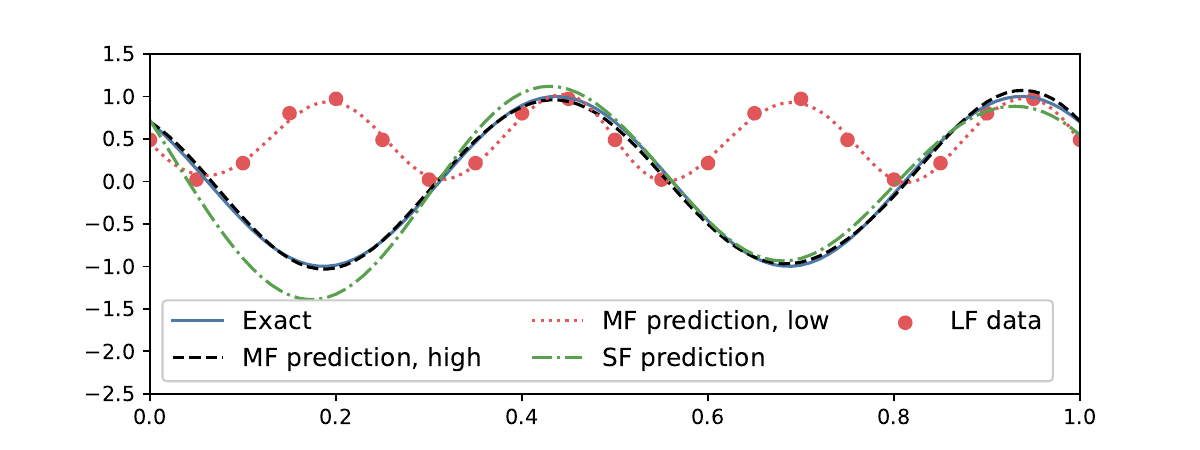}
\end{subfigure}
\begin{subfigure}{0.45\textwidth}
\centering
\caption{$a =4.4526315$}
\includegraphics[width=\textwidth]{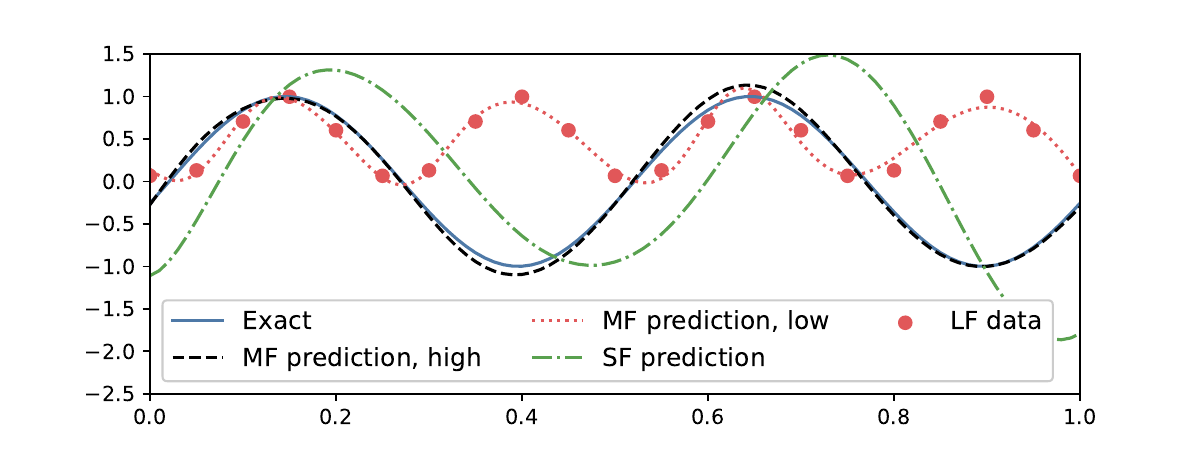}
\end{subfigure}
\begin{subfigure}{.3\textwidth}
\centering
\caption{SF error}
\includegraphics[width=\textwidth]{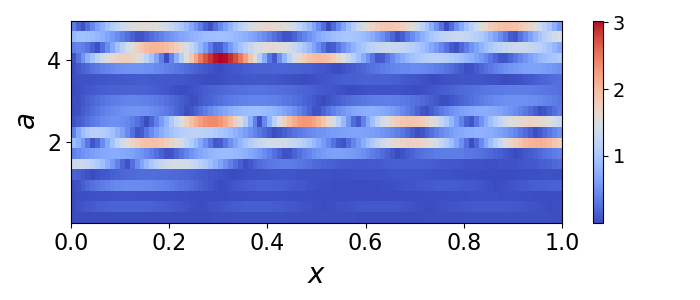}
\end{subfigure}
\begin{subfigure}{.3\textwidth}
\centering
\caption{MF error, high-fidelity data}
\includegraphics[width=\textwidth]{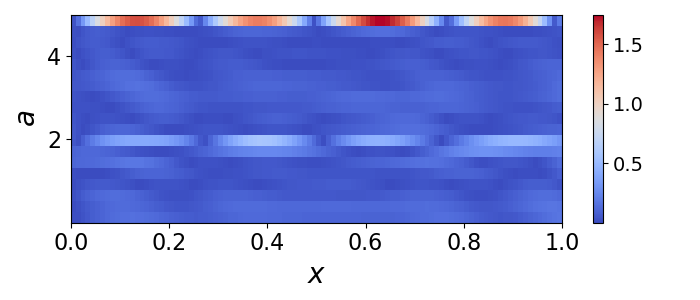}
\end{subfigure}
\begin{subfigure}{.3\textwidth}
\centering
\caption{MF error, low-fidelity data}
\includegraphics[width=\textwidth]{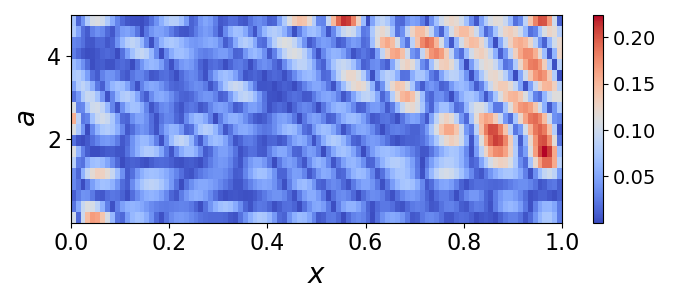}
\end{subfigure}
\caption{Physics-informed multifidelity: nonlinear correlation. (a) Results of the single fidelity and multifidelity predictions of the high- and low-fidelity data with $a =0.79605263$. (b) Results of the single fidelity and multifidelity predictions of the high- and low-fidelity data with $a =4.4526315$. (c) Single-fidelity error as a function of $a$ and $x$, (d) multifidelity high-fidelity prediction error as a function of $a$ and $x$, and (e) multifidelity low-fidelity prediction error as a function of $a$ and $x$.} \label{fig:1d_phys}
\end{figure}

We consider the case given by: 
\begin{align}
  &  y_L(a)(x) = \cos(4 \pi x + a)^2\\
  & \frac{\partial}{\partial x} y_H(a)(x) = -4 \pi\sin(4 \pi x + a)  \\
  & y_H(a)(0) = \cos(a)  
\end{align}
for $x \in [0, 1]$ and $a \in [0, 5]$. We note that the exact solution to the ODE is  $y_H(a)(x) = \cos(4 \pi x + a)$. Parameters are given in Tab. \ref{tab:train_params_all_PI} (see \ref{sec:training_params} and \ref{sec:1d_pi_app}) and results in Fig. \ref{fig:1d_phys}. The testing set consists of 20 input values of $a \in [.0125, 4.975]$ and 101 values of $x \in [0, 1]$. On the test set, the mean MSE of the single fidelity prediction is 0.61648 and the mean MSE of the high-fidelity MF prediction is 0.16232. 

\subsection{Viscous Burgers Equation}\label{sec:burgers}
We consider the viscous one-dimensional Burgers equation with periodic boundary conditions:
\begin{align}
  \textcolor{black}{  \frac{\partial s}{\partial t} + s \frac{\partial s}{\partial x} - \nu \frac{\partial ^2s}{\partial x^2}} &=    0, \; (x, t) \in (0, 1) \times (0, 1] \\
  s(x, 0) &= u(x), \; x\in(0, 1), \\
  s(0, t) &= s(1, t), \; t\in(0, 1), \\
  \frac{\partial s}{\partial x} &= \frac{\partial s}{\partial x}(1, t), \; t\in(0, 1)
\end{align}
where $\nu$ is the viscosity. The initial condition, $u(x)$, is generated from a Gaussian random field (GRF). To generate the training data we follow the procedure in \cite{wang2021learning, wang2021improved, predictiveintelligencelab} and generate $N = 1500$ samples from a Gaussian random field $\sim \mathcal{N}(0, 25^2(-\Delta + 5^2I)^{-4})$. For the physics-informed cases, the initial condition is sampled at $P_{IC} = 100$ uniformly spaced locations on $x = [0, 1]$. The boundary conditions are randomly sampled at $P_{BC} = 100$ locations on $(x, t) = (0, t)$ and $(x, t) = (1, t)$. The residual is evaluated on $P_p = 2,500$ randomly sampled collocation points from the interior of the domain.  The training is completed with $N_H= 1000$ samples, and the results are tested on the remaining 500 samples. We consider three cases: data-only, using the low-fidelity training set, physics-only, using the physics-informed training set, and multifidelity, which combines the low-fidelity training set and physics-informed training set. 

\textcolor{black}{The low-fidelity data for the multifidelity and data-only training are generated by solving Burgers equation with the Chebfun package \cite{driscoll2014chebfun} with a spectral Fourier discretization and a fourth-order stiff time-stepping scheme (ETDRK4) \cite{cox2002exponential} with a timestep $\Delta t = 5\times10^{-3}$ for $\nu = 10^{-2}$ and $\nu = 10^{-3}$ and $\Delta t = 10^{-4}$ for $\nu = 10^{-4}$.} The initial condition is sampled on $M_L = 21$ uniformly spaced locations on $x = [0, 1]$, and snapshots of the solution are saved every $\delta t = 0.05$, to give data on a $21 \times 21$ grid. We take the number of low-fidelity training sets to be either  $N_L = 200$ or $N_L = 1000.$ We also consider a second low-fidelity dataset, generated by adding Gaussian white noise generated by a normal distribution with variance $\sigma^2 = 4\times 10^{-4}$ and mean $0$ to each point in the $21 \times 21$ output. This is referred to as training ``with noise''. \textcolor{black}{To calculate the errors in each case, we generate a high-fidelity test dataset with the same numerical scheme as the low-fidelity data, but with time step $\Delta t = 10^{-4}$,  the initial condition  sampled at $P=101$ locations, and the time snapshots are taken every $\delta t = 0.01$, to give data on a $101 \times 101$ grid. We calculate the relative mean $L_2$ errors between the output of the high-fidelity numerical solver and the DeepONet models on the same grid.}

The multifidelity loss function is given by Eq. \ref{eq:loss_physics}. 
For the single fidelity data-only case, $\lambda_{1}=\lambda_{3}=\lambda_{4} = \lambda_{5} = \lambda_{6} = 0.$ In the single fidelity physics-only case,  $\lambda_{2} = \lambda_{3} = \lambda_{4} = 0.$ Across all cases, the non-zero weights are kept fixed: 
$\lambda_{1} = 10$,  $\lambda_{2} = 1$, $\lambda_{3} = 10^{-6}$,  $\lambda_{4} = 10^{-6}$, $\lambda_{5}=20$, and $\lambda_{6}=1$. The other hyperparameters are given in Tab. \ref{tab:train_params_all_PI} (see \ref{sec:training_params} and \ref{sec:burgers_app}). In the multifidelity training, the low-fidelity output is sampled on an $11 \times 11$ mesh to use as input to the physics-informed nonlinear branch net. 

\begin{figure}[h!]
\begin{subfigure}{\textwidth}
\centering
\caption{$\nu = 10^{-3}$}
\includegraphics[width=\textwidth]{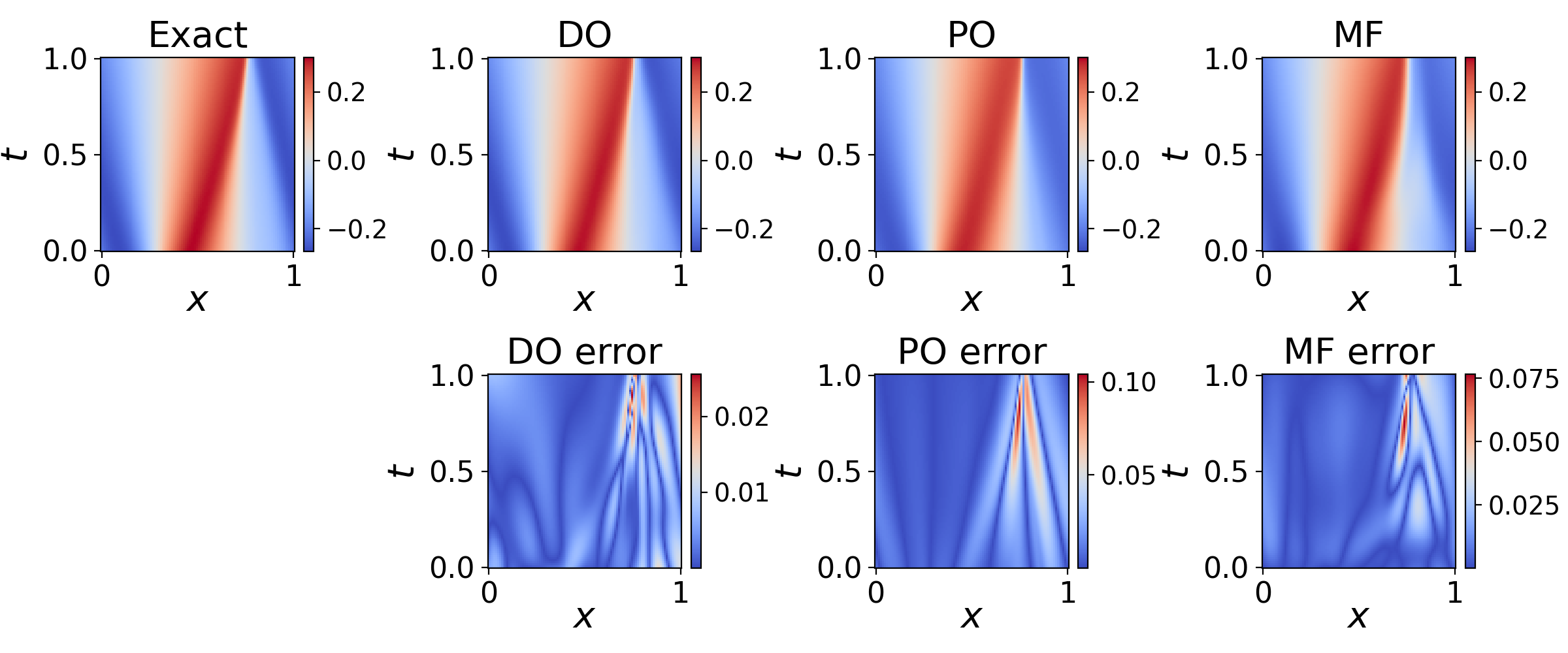}\\
\includegraphics[width=.8\textwidth]{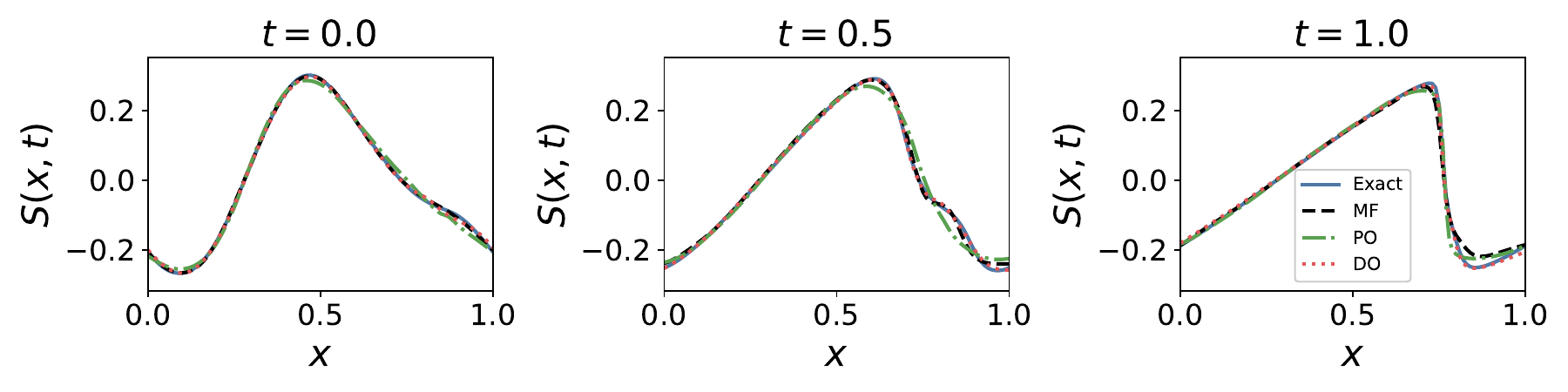}
\end{subfigure}
\begin{subfigure}{\textwidth}
\centering
\caption{$\nu = 10^{-4}$}
\includegraphics[width=\textwidth]{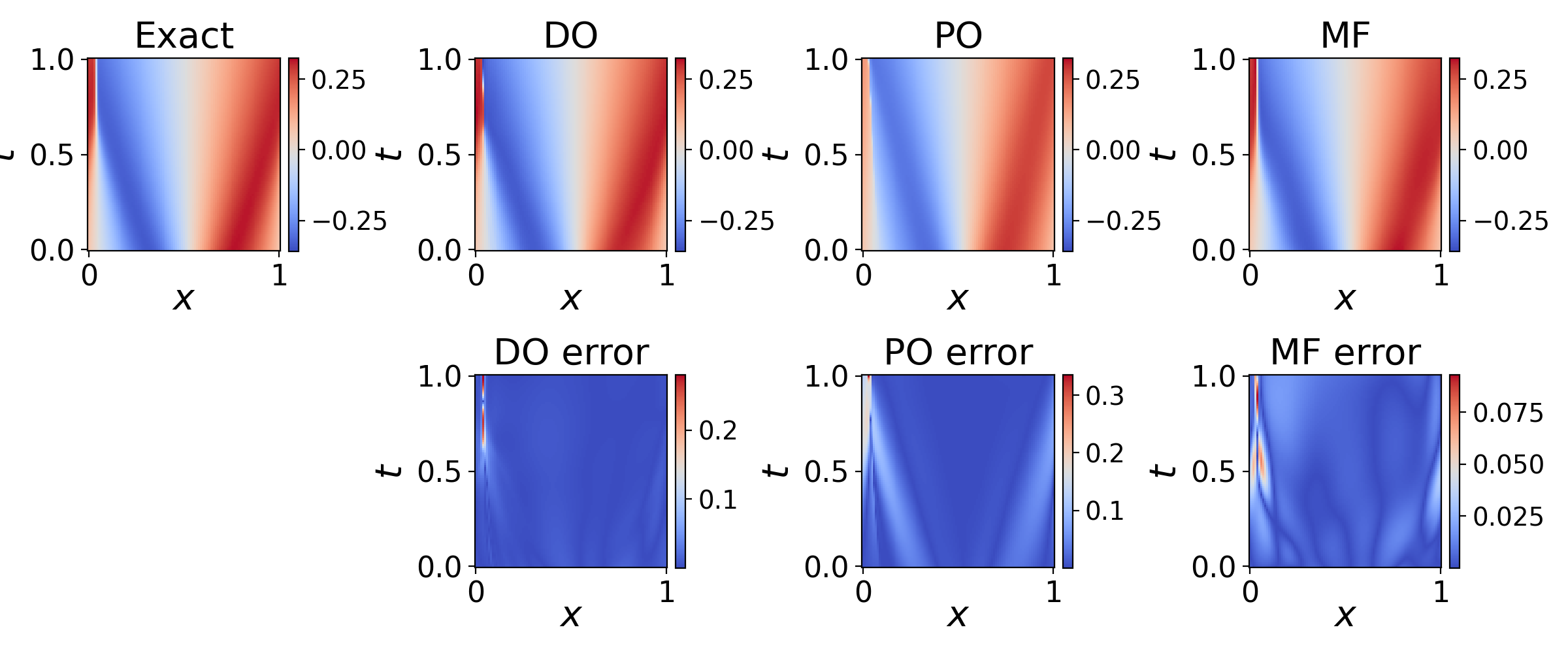}\\
\includegraphics[width=.8\textwidth]{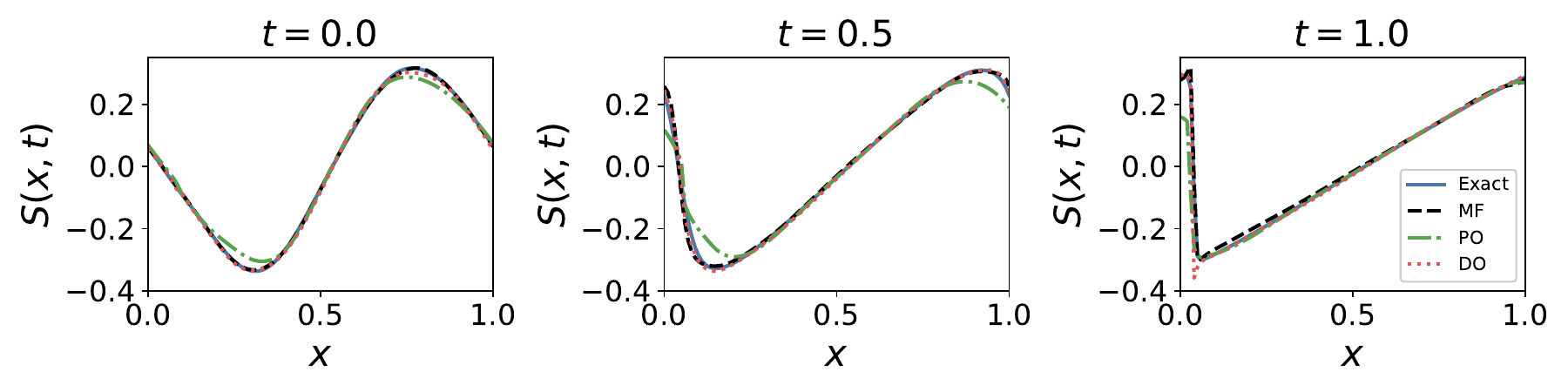}
\end{subfigure}
\caption{Physics-informed multifidelity: viscous Burgers equation. Exact solution and results from the data-only training (DO), physics-only training (PO), and multifidelity training (MF) for (a) $\nu = 10^{-3}$ and (b)$\nu = 10^{-4}$. The errors are the absolute errors between the exact solution and the method output. These examples use the low-fidelity dataset without noise.} \label{fig:burgers}
\end{figure}

\textcolor{black}{The results of the physics-only, data-only, and multifidelity training are given in Tab. \ref{tab:Burger-results}. We can see that the data-only predictions are very accurate, especially when no noise is included. With noise, the data-only predictions show larger relative errors for small values of the viscosity. As reported in \cite{wang2021learning, wang2021improved}, at the smallest value of the viscosity the physics-only training has a large relative error. This is also seen in our results in the large relative error in Table \ref{tab:Burger-results} for the physics-only training with small $\nu$, and in Fig. \ref{fig:burgers}, which gives results of typical outputs from the training for $\nu = 10^{-3}$ and $\nu = 10^{-4}$ with $N_L = 1000$. At the lowest viscosity, the physics-only method underpredicts the steep gradient while the data-only method overshoots the gradient.  The errors are concentrated in the area with the highest gradients. }

\textcolor{black}{In the multifidelity training, we train by enforcing physics as the high-fidelity model and data with and without noise as the low-fidelity data. In all cases, the multifidelity method is able to  reduce the relative error compared with physics-only training. Interestingly, the results of the multifidelity method with and without noise are quite similar, indicating that the multifidelity method can overcome noisy low-fidelity data by correcting the low-fidelity output by enforcing the PDE. }

\textcolor{black}{To further study the sensitivity of the multifidelity training with respect to the low-fidelity dataset, for the most difficult case, $\nu = 10^{-4}$, we reduce the size of the low-fidelity training set from $N_L = 1000$ to $N_L = 200$.  We find that, unsurprisingly, the data-only errors increase significantly due to the decrease in available training data. This is especially noticeable for the data-only training with noisy data, where the relative error increases from $10.65\%$ to $26.11\%$. In contrast, the relative error with the multifidelity method has a significantly smaller increase of about $3\%$. Even with the smaller and noisy low fidelity training set, the multifidelity method still represents a significant reduction in relative error compared with the physics-only training for $\nu = 10^{-4}$.}

\begin{figure}[h]
\begin{subfigure}{\textwidth}
\centering
\caption{$\nu = 10^{-4}$, $N_L = 1000$}
\includegraphics[width=\textwidth]{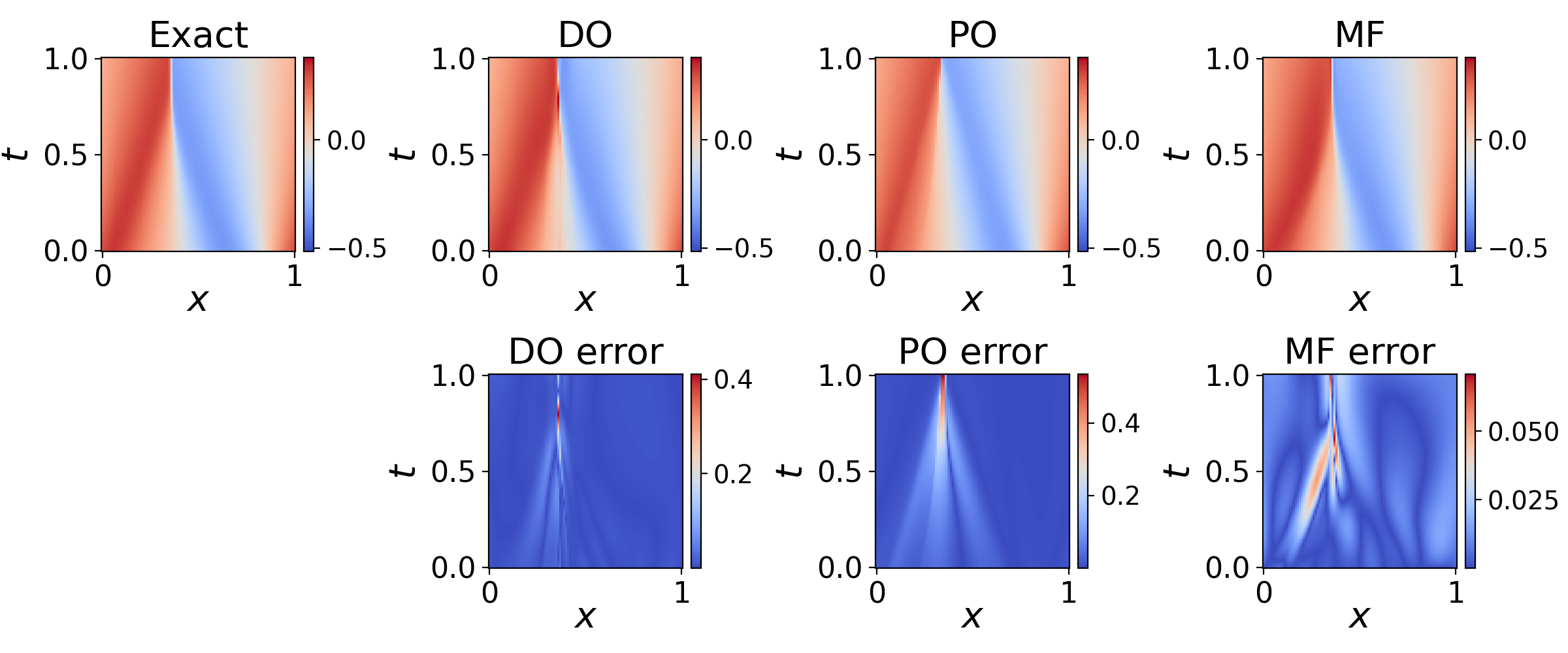}\\
\includegraphics[width=.8\textwidth]{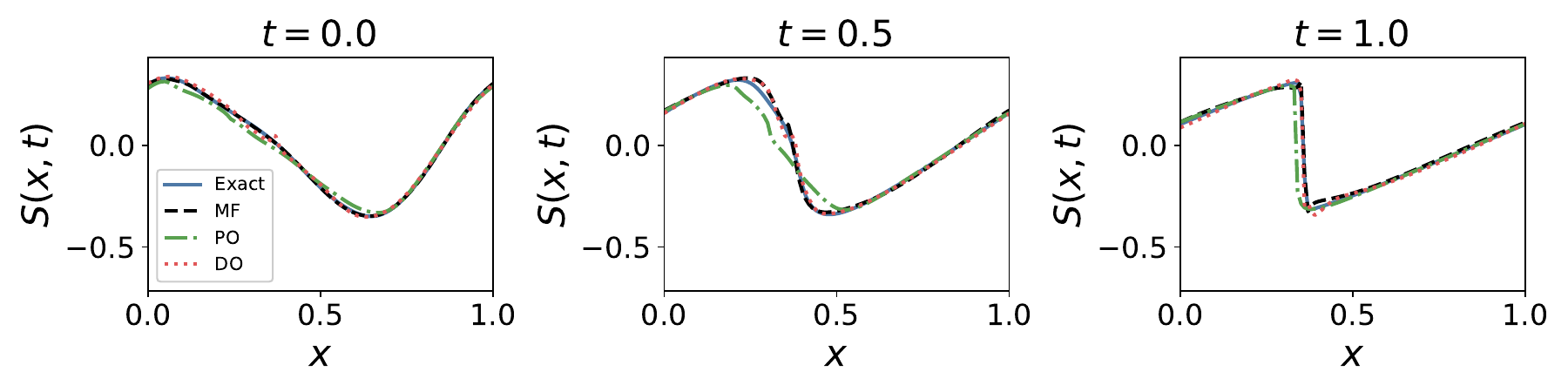}
\end{subfigure}
\begin{subfigure}{\textwidth}
\centering
\caption{$\nu = 10^{-4}$, $N_L = 200$}
\includegraphics[width=\textwidth]{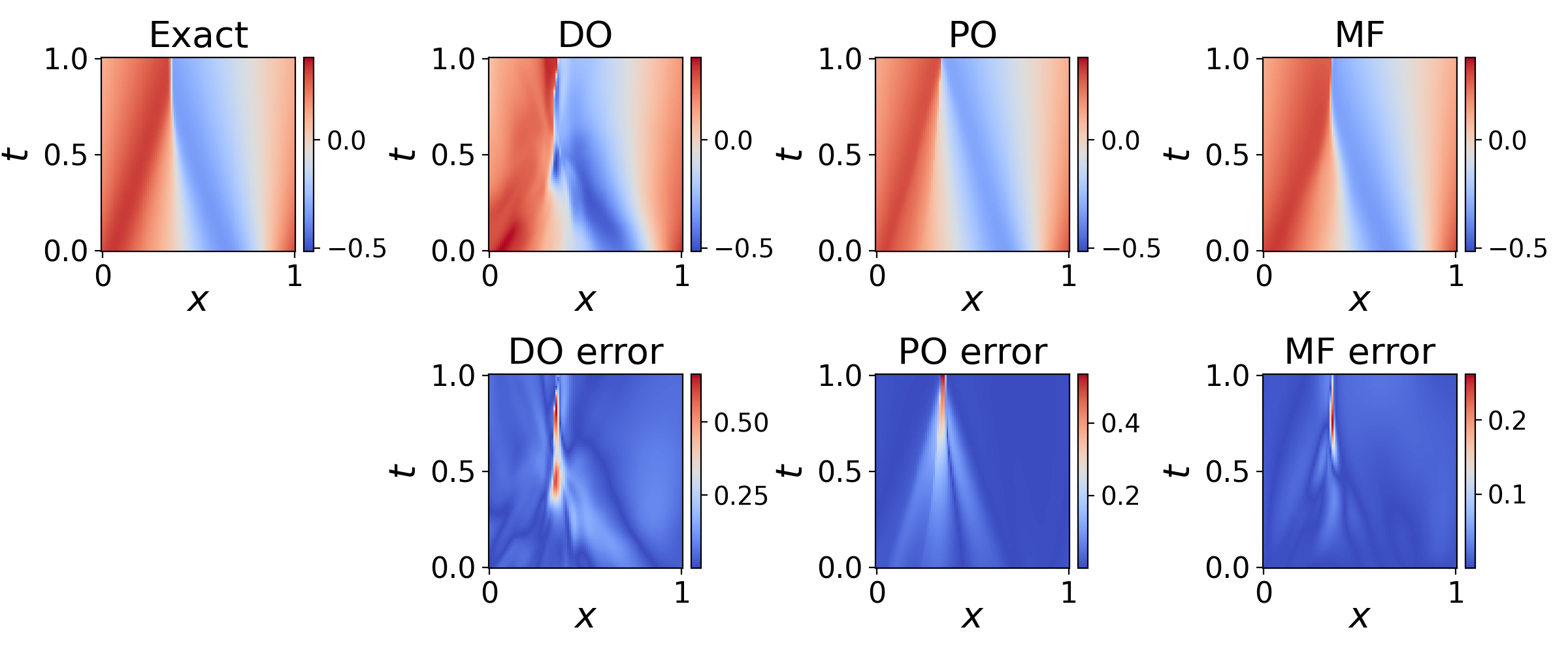}\\
\includegraphics[width=.8\textwidth]{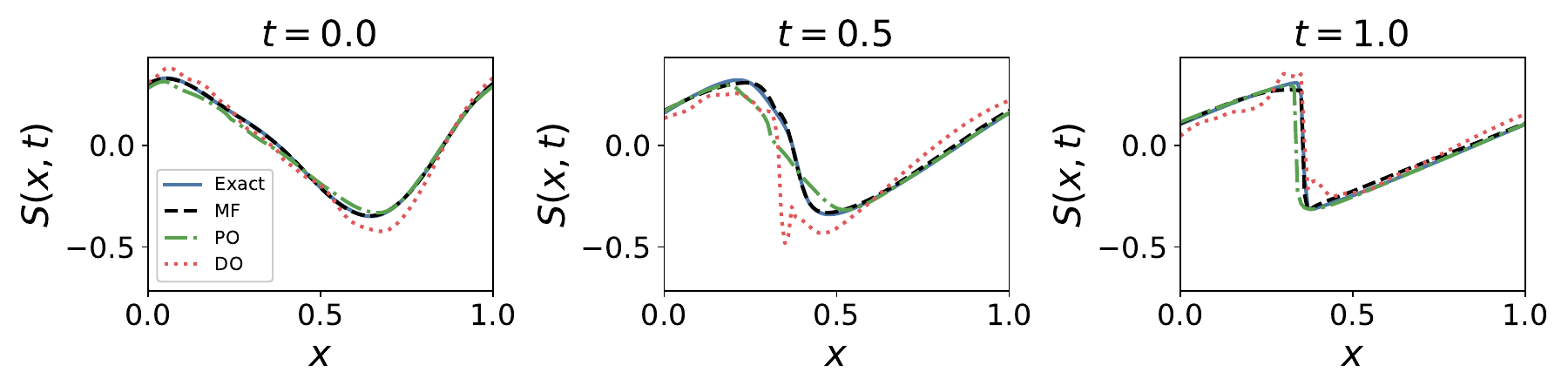}
\end{subfigure}
\caption{Physics-informed multifidelity: viscous Burgers equation. Exact solution and results from the data-only training (DO), physics-only training (PO), and multifidelity training (MF) for (a) $\nu = 10^{-4}$ with $N_L = 1000$ and (b) $\nu = 10^{-4}$ with $N_L = 200$. The errors are the absolute errors between the exact solution and the method output. These examples use the low-fidelity dataset with noise.} \label{fig:burgers_noisy}
\end{figure}

\textcolor{black}{In Fig. \ref{fig:burgers_noisy} we show the results of training with the noisy low-fidelity dataset with $N_L =200$ and $N_L =1000$ for $\nu = 10^{-4}.$ With the larger number of low-fidelity training sets in Fig. \ref{fig:burgers_noisy}a, the multifidelity output is comparable to the multifidelity DeepONet trained with data without noise in Fig. \ref{fig:burgers}b. The data-only error increases significantly with fewer low-fidelity noisy training sets.}

A table of the computational cost for the different methods is given in Tab. \ref{tab:Burger-training-time} (\ref{sec:burgers_app}). While the multifidelity training comes at increased computational cost over the physics-only model, due to the additional computations needed for the three subnetworks, the additional training cost results in higher accuracy.

\begin{table}[]
    \centering
    \begin{tabular}{c|P{25mm}|P{25mm}|P{25mm} | P{25mm}}
    \hline
     Parameters  & $\nu = 10^{-2}$ \newline $N_L = 1000$ & $\nu = 10^{-3}$ \newline $N_L = 1000$ & $\nu = 10^{-4}$ \newline $N_L = 1000$ & $\nu = 10^{-4}$ \newline $N_L = 200$ \\ \hline
    Data-only & $1.02\% \pm  0.81\%$ & $2.46 \% \pm  1.67\%$ & $7.64 \% \pm 2.66 \%$ & $ 13.57\% \pm  7.40\%$ \\
    Data-only with noise & $4.44\% \pm  3.48\%$ & $6.50 \% \pm 3.45 \%$ & $10.63 \% \pm 5.54 \%$ & $ 26.11\% \pm  15.38\%$\\
    Physics-only & $ 3.97\% \pm  5.71\%$ &$ 8.66\% \pm  6.47\%$ &$ 23.63\% \pm  10.22\%$ & --\\
    Multifidelity &$ 2.81\% \pm  1.81\%$ & $ 6.25\% \pm  2.20\%$&  $ 7.05\% \pm  3.01\%$ & $ 9.70\% \pm  4.60\%$ \\
    Multifidelity with noise & $ 2.89\% \pm  1.70\%$& $ 6.65\% \pm  2.48\%$&$ 7.03\% \pm  3.10\%$ & $ 10.16\% \pm  5.55\%$\\ \hline
    \end{tabular}
    \caption{Physics-informed multifidelity: viscous Burgers equation mean relative $L_2$ errors. The physics-only and multifidelity cases all use $N_H = 1000$. \textcolor{black}{Note that the physics-only case does not use any low-fidelity data, so the number of low fidelity samples, $N_L$, does not impact the results. Therefore, we do not report a result for the physics-only case with $N_L = 200$ in the last column, because it would be identical to the physics-only case with $N_L = 1000$.} }
    \label{tab:Burger-results}
\end{table}


\textcolor{black}{This section highlights several key features of our flexible multifidelity framework, and illustrates when you may choose to use single fidelity or multifidelity training. Not all problems require multifidelity, and a great deal of recent work has shown the potential of single fidelity physics-informed DeepONets, see, for example, \cite{wang2021learning, wang2021improved, goswami2022physics, goswami2022physics2, koric2023data}, among many others. However, in some cases physics-informed DeepONets can still struggle to train. In those cases, the multifidelity DeepONet can use a small amount of low fidelity data to improve the training. We have shown that the method is not sensitive to noise in the data, which adds additional flexibility if the only available low fidelity data is noisy. If a very accurate result is needed, a multifidelity DeepONet may increase accuracy over single fidelity physics-informed DeepONets. }

\section{Discussion}\label{sec:discussion} 
In this work, we presented a new composite DeepONet framework for learning operators with multifidelity data. The method uses a large set of inexpensive low-fidelity data and either a small amount of high-fidelity data or enforces physics on the system. In multiphysics and complex system problems it is common to have low-order numerical solvers available, but attaining high-order measurements or simulations is costly. This method allows for using both low- and high-fidelity data to achieve higher accuracy. When PDEs describing the system are known, the physics-informed multifidelity DeepONet allows for training with low-fidelity data and using a physics-informed DeepONet for the high-fidelity training. 

Many extensions of this work are possible to accommodate the data available for novel applications. 
 For example, if a solver for the low-fidelity data is available, either through a surrogate model such as a DeepONet or through a numerical method, this solver can be incorporated instead of the low-fidelity subnetwork presented here, as in \cite{ahmed2023multifidelity}. A discussion of this case is given in \ref{sec:noncomp}. The multifidelity DeepONet framework can also be used to include more than two fidelities of data by incorporating additional fidelity data into the branch net for training. If partial knowledge of the physics is available, the loss function could be modified to enforce the incomplete physical model on the low-fidelity output, then correct the incomplete model with data-driven high-fidelity DeepONets. This case is particularly applicable in cases such as fluid modeling, where sparse but high-fidelity experimental data is available. \textcolor{black}{Future extensions of this work include incorporating weighting schemes in the loss function to improve training, along the lines of \cite{mcclenny2020self, wang2021improved}. Currently, choosing the weights in the loss function can require some knowledge of the scales in the problem. This is not a problem unique to multifidelity DeepONet training, and is also an active area of research for single fidelity physics-informed DeepONets.}

\section{Acknowledgements}
The authors wish to thank P. Perdikaris, Q. He, and L. Lu for helpful discussions, K. C. Sockwell for co-developing the ice-sheet code, and T. Hillebrand for generating the Humboldt grid.

The work of AH, GEK and PS is supported by the U.S. Department of Energy, Advanced Scientific Computing Research program, under the Physics-Informed Learning Machines for Multiscale and Multiphysics Problems (PhILMs) project (Project No. 72627). Support for Mauro Perego was provided through the Scientific Discovery through Advanced Computing (SciDAC) program funded by the US Department of Energy (DOE), Office of Science, Advanced Scientific Computing Research and Biological and Environmental Research Programs.  Pacific Northwest National Laboratory (PNNL) is a multi-program national laboratory operated for the U.S. Department of Energy (DOE) by Battelle Memorial Institute under Contract No. DE-AC05-76RL01830. The computational work was performed using PNNL Institutional Computing at Pacific Northwest National Laboratory. Sandia National Laboratories is a multimission laboratory managed and operated by National Technology and Engineering Solutions of Sandia, LLC., a wholly owned subsidiary of Honeywell International, Inc., for the U.S. Department of Energy’s National Nuclear Security Administration under contract DE-NA-0003525.

\clearpage

\bibliographystyle{elsarticle-num} 
\bibliography{Amanda-bib}

\appendix

\section{Notation and abbreviations}

\begin{tabular}{l l }
$N_L$ & Number of low-fidelity input samples in the training data set\\
$M_L$ & Number of locations for evaluating input samples to the low-fidelity network \\
$P_L$ & Length of the low-fidelity output \\
$N_H$ & Number of high-fidelity input samples in the training data set\\
$M_H$ & Number of locations for evaluating input samples to the high-fidelity network \\
$P_H$ & Length of the high-fidelity output \\
$P_{BC}$ & Number of boundary points \\
$P_{p}$ & Number of collocation points for evaluating the PDE residual\\
$\theta$ & All trainable parameters of the multifidelity DeepONet system\\
$\mathbf{u}$ & An input function \\
$\mathcal{G}(\mathbf{u})(\mathbf{x})$ & An operator\\
$\mathcal{G}^\theta(\mathbf{u})(\mathbf{x})$ & A DeepONet representation of an operator\\
$\mathcal{F}_{LF}^\theta(\mathbf{u})(\mathbf{x})$ & The output of the low-fidelity DeepONet  \\
$\mathcal{F}_{nl}^\theta(\mathbf{u})(\mathbf{x})$ & The output of the nonlinear DeepONet \\
$\mathcal{F}_l^\theta(\mathbf{u})(\mathbf{x})$ & The output of the linear DeepONet \\
SSA &  Shallow Shelf Approximation\\
MOLHO & Mono-Layer Higher-Order model \\
MSE & Mean squared error\\
ODE & Ordinary differential equation \\
PDE & Partial differential equation\\
DO & Data-only \\
PO & Physics-only \\
LF & Low-fidelity \\
HF & High-fidelity \\
MF & Multifidelity \\
SF & Single fidelity
    \end{tabular}

\section{Additional data-driven examples} \label{sec:add_data-drive}

\subsection{One-dimensional, correlation with $u$}\label{1d_corr}
\begin{figure}[ht]
\centering
\begin{subfigure}{0.4\textwidth}
\centering
\caption{$a = 10.7316$}
\includegraphics[width=\textwidth]{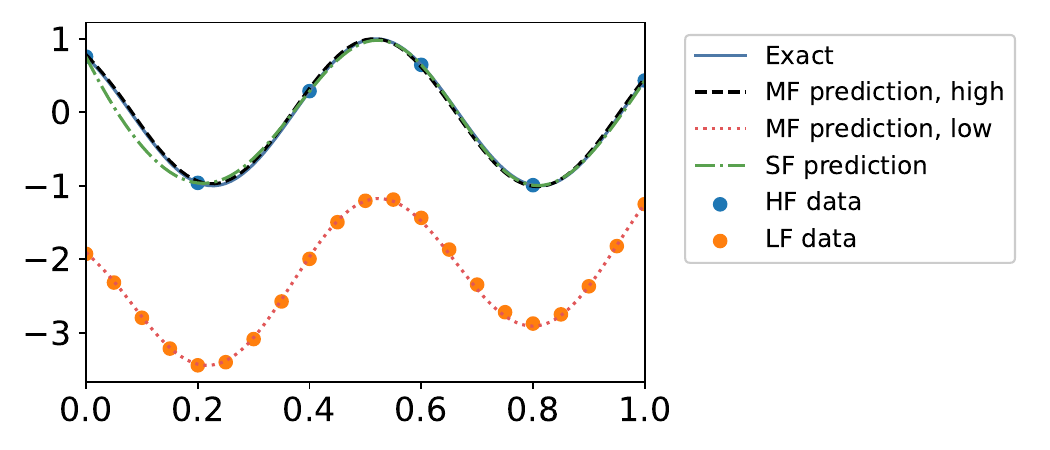}
\end{subfigure}\hspace{0.02\textwidth}
\begin{subfigure}{0.4\textwidth}
\centering
\caption{$a = 13.4684$}
\includegraphics[width=\textwidth]{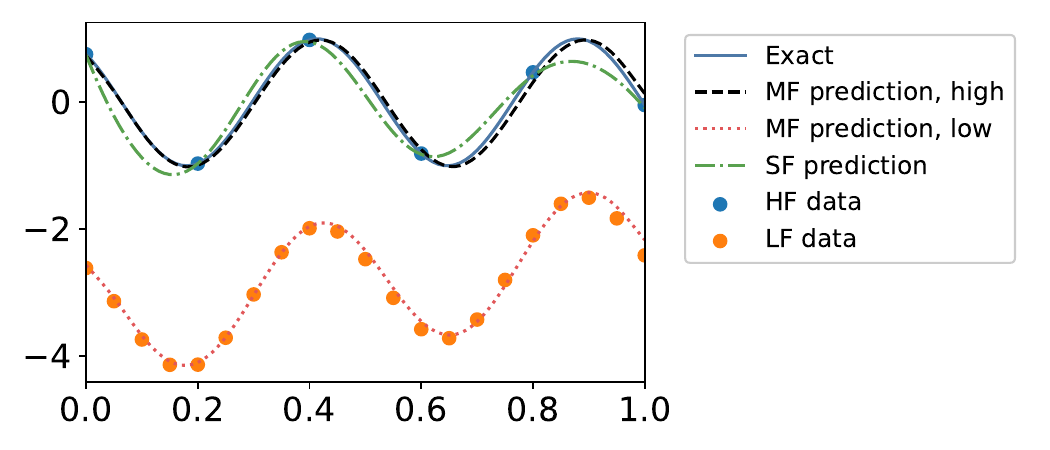}
\end{subfigure}\\
\begin{subfigure}{.3\textwidth}
\centering
\caption{SF error}
\includegraphics[width=\textwidth]{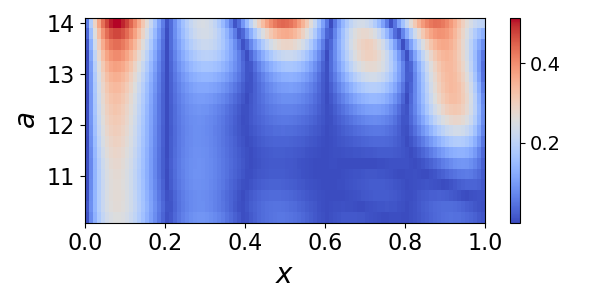}
\end{subfigure}
\begin{subfigure}{.3\textwidth}
\centering
\caption{MF error, high-fidelity data}
\includegraphics[width=\textwidth]{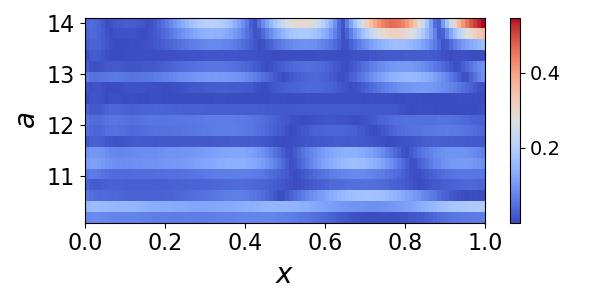}
\end{subfigure}
\begin{subfigure}{.3\textwidth}
\centering
\caption{MF error, low-fidelity data}
\includegraphics[width=\textwidth]{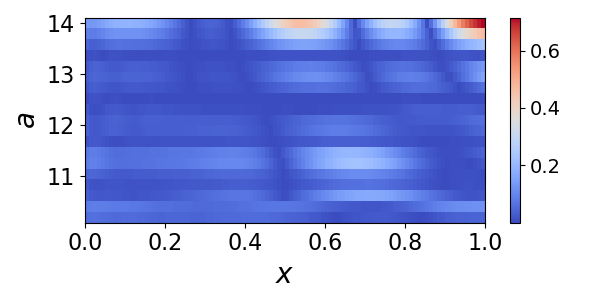}
\end{subfigure}
\caption{Data-driven multifidelity: one-dimensional, correlation with $u$. (a-b) Results of the single fidelity and multifidelity predictions of the high-fidelity data.  (c)  Single-fidelity relative error as a function of $a$ and $x$, (d) multifidelity high-fidelity prediction relative error as a function of $a$ and $x$, and (e) multifidelity low-fidelity prediction relative error as a function of $a$ and $x$.} \label{fig:1d_linwithb_results}
\end{figure}

\begin{figure}[H]
\begin{minipage}{0.55\textwidth}
\begin{table}[H]
    \centering
        \begin{tabular}{c|c }
    \hline
    Parameter & Value \\
    \hline
      LF data   &  $M_L = P_L= 21$ \\
      HF data  & $M_H = P_H = 6$  \\
      Number of datasets & $N_L = N_H = 5$ values of $a$ \\
      $\lambda_1$ & 0.1 \\
      $\lambda_2$ & 1 \\
      $\lambda_3$ & $1\times 10^{-1}$ \\
      $\lambda_4$ & $1\times 10^{-3}$ \\
      SF learning rate & (1e-3, 2000, 0.9)\\
      MF learning rate & (1e-3, 5000, 0.97)\\
      SF network size & 3 layers, 30 neurons \\
      MF low-fidelity network size &  3 layers, 30 neurons\\
      MF linear network size & 1 layer, 5 neurons \\
      MF nonlinear network size & 2 layers, 20 neurons\\
          \hline
    \end{tabular}
    \caption{Training parameters for the one-dimensional problem with correlation as a function of $u$.}
    \label{tab:1d_linwithb}
\end{table}
\end{minipage}
\begin{minipage}{0.42\textwidth}
\centering
\includegraphics[width=\textwidth]{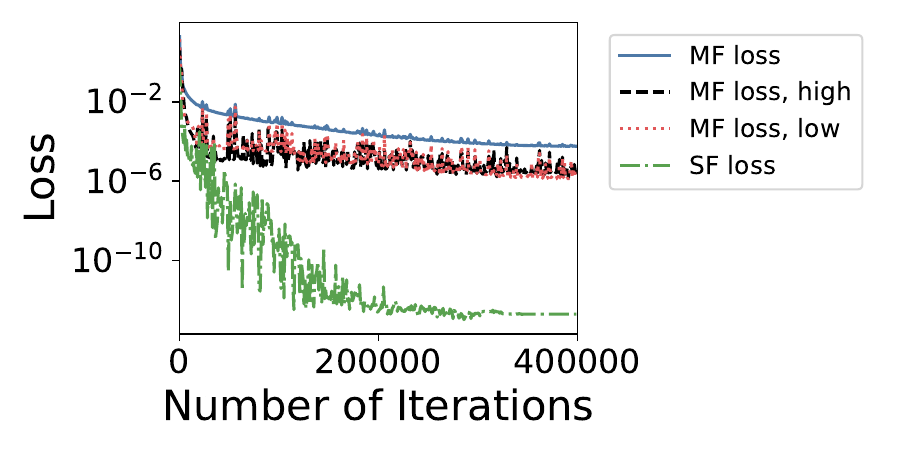}
\caption{Training loss for the one-dimensional problem with correlation as a function of $u$.} \label{fig:1d_linwb_loss}
\end{minipage}
\end{figure}

The multifidelity data-driven training is able to capture complex, nonlinear correlations between the low- and high-fidelity datasets. To illustrate this, we consider a case where the correlation depends on the input function, $u$: 
\begin{align}
    y_L(u)(x) &= \sin(u) + x - 0.25u\\
    y_H(u)(x) &= \sin(u)  \\
    u &= ax-4
\end{align}
for $x \in [0, 1]$ and $a \in [10, 14]$. We have $y_H(u)(x) = y_L(u)(x)-x + 0.25u =  y_L(u)(x)-x + 0.25(ax-4)$. Parameters are given in Tab. \ref{tab:1d_linwithb} and results in Fig. \ref{fig:1d_linwithb_results}. The single fidelity case has large errors across the $(x, a)$ domain. The error in the high-fidelity output from the multifidelity training is concentrated where the low-fidelity error is the highest.

\subsection{Two-dimensional, linear correlation}\label{2d_linear}

\begin{figure}[ht]
\begin{subfigure}{\textwidth}
\centering
\caption{$a = 8.5211$}
\includegraphics[width=.25\textwidth]{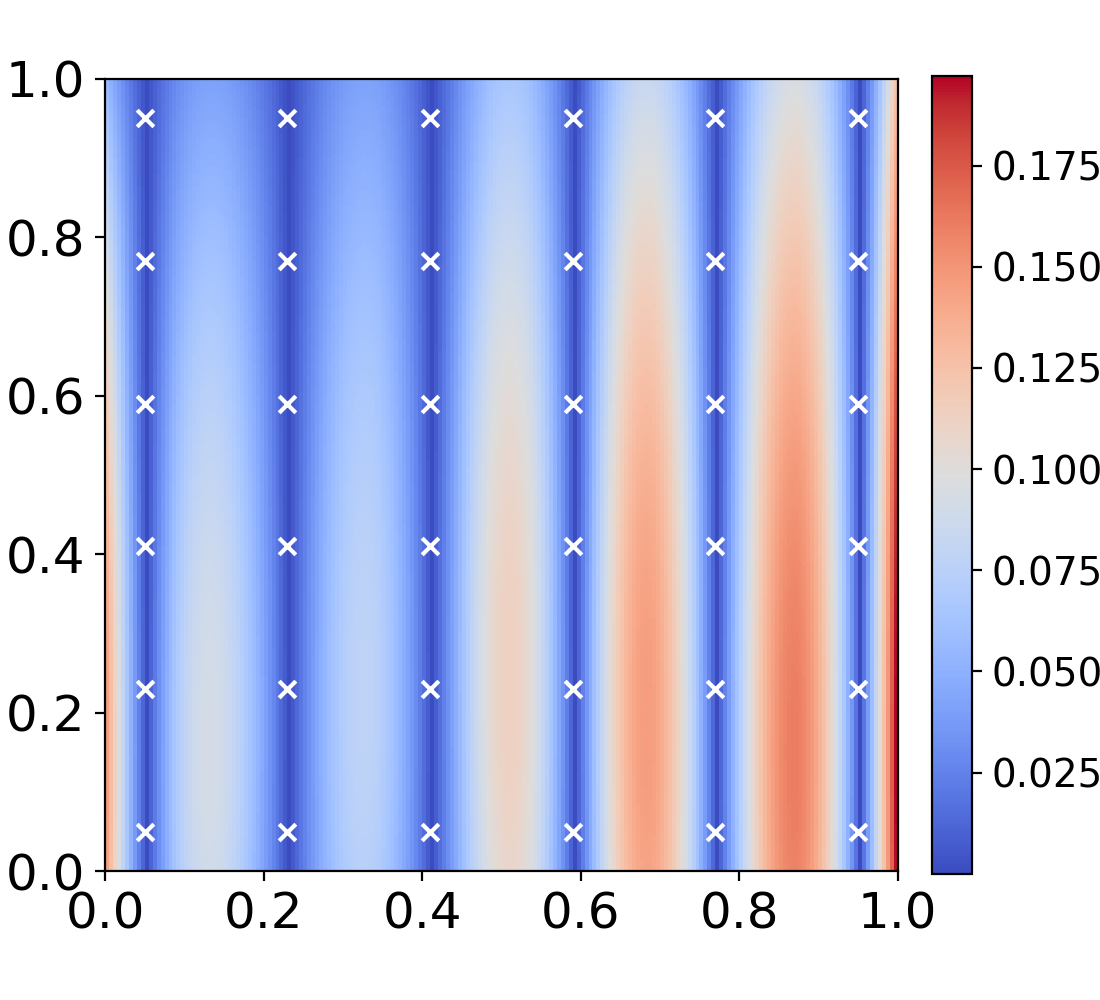}
\hspace{.05\textwidth}
\includegraphics[width=.25\textwidth]{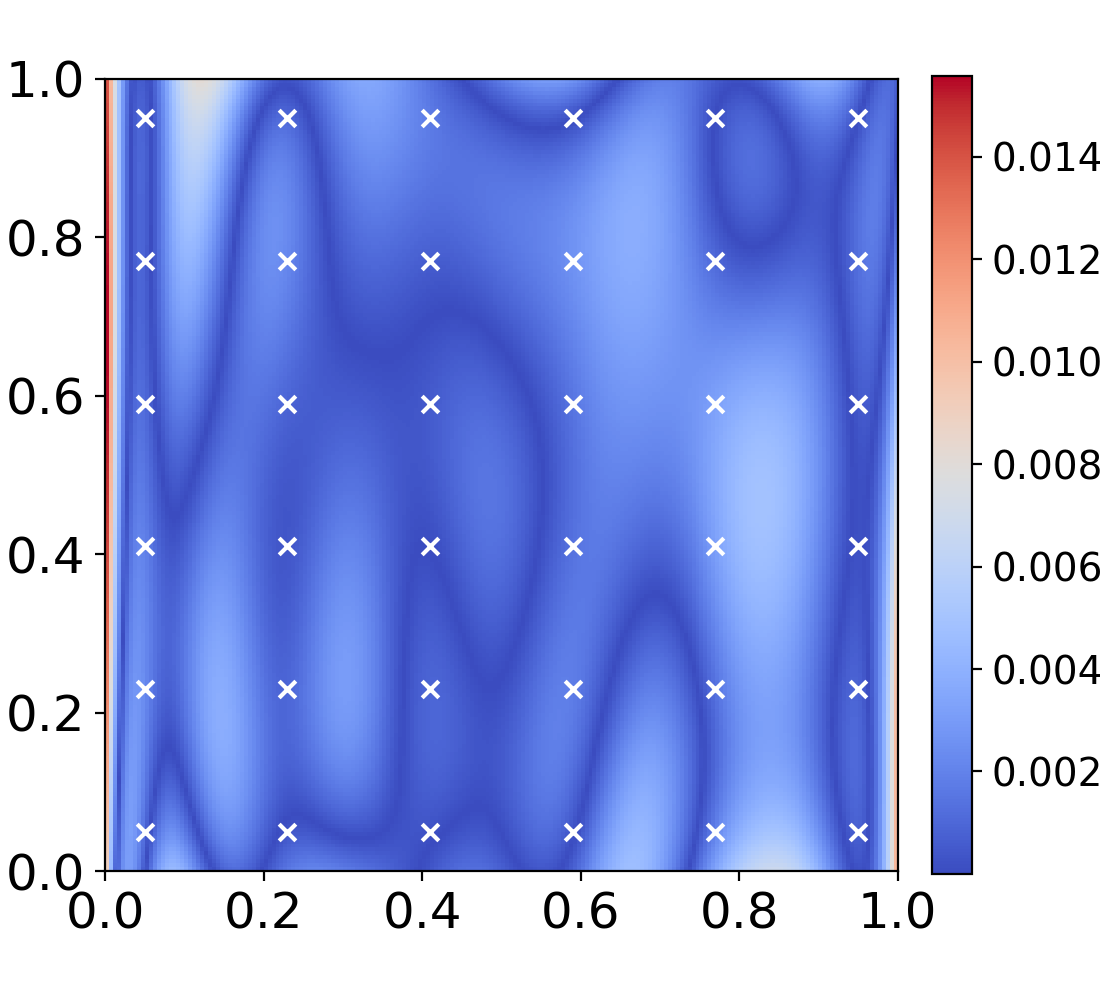}
\hspace{.05\textwidth}
\includegraphics[width=.25\textwidth]{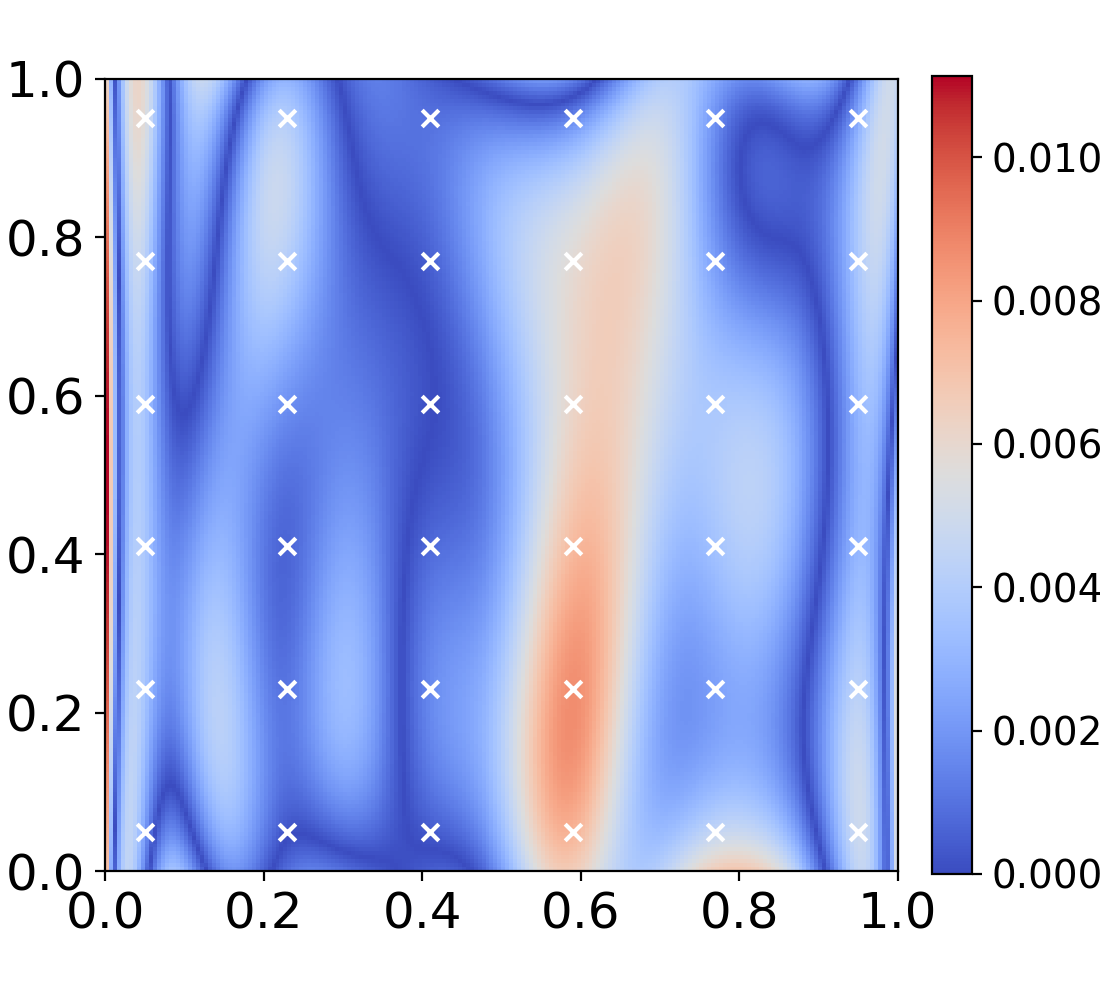}
\end{subfigure}
\begin{subfigure}{\textwidth}
\centering
\caption{$a = 9.5737$}
\includegraphics[width=.25\textwidth]{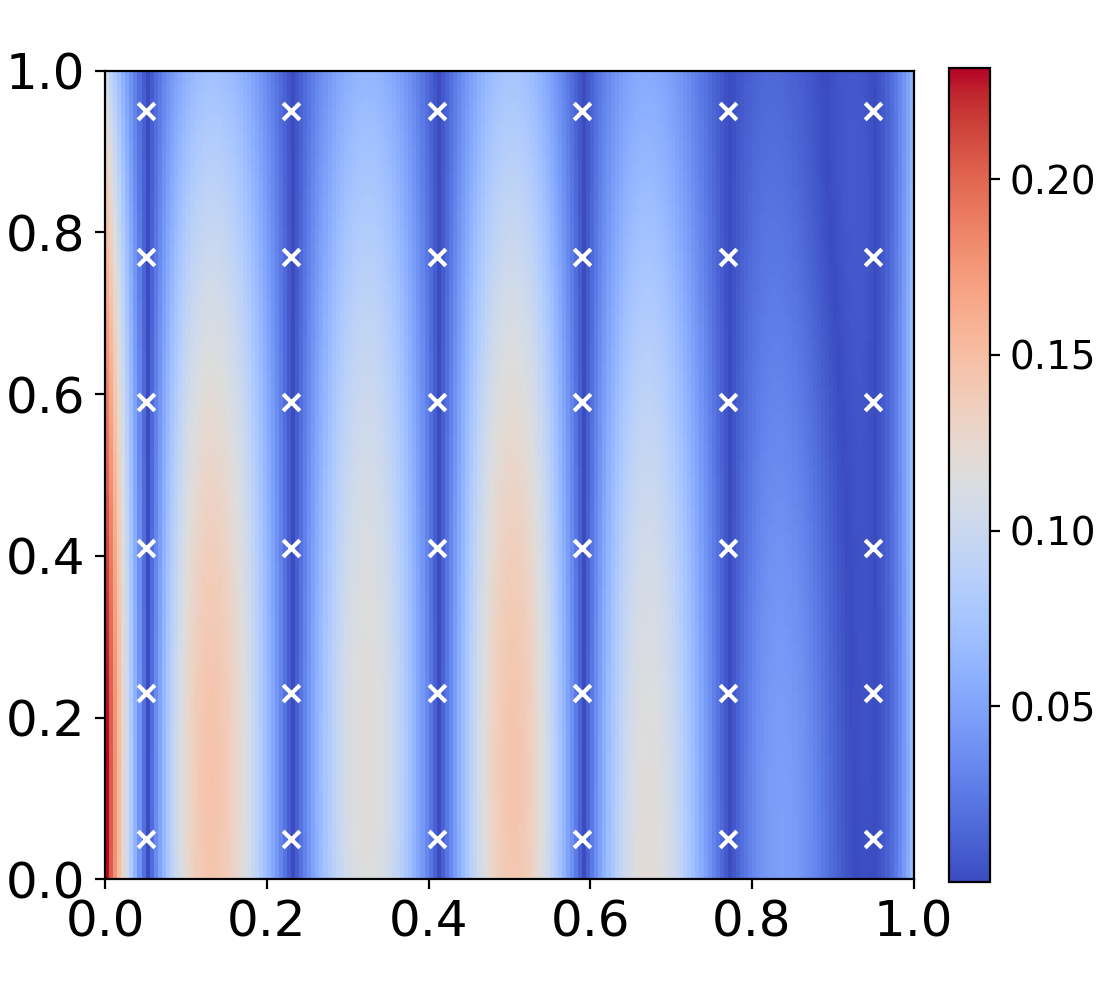}
\hspace{.05\textwidth}
\includegraphics[width=.25\textwidth]{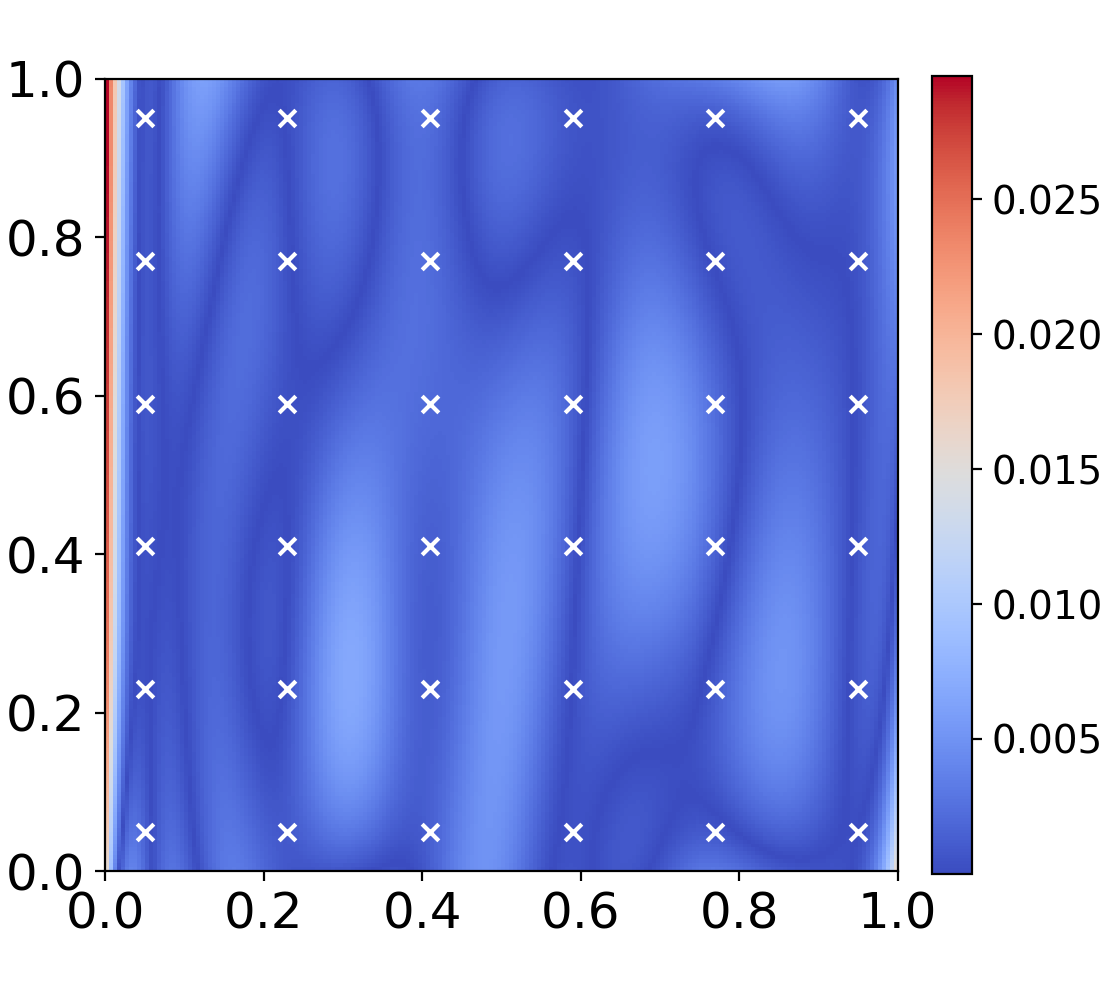}
\hspace{.05\textwidth}
\includegraphics[width=.25\textwidth]{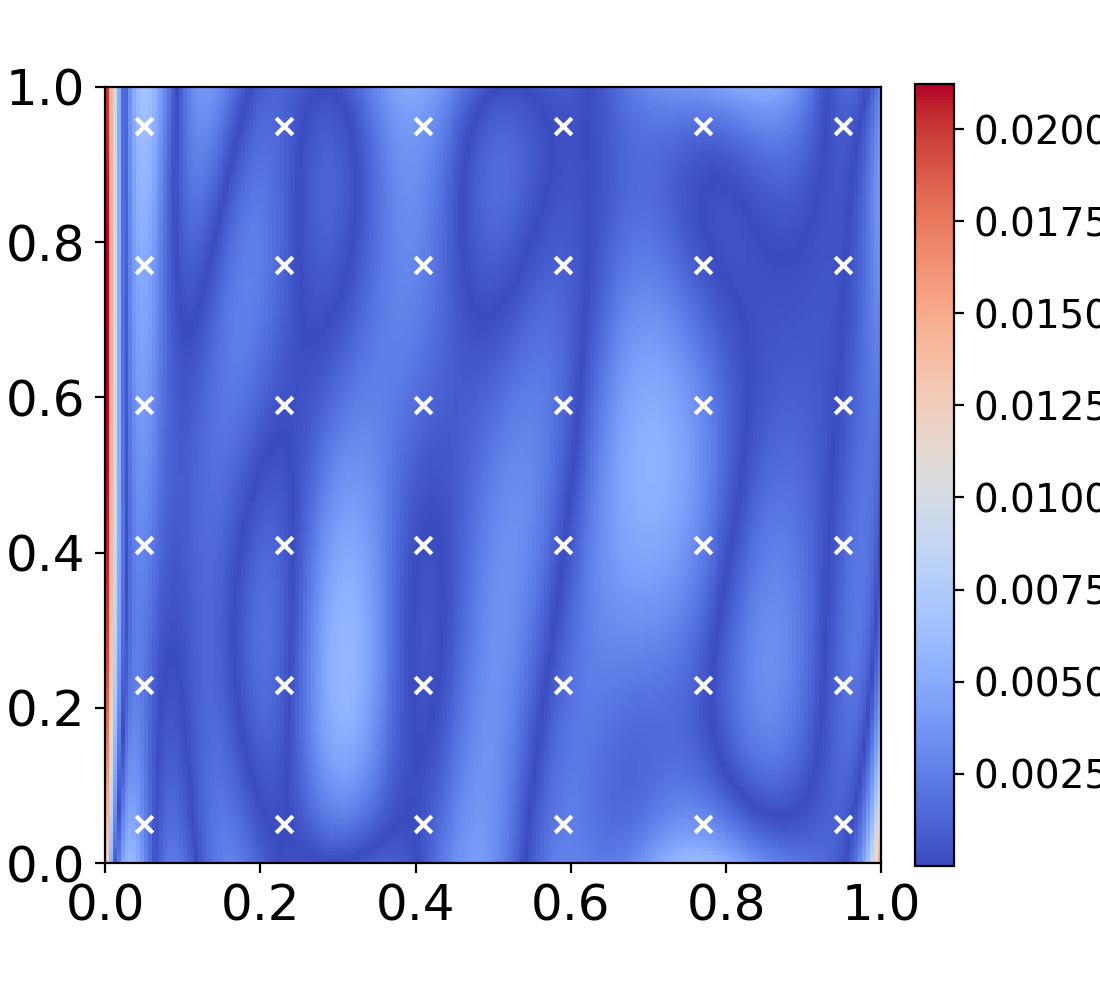}
\end{subfigure}
\caption{Data-driven multifidelity: two-dimensional, linear correlation. (a) Absolute error of the high-fidelity prediction, multifidelity prediction of the high-fidelity data, and multifidelity prediction of the low-fidelity data for $a = 8.5211$. (b) Absolute error of the high-fidelity prediction, multifidelity prediction of the high-fidelity data, and multifidelity prediction of the low-fidelity data for $a = 9.5737$.   The high-fidelity data points are shown in white for illustration.} \label{fig:2d_lin_results}
\end{figure}

\begin{figure}[H]
\centering
\begin{subfigure}{0.3\textwidth}
\centering
\caption{}
\includegraphics[width=\textwidth]{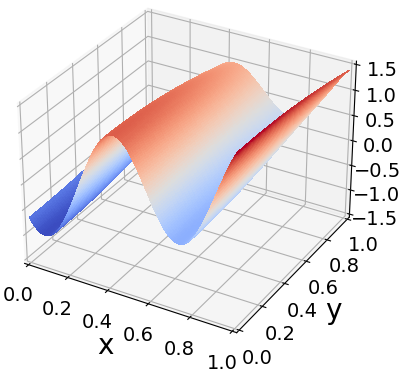}
\end{subfigure}\hfill
\begin{subfigure}{0.3\textwidth}
\centering
\caption{}
\includegraphics[width=\textwidth]{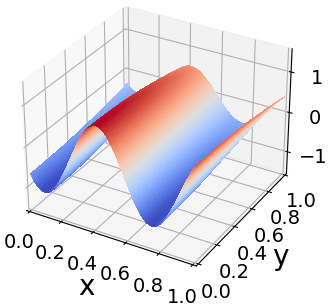}
\end{subfigure}
\caption{(a) Exact low-fidelity function, $z_L$, and (b) exact high-fidelity function, $z_H$, for $ a = 9.5737$. } \label{fig:2d_lin_exact}
\end{figure}

\begin{figure}[H]
\begin{minipage}{0.6\textwidth}
\begin{table}[H]
    \centering
            \begin{tabular}{c|c }
    \hline
    Parameter & Value \\
    \hline
      LF data   &  $M_L = P_L= 21^2$, $x, y \in [.02, .98]$ \\
      HF data  & $M_H = P_H = 6^2$, $x, y \in [.05, .95]$  \\
      Number of datasets & $N_L = N_H = 20$ values of $a$ \\
      $\lambda_1$ & 1 \\
      $\lambda_2$ & 1 \\
      $\lambda_3$ & $1\times 10^{-3}$ \\
      $\lambda_4$ & $1\times 10^{-4}$ \\
      SF learning rate & (1e-3, 5000, 0.9)\\
      MF learning rate & (1e-3, 5000, 0.97)\\
      SF network size & 3 layers, 40 neurons \\
      MF low-fidelity network size &  3 layers, 30 neurons\\
      MF linear network size & 1 layer, 5 neurons \\
      MF nonlinear network size & 2 layers, 20 neurons\\
          \hline
    \end{tabular}
    \caption{Training parameters for the two-dimensional problem with linear correlation.}
    \label{tab:2d_lin}
\end{table}
\end{minipage}
\begin{minipage}{0.35\textwidth}
\centering
\includegraphics[width=\textwidth]{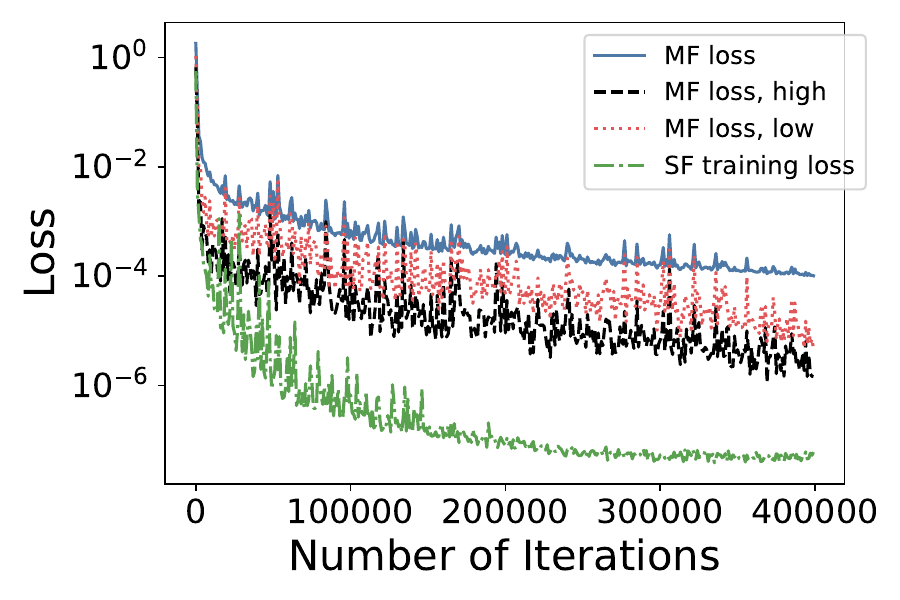}
\caption{Training loss for the two-dimensional problem with linear correlation.} \label{fig:2d_lin_loss}
\end{minipage}
\end{figure}

We consider a two-dimensional problem with a linear correlation between the low-fidelity and high-fidelity data: 
\begin{align}
    z_L(u)(x, y) &= \cos(u)\cos(y) + x \\
    z_H(u)(x, y) &= \cos(u)\cos(y)  \\
    u &= ax-4
\end{align}
for $x, y \in [0, 1]$ and $a \in [8, 10]$. We have $z_H(u)(x, y) = z_L(u)(x, y)-x$ (see Fig. \ref{fig:2d_lin_exact}.) The training parameters are given in Tab.\ref{tab:2d_lin} and the results are given in Fig. \ref{fig:2d_lin_results}. Note that the multifidelity prediction results in absolute errors approximately one order of magnitude smaller than the single fidelity prediction. 
The linear correlation found is: 
\begin{equation}
        \mathcal{F}_l(u)(x, y) = 0.9973\mathcal{F}_{LF}(u)(x, y) -0.9392x-0.0039y -0.0032x\mathcal{F}_{LF}(u)(x, y)-0.0048y\mathcal{F}_{LF}(u)(x, y)-0.0256.
\end{equation}

\section{Training parameters} \label{sec:training_params}
In this section we provide the training parameters used in each test case in the main text for reproducibility. 

All test cases are trained with the \texttt{adam} optimizer. Unless noted, we use the hyperbolic tangent activation function. \textcolor{black}{The hyperbolic tangent activation function was chosen because it had the best performance in our tests, and is a standard choice for physics-informed DeepONet training. A study of the impact of the action function is outside the scope of this work.}

\begin{landscape}
\begin{table}[H]
    \centering
        \begin{tabular}{c|c | c| c | c}
    \hline
    Case & Sec. \ref{1d_jump} & Sec. \ref{2d_nonlinear} & Sec. \ref{sec:multiresolution} & Sec. \ref{sec:multiorder}\\
    \hline
      LF data   &  $M_L = P_L= 38$ & $M_L = P_L= 21^2$ & $M_L = P_L= 2 \cdot 15^2$    & $M_L = P_L= 2\cdot 1426$ \\
      HF data  & $M_H = P_H = 5$ & $M_H = P_H = 6^2$ & $M_H = P_H = 2\cdot 41^2$ & $M_H = P_H = 2\cdot 1426$  \\
      Number of SF datasets & $N_H = 10$& $N_H = 20$  & $N_H = 20$ or $50$ & $N_H = 20$\\
     Number of LF datasets & $N_L = 20$  & $N_L = 20$  & $N_L = 100$ & $N_L = 80$\\
      Number of HF datasets & $N_H = 10$  & $N_H = 20$  & $N_H = 20$ & $N_H = 20$\\
      $\lambda_1$ & 0.1 & 1 & 1 & 1  \\
      $\lambda_2$ & 1 & 1 & 10  &1 \\
      $\lambda_3$ & $1\times 10^{-1}$ & $1\times 10^{-3}$& $1\times 10^{-3}$& $1\times 10^{-3}$ \\
      $\lambda_4$ & $1\times 10^{-3}$ & $1\times 10^{-4}$& $1\times 10^{-4}$ &$1\times 10^{-4}$\\
      SF learning rate & (1e-3, 2000, 0.9) & (1e-3, 5000, 0.9)& (1e-4, 5000, 0.9)& (5e-5, 5000, 0.9)\\
      MF learning rate & (5e-4, 2000, 0.99)& (1e-3, 5000, 0.97)& (1e-4, 5000, 0.97)& (5e-4, 5000, 0.97)\\
      SF network size & 3 layers, 30 neurons  & 3 layers, 40 neurons & 5 layers, 200 neurons & 4 layers, 150 neurons \\
      MF low-fidelity network size &  3 layers, 40 neurons & 3 layers, 30 neurons & 5 layers, 200 neurons & 4 layers, 150 neurons\\
      MF linear network size & 1 layer, 5 neurons & 1 layer, 5 neurons & 1 layer, 10 neurons & 1 layer, 10 neurons  \\
      MF nonlinear network size & 2 layers, 30 neurons& 2 layers, 20 neurons& 5 layers, 200 neurons& 3 layers, 150 neurons\\
          \hline
    \end{tabular}
    \caption{Training parameters for the data-driven cases presented in the main text. LF denotes low-fidelity, HF denotes high-fidelity, and SF denotes single-fidelity. The learning rate is set with the \texttt{optimizers.exponential\_decay} function in \texttt{JAX} \cite{jax2018github}, which requires three parameters given as: (initial value, decay steps, decay rate). In all cases, $M_L = P_L$ for the low-fidelity dataset and $M_H=P_H$ for the high-fidelity dataet}
    \label{tab:train_params_all}
\end{table}
\end{landscape}

\begin{table}[H]
    \centering
                \begin{tabular}{c|c | c}
    \hline
    Case & Sec. \ref{sec:1d_pi} & Sec. \ref{sec:burgers} \\
    \hline
      LF data   &  $M_L = P_L=21$&  $M_L=21$, $P_L = 21^2$ \\
      HF BC  & $P_{BC} =1$ at $x=0$ & $P_{BC} =100$ at $x=0$ and $x=1$\\
      HF collocation points  & $P_p = 101$& $P_p = 2500$  \\
      Number of HF datasets  & $N_H=10$ values of $a$ & $N_H = 1000$\\
      Number of LF datasets  & $N_L=20$ values of $a$ & $N_L = 1000$ or $200$ \\
      $\lambda_1$ & $1\times 10^{-1}$ & 10\\
      $\lambda_2$ & 1 & 1\\
      $\lambda_3$ & $1\times 10^{-2}$ & $1\times 10^{-6}$\\
      $\lambda_4$ & $1\times 10^{-4}$ & $1\times 10^{-6}$\\
      $\lambda_5$ & 0 & $20$\\
      $\lambda_6$ & $1\times 10^{-2}$& $1$ \\
      SF learning rate & (1e-3, 2000, 0.9) & --\\
      DO learning rate & -- &(1e-3, 2000, 0.9) \\
      PO learning rate & -- & (1e-3, 2000, 0.9) \\
      MF learning rate & (1e-3, 5000, 0.95) & (1e-3, 2000, 0.9)\\
    SF network size & 3 layers, 20 neurons  & --\\
      DO network size & --& 7 layers, 100 neurons \\
      PO network size & -- & 7 layers, 100 neurons \\
      MF low-fidelity network size &  3 layers, 30 neurons &  7 layers, 100 neurons\\
      MF linear network size & 1 layer, 5 neurons & 1 layer, 10 neurons\\
      MF nonlinear network size & 2 layers, 20 neurons & 4 layers, 100 neurons\\
          \hline
    \end{tabular}
    \caption{Training parameters for the physics-informed cases presented in the main text. For the viscous Burgers case from Sec. \ref{sec:burgers}, DO refers to data-only and PO refers to physics-only. Both represent single fidelity training. The learning rate is set with the \texttt{optimizers.exponential\_decay} function in \texttt{JAX} \cite{jax2018github}, which requires three parameters given as: (initial value, decay steps, decay rate). }
    \label{tab:train_params_all_PI}
\end{table}

\subsection{One-dimensional, jump function}\label{1d_jump_app}

\begin{figure}[H]
\centering
\includegraphics[width=0.42\textwidth]{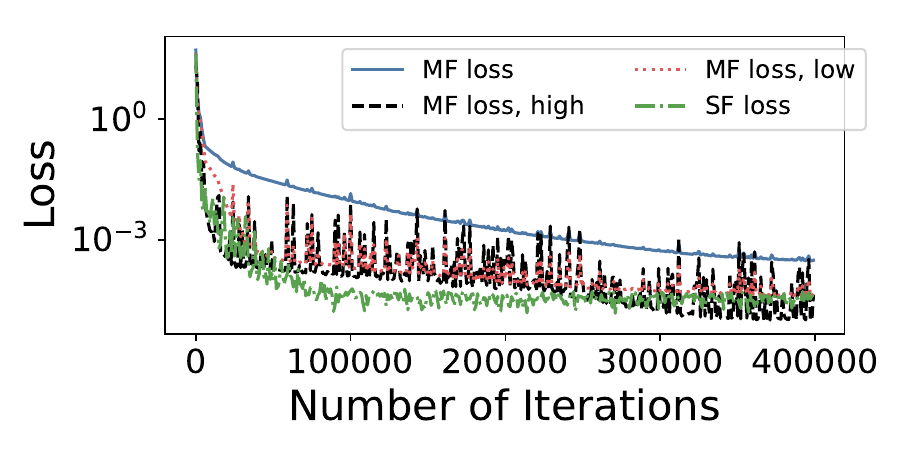}
\caption{One-dimensional jump function training loss.} \label{fig:1d_jump_loss}
\end{figure}

\begin{table}[h!]
    \centering
    \begin{tabular}{c |c|c|c| c| c}
    \hline
 $\lambda_1/\lambda_2 $ & $\lambda_3/\lambda_2 $        & $\mathcal{F}_l(u)(x)$ &$ x$ & $1$ & $x\mathcal{F}_l(u)(x)$  \\ \hline
$0.001$& $1\times 10^{-1}$ &  2.00 & -19.99 & 19.99 & 0.002 \\
$ 0.01$& $ 1\times 10^{-1}$ &2.00 & -20.00 & 20.01 & 0.0003 \\
$  0.1$& $ 1\times 10^{-1}$ &1.95 & -19.15 & 19.34 & 0.048 \\
$1$& $ 1\times 10^{-1}$ & 0.58 & 2.72 & 2.13 & 1.43 \\
$  10$& $1\times 10^{-1}$ &  0.05 & 0.97 & 0.03 & 1.99 \\
$  100$& $1\times 10^{-1}$ & 0.11 & 0.20 & 0.01 & 2.07 \\
$0.1$& $ 1\times 10^{-4}$ & 0.28 & 0.32 & 0.05 & 1.93 \\
$0.1$& $ 1\times 10^{-3}$ &0.24 & 1.43 & 0.25 & 1.78 \\
$0.1$& $ 1\times 10^{-2}$& 0.60 & 2.29 & 2.45 & 1.42 \\
$ 0.1$& $ 1\times 10^{0}$ & 2.00 & -19.96 & 19.97 & 0.002 \\
$ 0.1$ & $ 1\times 10^{1}$ & 2.00 & -19.97 & 19.97 & 0.001 \\
        \hline
    \end{tabular}
    \caption{Mean square errors for the trained multifidelity model with varying hyperparameters in the loss function, Eq. \ref{eq:loss_data}. We fix $\lambda_2 = 1$ and $\lambda_4 = 10^{-4}$. The results in Sec. \ref{1d_jump} take $\lambda_1/\lambda_2 = 0.1$, and $\lambda_3/ \lambda_2 = 1\times 10^{-1}$. The exact learned correlation should be $
    2 \mathcal{F}_{LF}(u)(x)-20x+20 - 0x \mathcal{F}_{LF}(u)(x)$.}
    \label{tab:lin_eq}
\end{table}

\subsection{Two-dimensional, nonlinear correlation}\label{2d_nonlinear_app}

\begin{figure}[H]
\centering
\begin{subfigure}{0.3\textwidth}
\centering
\caption{}
\includegraphics[width=\textwidth]{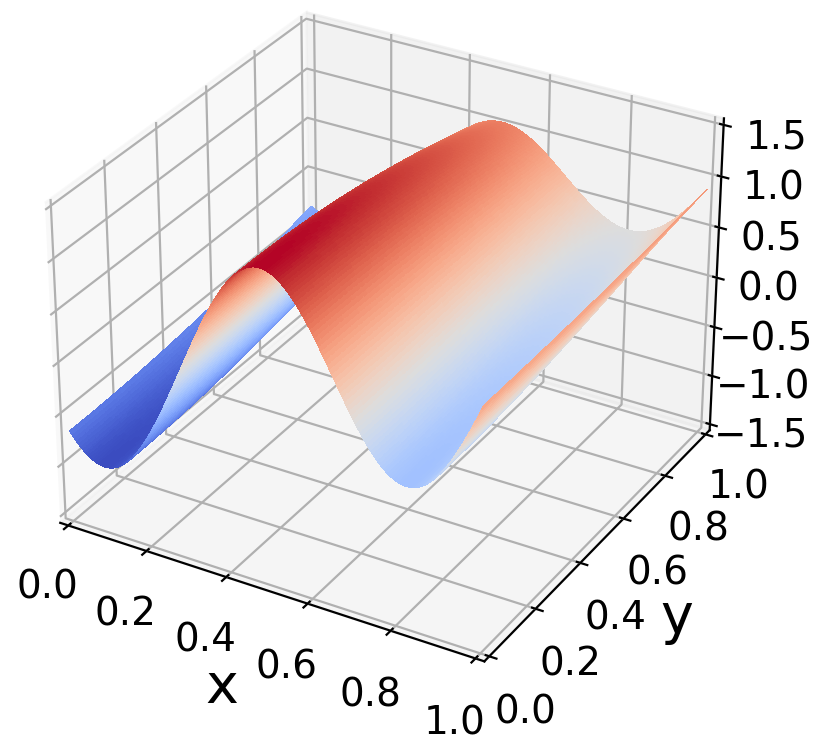}
\end{subfigure}\hfill
\begin{subfigure}{0.3\textwidth}
\centering
\caption{}
\includegraphics[width=\textwidth]{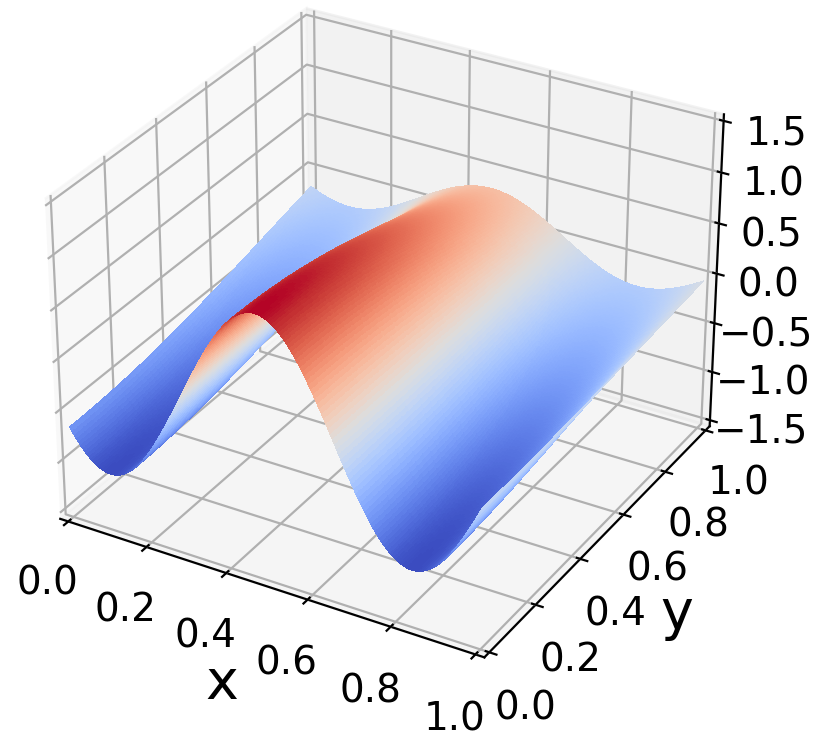}
\end{subfigure}
\caption{(a) Exact low-fidelity function, $z_L$, and (b) exact high-fidelity function, $z_H$, for $ a = 8.5211$. } 
\label{fig:2d_nonlin_exact}
\end{figure}

\begin{figure}[H]
\centering
\includegraphics[width=0.42\textwidth]{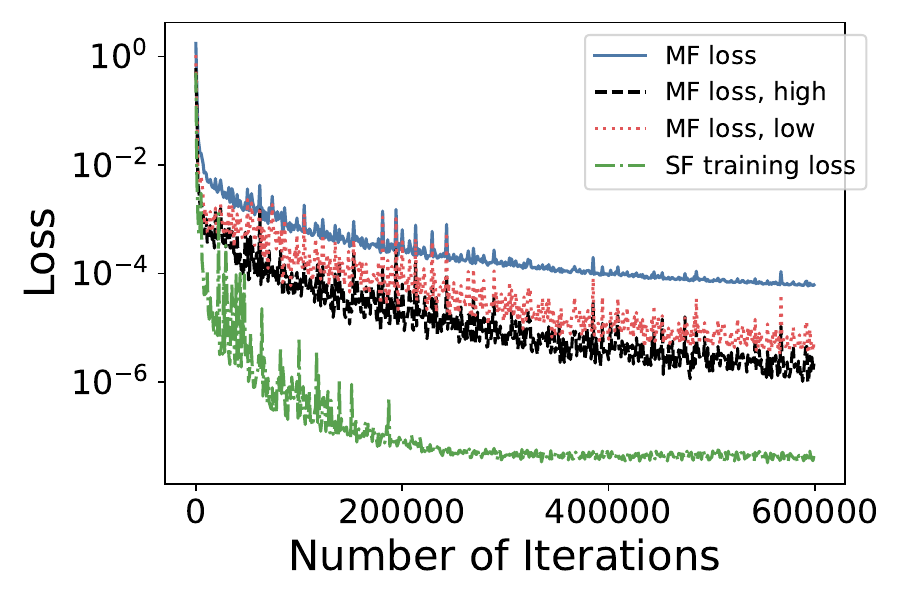}
\caption{Training loss for the two-dimensional problem with nonlinear correlation.} \label{fig:2d_nonlin_loss}
\end{figure}

\subsection{Multiresolution ice-sheet modeling}\label{sec:multiresolution_app}

\begin{figure}[H]
\begin{minipage}{0.55\textwidth}
\begin{table}[H]
    \centering
    \begin{tabular}{c|c}
    \hline
      Method   &  Computational cost (hours) \\ \hline
    Single fidelity, $N_H=20$ & 1.345 \\
    Single fidelity, $N_H=50$ & 1.380 \\
    Multifidelity, $N_H=20$, $N_L=100$ & 3.690 \\
        \hline
    \end{tabular}
    \caption{Computational cost for \textcolor{black}{completing the full training of the  single fidelity and multifidelity DeepONets, not including data generation, for} the multiresolution ice-sheet problem (hours). For the single fidelity training the batch size is $10$ values of $u$ and for the multifidelity training the batch size is $10$ values of $u$ for the low-fidelity data and $5$ values of $u$ for the high-fidelity data. Reported times are on one NVIDIA P100 GPU.}
    \label{tab:Ice-sheets-comp-cost}
\end{table}
\end{minipage}
\begin{minipage}{0.03\textwidth}
\hspace{0.05\textwidth}
\end{minipage}
\begin{minipage}{0.35\textwidth}
\centering
\includegraphics[width=\textwidth]{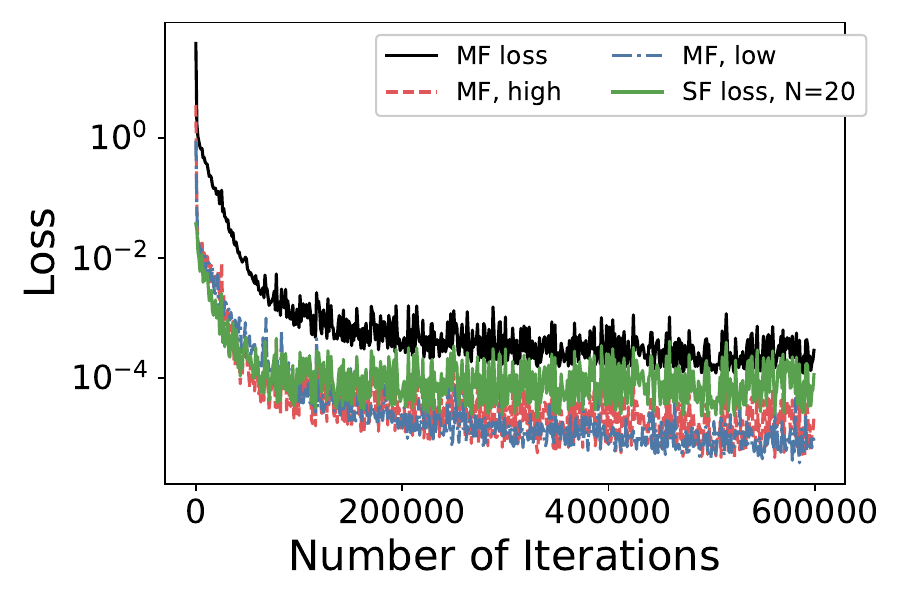}
\caption{Training loss for the multiresolution ice-sheet problem} \label{fig:ice_multires}
\end{minipage}
\end{figure}

\subsection{Multiorder ice-sheet modeling}\label{sec:multiorder_app}

\begin{figure}[H]
\begin{minipage}{0.55\textwidth}
\begin{table}[H]
    \centering
     \begin{tabular}{c|c}
    \hline
      Method   &  Computational cost (hours) \\ \hline
    Single fidelity & 2.68 \\
    Multifidelity &10.32 \\
        \hline
    \end{tabular}
    \caption{Computational cost for \textcolor{black}{completing the full training of the  single fidelity and multifidelity DeepONets, not including data generation, for}  for the multiorder ice-sheet example (hours).  For the single fidelity training the batch size is $20$ values of $u$ and for the multifidelity training the batch size is $20$ values of $u$ for both the low- and high- fidelity data. Reported times are on one NVIDIA P100 GPU.}
    \label{tab:ice_multiord-comp-cost}
\end{table}
\end{minipage}
\begin{minipage}{0.03\textwidth}
\hspace{0.05\textwidth}
\end{minipage}
\begin{minipage}{0.35\textwidth}
\centering
\includegraphics[width=\textwidth]{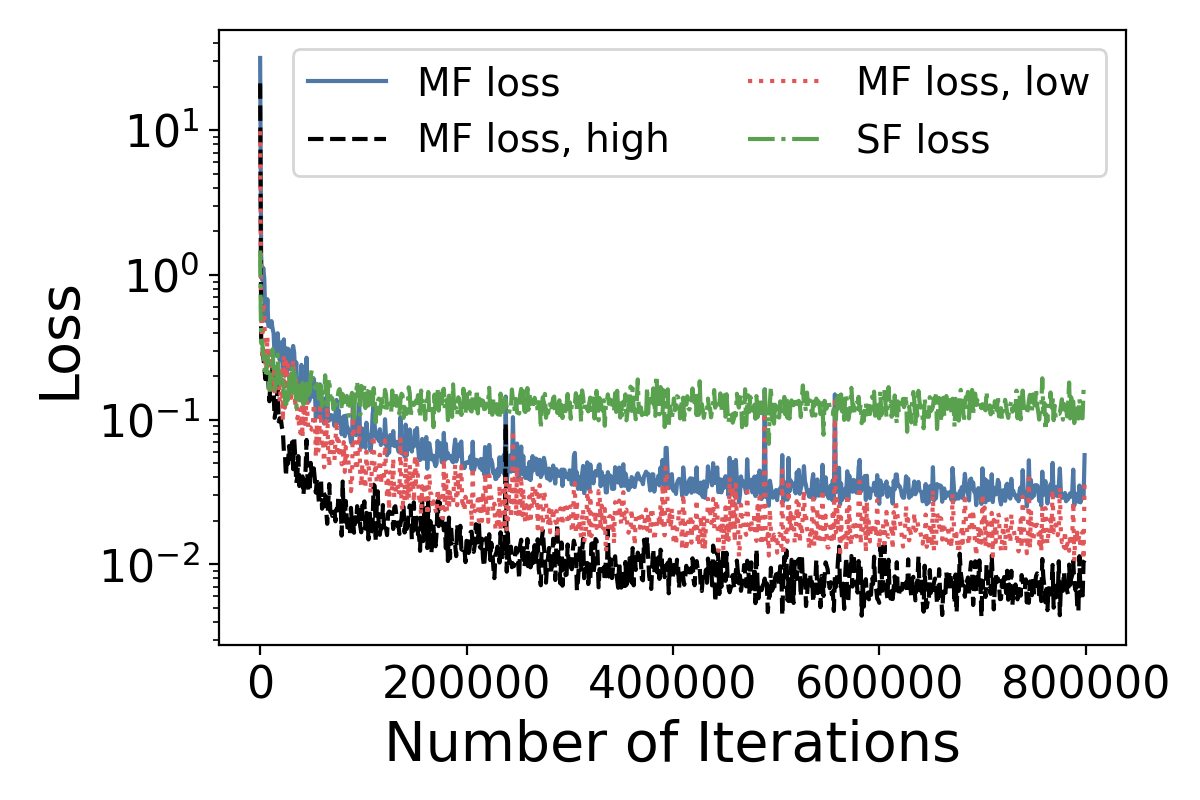}
\caption{Training loss for the multiorder ice-sheet example.} \label{fig:ice_multiord}
\end{minipage}
\end{figure}

\subsection{One-dimensional physics-informed}\label{sec:1d_pi_app}

\begin{figure}[H]
\centering
\includegraphics[width=0.42\textwidth]{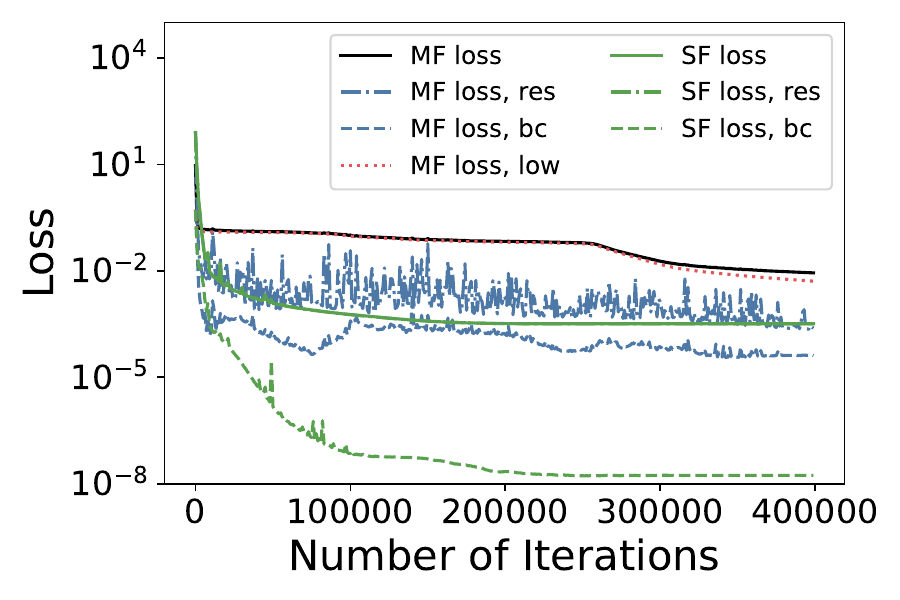}
\caption{Training loss for the one-dimensional physics-informed problem.} \label{fig:1d_phys_loss}
\end{figure}

\subsection{Two-dimensional physics-informed (Viscous Burgers equation)}\label{sec:burgers_app}

\begin{figure}[H]
\centering
\begin{subfigure}{0.3\textwidth}
\caption{}
\includegraphics[width=\textwidth]{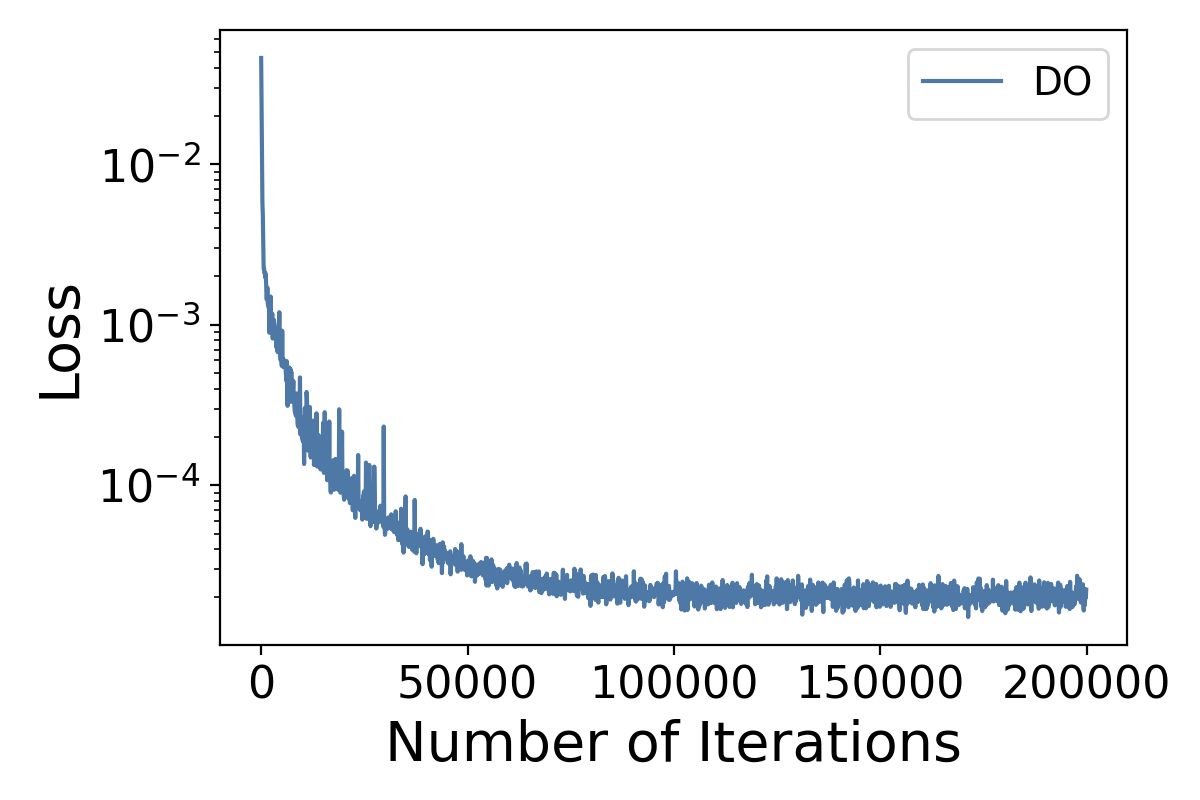}
\end{subfigure}
\begin{subfigure}{0.3\textwidth}
\caption{}
\includegraphics[width=\textwidth]{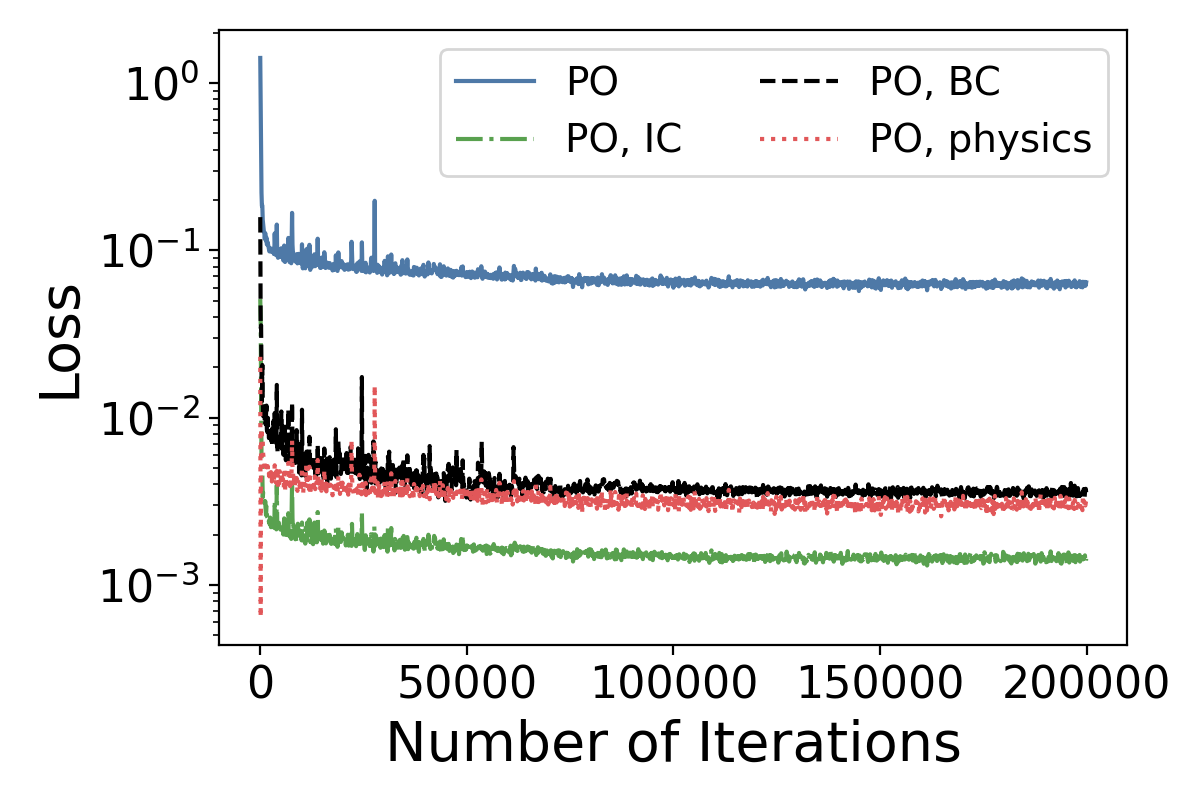}
\end{subfigure}
\begin{subfigure}{0.3\textwidth}
\caption{}
\includegraphics[width=\textwidth]{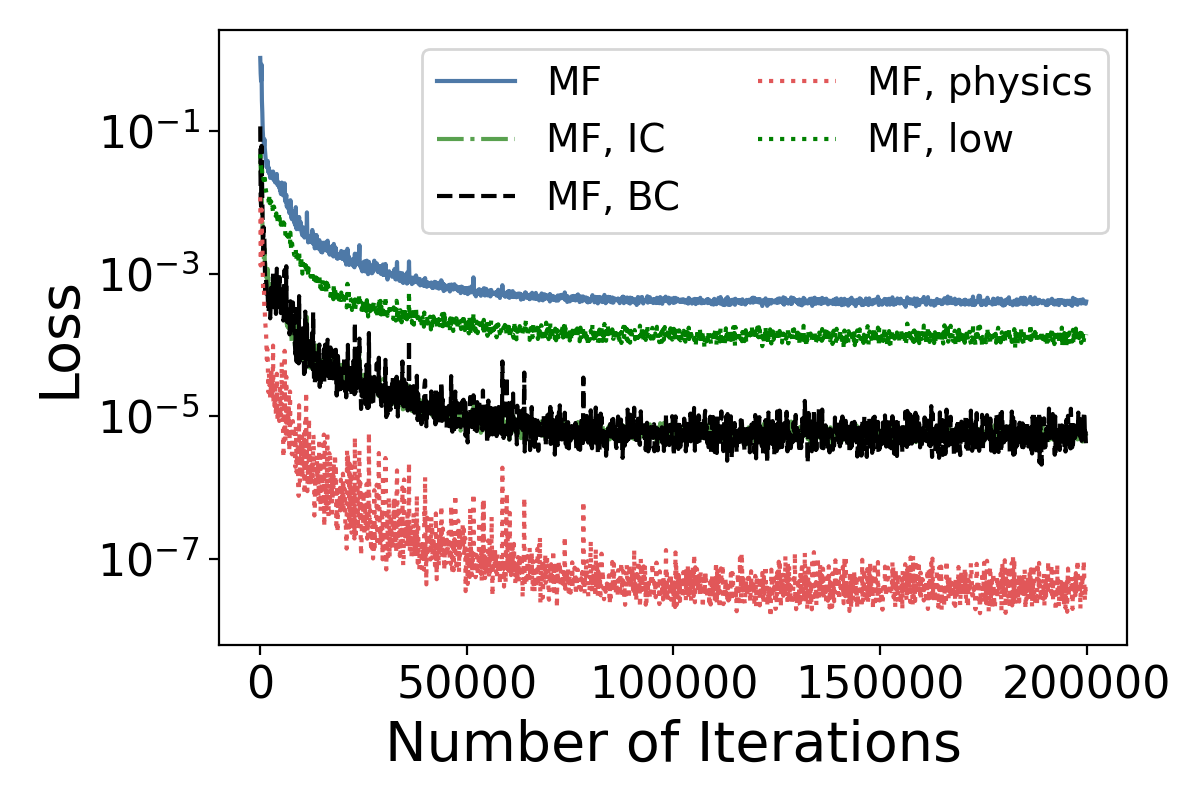}
\end{subfigure}
\caption{Example of training loss for Burgers equation with $\nu = 10^{-4}$, $N_L = 1000$, and non-noisy data. (a) Data-only losses. (b) Physics-only losses. (c) Multifidelity losses. } \label{fig:2d_phys_loss}
\end{figure}

\begin{table}[H]
    \centering
    \begin{tabular}{c|P{25mm}|P{25mm}|P{25mm} | P{25mm}}
    \hline
     Parameters  & $\nu = 10^{-2}$ \newline $N_L = 1000$ & $\nu = 10^{-3}$ \newline $N_L = 1000$ & $\nu = 10^{-4}$ \newline $N_L = 1000$ & $\nu = 10^{-4}$ \newline $N_L = 200$ \\ \hline
    Data-only & 0.995 & 0.995 & 0.994 & 0.994  \\
    Data-only with noise & 0.995 &0.994 &  0.996 & 0.994 \\
    Physics-only & 9.92 & 9.96 &9.88& --- \\
    Multifidelity &  14.26 & 14.46 & 14.21 & 14.21 \\
    Multifidelity with noise &14.55 & 14.19 & 14.56 & 14.25\\ \hline
    \end{tabular}
    \caption{Physics-informed multifidelity: viscous Burgers equation. Computational cost for each case (hours). For the physics-only and multifidelity training, $N_H = 1000.$ The high-fidelity batch size is $25$ values of $u$ and the low-fidelity batch size is $100$ values of $u$. Reported times are on one NVIDIA P100 GPU. \textcolor{black}{Note that the physics-only case does not use any low-fidelity data, so the number of low fidelity samples, $N_L$, does not impact the results. Therefore, we do not report a result for the physics-only case with $N_L = 200$ in the last column, because it would be identical to the physics-only case with $N_L = 1000$.}}
    \label{tab:Burger-training-time}
\end{table}

\section{Non-composite Multifidelity DeepONets}\label{sec:noncomp}

\begin{figure}[ht]
\centering
\includegraphics[width=0.5\textwidth]{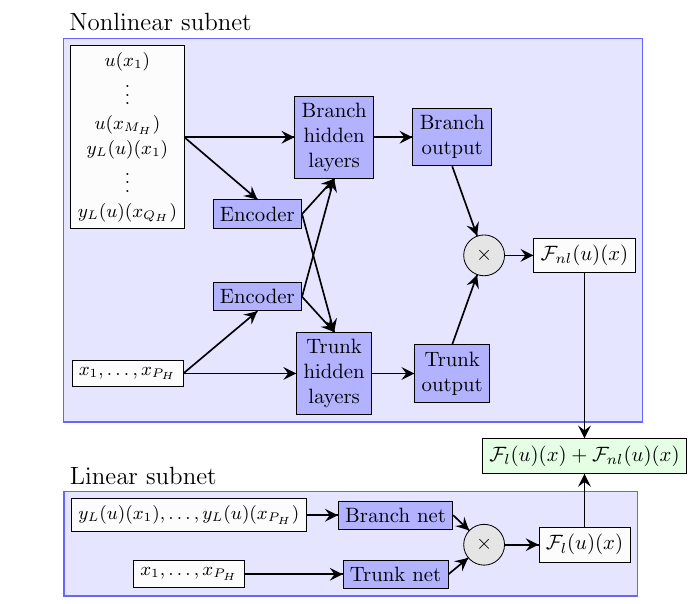}
\caption{Non-composite multifidelity DeepONet setup.} \label{fig:Non-comp-setup}
\end{figure}

In this section we consider the case where we have access to the low-fidelity data at every point in the high-fidelity training set and any location where we wish to output the operator values, as may be true if we know the exact function that generates the low-fidelity data. We call this the ``non-composite'' framework, and refer to the setup in the main text as the ``composite'' framework. The non-composite setup is shown in Fig. \ref{fig:Non-comp-setup}. We assume that we have high-fidelity data as described in Sec. \ref{sec:DD_note}. We also have a low-fidelity operator $y_L(u)(x)$ that we can use to generate values on the high-fidelity data points for input into the linear and nonlinear branch networks. 

We simultaneously train two DeepONets to learn the linear and nonlinear correlations between the low-fidelity operator and the high-fidelity operator. The loss function is given by
\begin{equation}
\mathcal{L}_{non-composite}(\theta) = \lambda_1 \mathcal{L}_{HF}(\theta_{nl}, \theta_{l}) + \lambda_3 \left( \sum w_{nl}^2 + \sum b_{nl}^2 \right) . \label{eq:loss_noncomp}
\end{equation}

\begin{figure}[]
\centering
\begin{subfigure}{0.45\textwidth}
\centering
\includegraphics[width=\textwidth]{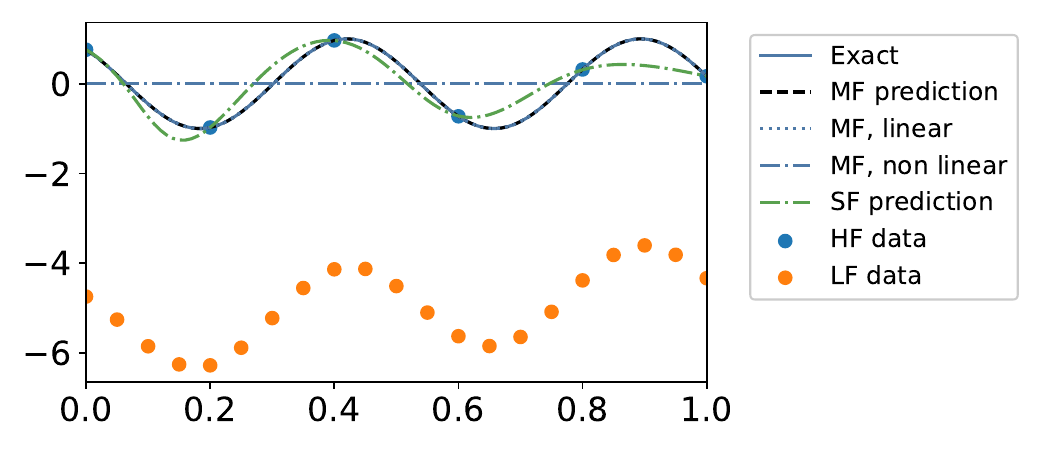}
\caption{Results of the single fidelity and multifidelity predictions of the high-fidelity data for a=13.2579 with the non-composite setup. } \label{fig:1d_noncomp_lin_results}
\end{subfigure}\hfill
\begin{subfigure}{0.45\textwidth}
\centering
\includegraphics[width=\textwidth]{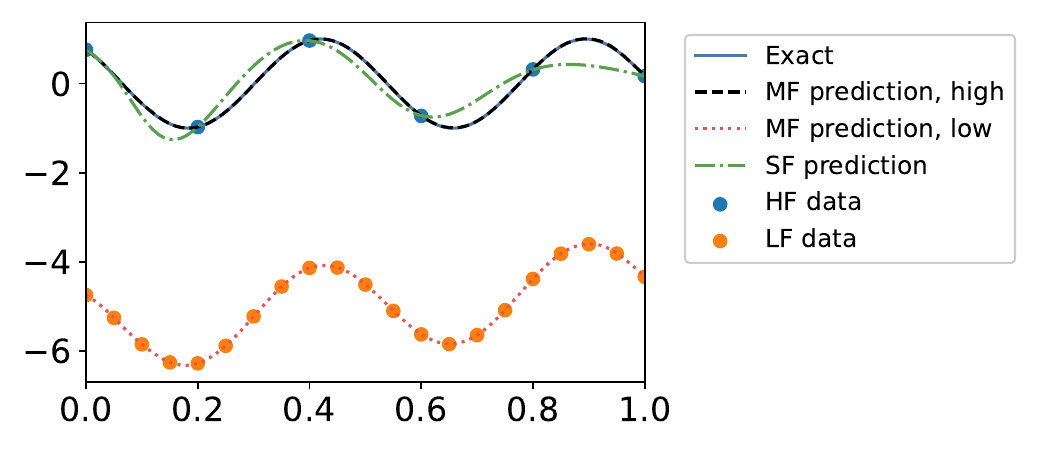}
\caption{Results of the single fidelity and multifidelity predictions of the high-fidelity and low-fidelity data for $a=13.2579$ with the composite setup. } \label{fig:1d_comp_lin_results}
\end{subfigure}
\caption{Comparison between the non-composite (a) and composite (b) setups.}
\label{fig:1d_lin_results}
\end{figure}

\begin{figure}[H]
\begin{minipage}{0.55\textwidth}
\begin{table}[H]
    \centering
        \begin{tabular}{c|c }
    \hline
    Parameter & Value \\
    \hline
      LF data   &  $M_L = P_L= 21$ \\
      HF data  & $M_H = P_H = 6$  \\
      Number of datasets & $N_L = N_H = 5$  values of $a$ \\
      $\lambda_1$ & $1\times 10^{-1}$ \\
      $\lambda_2$ & 1 \\
      $\lambda_3$ & $1\times 10^{-1}$ \\
      $\lambda_4$ & $1\times 10^{-3}$ \\
      SF learning rate & (5e-3, 2000, 0.9)\\
      MF learning rate & (1e-3, 5000, 0.97)\\
      SF network size & 3 layers, 30 neurons \\
      MF low-fidelity network size &  3 layers, 30 neurons\\
      MF linear network size & 1 layer, 5 neurons \\
      MF nonlinear network size & 2 layers, 20 neurons\\
          \hline
    \end{tabular}
    \caption{Parameters for the composite case.}
    \label{tab:1d_comp_lin_ex}
\end{table}
\end{minipage}
\begin{minipage}{0.42\textwidth}
\centering
\includegraphics[width=\textwidth]{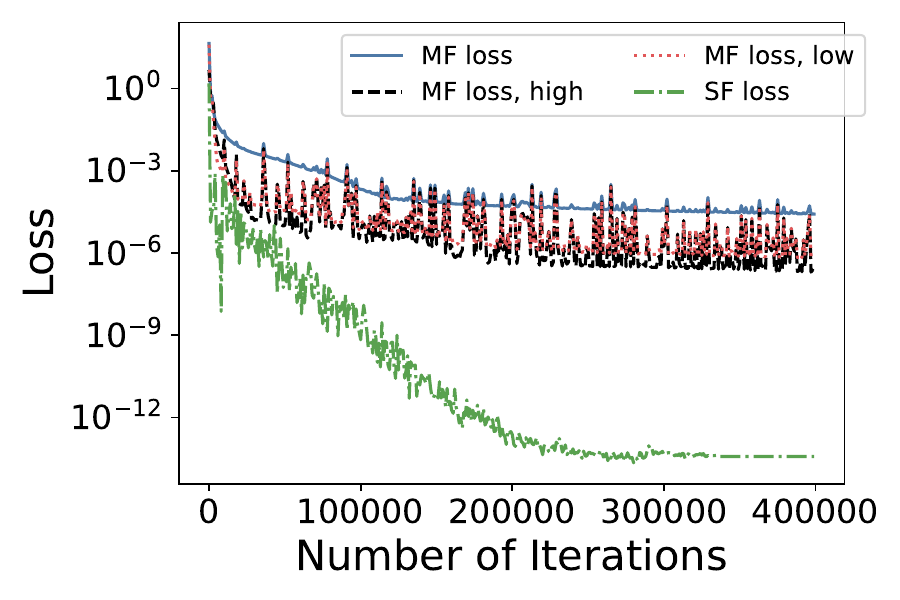}
\caption{Composite loss.} \label{fig:1d_comp_loss}
\end{minipage}
\end{figure}

\begin{figure}[H]
\begin{minipage}{0.55\textwidth}
\begin{table}[H]
    \centering
            \begin{tabular}{c|c }
    \hline
    Parameter & Value \\
    \hline
      LF data   &  $P_L= 21$ \\
      HF data  & $M_H = P_H = 6$  \\
      Number of datasets & $N_H = 5$  values of $a$ \\
      $\lambda_1$ & $1\times 10^{-1}$ \\
      $\lambda_3$ & $1\times 10^{-1}$ \\
      SF learning rate & (5e-3, 2000, 0.9)\\
      MF learning rate & (5e-5, 5000, 0.9)\\
      SF network size & 3 layers, 30 neurons \\
      MF linear network size & 1 layer, 7 neurons \\
      MF nonlinear network size & 2 layers, 20 neurons\\
          \hline
    \end{tabular}
    \caption{Parameters for the non-composite case.}
    \label{tab:1d_noncomp_lin_ex}
\end{table}
\end{minipage}
\begin{minipage}{0.42\textwidth}
\centering
\includegraphics[width=\textwidth]{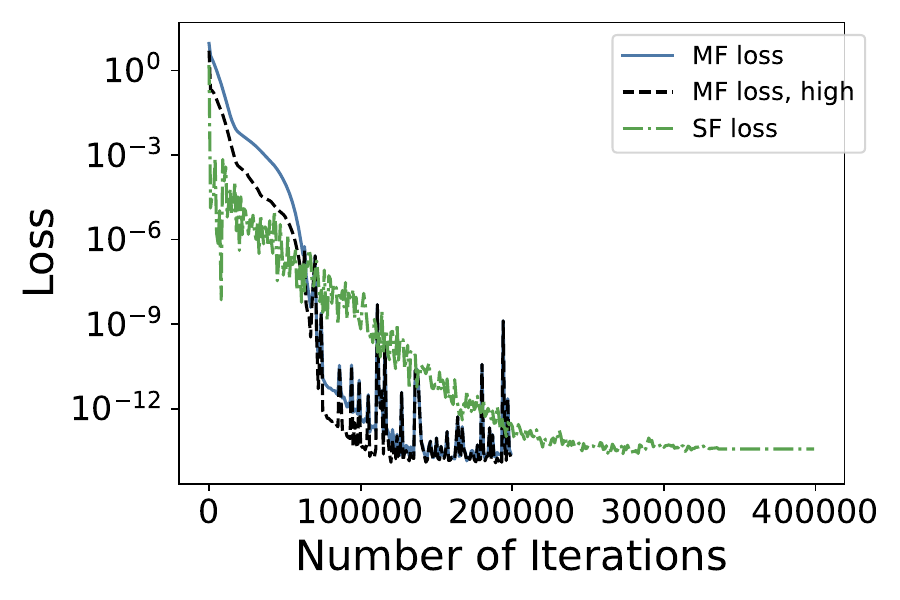}
\caption{Non-composite loss.} \label{fig:1d_noncomp_loss}
\end{minipage}
\end{figure}

We illustrate the impact of including an exact solver for the low-fidelity data instead of a DeepONet by considering the same example problem for both cases. We take 
\begin{align}
    y_L(u)(x) &= \sin(u)+x-5.5\\
    y_H(u)(x) &=  \sin(u) \\
    u &= ax-4
\end{align}
for $x \in [0, 1]$ and $a \in [10, 14]$. Note that $y_H(u)(x) = y_L(u)(x)-x+5.5$.  
The training parameters for the composite training are given in Tab. \ref{tab:1d_comp_lin_ex}, results are given in Fig. \ref{fig:1d_comp_lin_results}. The learned linear correlation is: 
\begin{equation}
    \mathcal{F}_l(u)(x) = 0.9998\mathcal{F}_{LF}(u)(x)-0.9989x+5.4989 - 0.00028 x \mathcal{F}_{LF}(u)(x) \label{eq:comp_lin}. 
\end{equation}
The training parameters for the non-composite training are given in \ref{tab:1d_noncomp_lin_ex}, and the results are given in Fig. \ref{fig:1d_noncomp_lin_results}. For this problem, we can recover the learned linear correlation as: 
\begin{equation}
    \mathcal{F}_l(u)(x) = 1.0000003 y_L(u)(x)-1.000002x+5.500004 - 4.0047\times 10^{-7} x y_L(u)(x)  \label{eq:NONcomp_lin}. 
    \end{equation}
Because the non-composite framework does not introduce errors from the low-fidelity DeepONet output, Eq. \ref{eq:NONcomp_lin} is more accurate than Eq. \ref{eq:comp_lin}, although both agree well with the exact equation. 

\section{Mathematical models for ice-sheets} \label{sec:ice_models}
We consider an ice-sheet with domain $\Sigma \in \mathbb{R}^2$ that is fixed in time (the ice-sheet can change thickness, but its domain does not change.) Denote the spatial coordinates by $(x, y, z) \in \mathbb{R}^3$, where $z=0$ corresponds to sea level. The ice domain at time $t$ is given by:  
$$
\Omega(t) := \{(x,y,z) \; \text{s.t.} \; (x,y) \in \Sigma, \; \text{and}\; l(x,y,t) < z < s(x,y,t) \},
$$
where $\Gamma_l(t)  := \{(x,y,z) \; \text{s.t.} \; z = l(x,y,t)\}
$ denotes the lower surface of the ice at time $t$ and $\Gamma_s(t)  := \{(x,y,z) \; \text{s.t.} \; z = s(x,y,t)\}$ denotes the upper surface of the ice. The ice bed topography is  assumed to be constant in time and is given by $\Gamma_b := \{(x,y,z) \; \text {s.t.} \; z = b(x,y)\}$. Denote the ice velocity by $\mathbf{u} = (u, w, v)$. 

\subsection{Ice-sheet flow model} \label{sec:first_order}

Many methods can be used to model the evolution of the ice-sheet over time. At the highest order, land ice can be modeled as a shear-thinning Stokes flow driven by gravity. 

The Stokes equation gives:
\begin{align}
    -\nabla \cdot \sigma & = \rho \mathbf{g} \\
    \nabla \cdot \mathbf{u} &= 0
\end{align}
with pressure $p$, ice density $\rho$, stress tensor $\sigma = 2\mu \mathbf{D} - pI$, and strain rate tensor  $\mathbf{D}_{ij}(\mathbf{u}) = \frac{1}{2} \left(\frac{\partial u_i}{\partial x_j}+\frac{\partial u_j}{\partial x_i} \right)$.
The nonlinear viscosity is given by 
\begin{equation} \label{viscosity}
    \mu = \frac12 A(T)^{-\frac{1}{n}} \, D_e(\mathbf u)^{\frac{1}{n}-1}
\end{equation}
with $n\geq 1$, typically $n = 3$. $A$ is the ice flow factor that depends on the ice temperature $T$. The \emph{effective strain rate} $D_e(\mathbf u) = \frac{1}{\sqrt{2}}|\mathbf{D}(\mathbf{u})|$, where $| \cdot |$ denotes the Frobenius norm.
Stokes equations are closed by the following boundary conditions:
$$
\left \{ \begin{array}{lll}
\sigma \mathbf{n} = 0 & \text{on } \Gamma_s & \text{stress free, atmospheric pressure neglected} \\
\sigma \mathbf{n} = \rho_w \, g\,  \min(z,0) \bf{n} & \text{on } \Gamma_{m}  & \text{boundary condition at the ice margin} \\
\mathbf u = \mathbf u_d & \text{on } \Gamma_{d}  & \text{Dirichlet condition at internal boundary} \\
\bf{u} \cdot \bf{n} = 0,\; (\sigma \bf{n})_{\parallel} = \beta \bf{u}_{\parallel} & \text{on } \Gamma_g & \text{impenetrability + sliding condition} \\
\sigma \mathbf{n} = \rho_w \, g\,  z\, \bf{n} & \text{on } \Gamma_f & \text{back pressure from ocean under ice shelves} \\
\end{array} \right.
$$

$\beta(x, y)$ is the friction (or sliding) coefficient, which is typically difficult to measure and unknown and $\rho_w$ is the density of the ocean water. The boundary condition at the margin includes the ocean back-pressure term, when the margin is partially submerged ($z<0$). 

The thickness of the ice is given by $H(x,y,t) = s(x,y,t)-l(x,y,t)$ on $\Sigma \times [0, t_f]$ and evolves according to
\begin{equation} \label{thickness}
  \partial_t H + \nabla \cdot (\mathbf{\bar u} H) = f_H 
\end{equation}
where $\displaystyle \mathbf{\bar u} := \frac{1}{H} \int_l^s \mathbf{u}\, dz$ is the depth-averaged velocity and $f_H$ is an accumulation rate, accounting for accumulation and melting at the upper surface and at the base of the ice-sheet. Typically, the ice margin is an outlet ($\mathbf{\bar u} \cdot \mathbf{n} >0$) so no boundary condition are needed. We constrain $H$ to be non-negative. 

Because the Stokes equations are difficult and computationally intensive to solve, a series of simplified models have been derived. These models exploit the shallow nature of ice-sheets to introduce approximations that decrease the computational intensity. We now introduce the Mono-layer higher-order (MOLHO) and Shallow Shelf Approximation (SSA) models. 

\subsection{Mono-layer higher-order (MOLHO)}
MOLHO~\cite{dosSantos2022} model is based on the Blatter-Pattyn~\cite[e.g., ][]{Dukowicz2010} approximation that can be derived neglecting the terms $w_x$ and $w_y$ in the strain-rate tensor $D$ and, using the continuity equation, replacing $w_z$ with $-(u_x+v_y)$:
\begin{equation}
    \mathbf{D} = \begin{bmatrix} u_x  & \frac12(u_y + v_x) & \frac12 u_z \\[2mm]
    \frac12(u_y + v_x) & v_y & \frac12 u_z \\[2mm]
        \frac12 u_z & \frac12 v_z & -(u_x+v_y) \end{bmatrix}.
\end{equation}
This leads to the following elliptic equations in the horizontal velocity $(u,v)$
\begin{equation}
    -\nabla \cdot (2\mu \hat{\mathbf{D}} ) = - \rho g \nabla s
\end{equation}
with 
\begin{equation}
    \hat{\mathbf{D}} = \begin{bmatrix} 2u_x + v_y & \frac12 (u_y + v_x) & \frac12 u_z \\
    \frac12 (u_y + v_x) & u_x + 2v_y & \frac12 v_z \end{bmatrix}
\end{equation}
Here the gradient is two-dimensional: $\nabla = [\partial_x, \partial_y]^T$. The viscosity $\mu$ is given by \eqref{viscosity} with the effective strain rate 
$$
D_e = \sqrt{u_x^2 + v_y^2 + u_x v_y + \frac{1}{4} (u_y + v_x)^2 + \frac{1}{4} u_z^2 + \frac{1}{4} v_z^2}.
$$
The boundary conditions read
$$
\left \{ \begin{array}{lll}
2\mu \hat{\mathbf{D}}\, \mathbf{n} = 0 & \text{on } \Gamma_s & \text{stress free, atmospheric pressure neglected} \\
2\mu \hat{\mathbf{D}}\, \mathbf{n} = \psi \bf{n} & \text{on } \Gamma_{m}  & \text{boundary condition at at ice margin} \\
\mathbf u = \mathbf u_d & \text{on } \Gamma_{d}  & \text{Dirichlet condition at internal boundary} \\
2\mu \hat{\mathbf{D}}\, \mathbf{n} = \beta \bf{u}_{\parallel} & \text{on } \Gamma_g & \text{sliding condition} \\
2\mu \hat{\mathbf{D}}\, \mathbf{n} = 0 & \text{on } \Gamma_f & \text{free slip under ice shelves} \\
\end{array} \right.
$$
where $\psi = \rho g (s-z) \mathbf{n} + \rho_w \, g\,  \min(z,0) \bf{n}$, which can be approximated with its depth-averaged value  $\bar{\psi} = \frac12 g H (\rho - r^2 \rho_w)$, $r$ being the  submerged ratio $r=\max\left(1-\frac{s}{H},0\right)$. 

The MOLHO model consists in solving the weak form of the Blatter-Pattyn model, with the ansatz that the velocity can be expressed as :
$$
\mathbf u(x,y,z) = \mathbf u_b(x,y) + \mathbf u_v(x,y) \left( 1 - \left(\frac{s-z}{H}\right)^{n+1}\right)
$$
The problem is then formulated as a system of two 2d PDEs in $\mathbf u_b$ and $\mathbf u_v$ (for a detailed derivation see~\cite{dosSantos2022}). Note that the depth-averaged velocity is given by
$\mathbf{\bar u} = \mathbf u_b + \frac{(n+1)}{(n+2)} \;\mathbf u_v$.

\subsection{Shallow Shelf Approximation (SSA)}
The shallow shelf approximation~\cite{morland_johnson_1980} is a simplification of the Blatter-Pattyn model, assuming that the velocity is uniform in $z$, so $\mathbf{u}=\mathbf{\bar u}$. It follows that $u_z=0$ and $v_z=0$, giving:

\begin{equation}
    \mathbf{D} = \begin{bmatrix} u_x  & \frac12 (u_y + v_x) & 0 \\
    \frac12 (u_y + v_x) & v_y & 0 \\ 
        0 & 0 & -(u_x+v_y) \end{bmatrix}, \quad 
    \hat{\mathbf{D}} = \begin{bmatrix} 2u_x + v_y & \frac12 (u_y + v_x) & 0 \\
    \frac12 (u_y + v_x) & u_x + 2v_y & 0 \end{bmatrix}, 
\end{equation}
and $D_e = \sqrt{u_x^2 + v_y^2 + u_x v_y + \frac{1}{4} (u_y + v_x)^2}$.
The problem simplifies to a 2D equation in $\Sigma$
$$
- \nabla \cdot \left( 2 \mu H  \hat {\mathbf D }(\mathbf{\bar u})\right) + \beta \mathbf{\bar u} =  -  \rho g H  \nabla s, \quad \text{in } \Sigma
$$
with $\bar \mu = \frac12 \bar{A}(T)^{-\frac{1}{n}} \, D_e(\mathbf{\bar u})^{\frac{1}{n}-1}$, where $\bar{A}$ is the depth-averaged flow factor and with boundary conditions:
$$
\left \{ \begin{array}{lll}
2\mu \hat{\mathbf{D}}(\mathbf{\bar u})\, \mathbf{n} = \bar \psi \bf{n} & \text{on } \Gamma_{m}  & \text{boundary condition at ice margin} \\
\mathbf{\bar u} = \mathbf{\bar u}_d & \text{on } \Gamma_{d}  & \text{Dirichlet condition at internal boundary} \\
\end{array} \right .
$$
Recall that $\bar{\psi} = \frac12 g H (\rho - r^2 \rho_w)$, $r$ being the the submerged ratio $r=\max\left(1-\frac{s}{H},0\right)$. With abuse of notation, here $\Gamma_{m}$ and $\Gamma_d$ are intended to be subsets of $\partial \Sigma$.

\end{document}